\documentclass[%
 aip,
 amsmath,amssymb,
 reprint,%
]{revtex4-2}

\usepackage{graphicx}
\usepackage{dcolumn}
\usepackage{bm}

\usepackage[utf8]{inputenc}
\usepackage[T1]{fontenc}
\usepackage{mathptmx}
\usepackage{etoolbox}
\usepackage{geometry}
\usepackage{color}
\usepackage{float}
\usepackage{mathrsfs}
\usepackage{mathtools}
\usepackage{amsthm}
\usepackage{graphicx}
\usepackage{amssymb}
\usepackage{amsfonts}
\usepackage{amstext}
\usepackage{dsfont}
\usepackage{mathdots}
\usepackage{xcolor}
\usepackage{lipsum}
\usepackage{array}
\usepackage{multirow}
\usepackage{babel}
\usepackage{subfig}
\usepackage[utf8]{inputenc}

\newtheorem{proposition}{Proposition}

\newtheorem*{remark}{Remark}

\renewcommand{\arraystretch}{1.1}

\makeatletter
\def\@email#1#2{%
 \endgroup
 \patchcmd{\titleblock@produce}
  {\frontmatter@RRAPformat}
  {\frontmatter@RRAPformat{\produce@RRAP{*#1\href{mailto:#2}{#2}}}\frontmatter@RRAPformat}
  {}{}
}%
\makeatother

\begin{document}

\preprint{AIP/123-QED}

\title[Notes on Resonant and Synchronized States]{Notes on Resonant and Synchronized States in Complex Networks}
\author{Paolo Bartesaghi}
\email{paolo.bartesaghi@unimib.it}
\affiliation{University of Milano - Bicocca, Via Bicocca degli Arcimboldi 8, 20126 Milano, Italy.}

\date{\today}

\begin{abstract}
Synchronization and resonance on networks are some of the most remarkable collective dynamical phenomena. The network topology, or the nature and distribution of the connections within an ensemble of coupled oscillators, plays a crucial role in shaping the local and global evolution of the two phenomena. This article further explores this relationship within a compact mathematical framework and provides new contributions on certain pivotal issues, including a closed bound for the average synchronization time in arbitrary topologies; new evidences of the effect of the coupling strength on this time; exact closed expressions for the resonance frequencies in terms of the eigenvalues of the Laplacian matrix; a measure of the effectiveness of an \textit{influencer node}'s impact on the network; and, finally, a discussion on the existence of a resonant synchronized state. Some properties of the solution of the linear swing equation are also discussed within the same setting. Numerical experiments conducted on two distinct real networks - a social network and a power grid - illustrate the significance of these results and shed light on intriguing aspects of how these processes can be interpreted within networks of this kind.
\end{abstract}

\maketitle

\begin{quotation}
Why do pink flamingos perform a pre-breeding dance in perfect unison? How does the entire Australian coral reef manage to synchronize to release millions of gametes at the same time each year? How does it happen that a number of individuals in a social network take part with increasing involvement in collective ritual events? These questions concern the problem of how multiple agents within a network can achieve an identical synchronized dynamic behavior, possibly with an amplification of the intensity of the phenomenon over time. Prompted by the ubiquity and transversality of these processes, the paper offers some novel insights into the already immense research field on the dynamics of networks of coupled identical harmonic oscillators.
\end{quotation}

\section{Introduction}
\label{section1}

Up to several thousand great flamingos, during the pre-breeding period, gather in southern France to perform in synchrony a variety of stereotyped movements. Such group performances, involving elaborate courtship rituals and synchronized dances, are believed to play a role in enabling reproductive synchronization\cite{Perrot2016}.
In the Australia's Great Barrier Reef, many coral species tune their reproductive cycles with incredible precision, synchronizing each year all on the same day, just after the full moon in November, to collectively release millions of gametes in one of the largest scale mass reproduction events on the planet \cite{Stone2002}.
In the epidemiological context, synchronization processes govern the dynamics of disease spread through complex networks in anthropic environments. Examples are the dynamics of childhood diseases such as measles, mainly determined by the annual seasonal cycle, or the pulse vaccination-induced synchronization in simple epidemic models. It is a peculiar feature that ecological communities tend to adjust to periodic environmental or climatic driving variables. Although ecologists and epidemiologists have essentially opposite goals, the mathematical laws underlying the different population dynamics are very similar and frequently involve internal synchronization or resonance with external driving factors \cite{Earn1998}.
We find countless examples in different settings. In the biomechanical context, \citet{Brumley2012} study a mechanism of self-propulsion in ancestral eukaryotes, unicellular organisms equipped with tiny mobile flagella that allow them to move in a liquid medium. Correlated changes in the directions of the flagellar beating can be explained by the raise of some kind of large scale coordination and synchronization, even in a weakly coupled system as is the case for the green alga \textit{Volvox}. In the human social setting, the emergence of synchronization processes among individuals, on a physical or mental level, can be induced by collective dynamics and participation in repetitive rituals or dances, as occurs in Zikr, a Sufi mystical dance performed by Chechen Sufis and involving hundreds of men at once.
However, it is well known that synchronization phenomena do not always have positive implications. Strong synchronization may indeed be related to pathological activities as is the case for the epileptic seizures and Parkinson's disease. \citet{Hammond2007} have observed an abnormally synchronized oscillatory activity at multiple levels of the basal ganglia-cortical loop, in humans with Parkinson's disease. This excessive synchronization correlates with motor deficit, and its suppression by specific therapies could suggest an effective approach to improve the living conditions of patients with severe motility disorders.

These examples show that synchronization emerges, in general, from the collaboration, competition and interaction of multiple agents when they undergo periodic and cyclic dynamics. A first interpretive model, albeit oversimplified, reduces such systems to a set of coupled identical harmonic oscillators all synchronizing on the same frequency.

The two pioneering papers on synchronization in coupled systems \cite{Yamada1983} and synchronization in chaotic systems \cite{Pecora1990} have stimulated a great deal of interest in the study of complete synchronization of coupled nonlinear dynamical processes. The first paradigmatic model of diffusively coupled identical oscillators on complex network is given in \citet{Pecora2000}, where the master stability function is introduced for the first time. Such a function accomplishes two goals: decrease the computational load and provide a guidance in designing coupling configurations that conform to the required stability. This approach is currently the most widely applied one to determine the stability of the fully synchronized state in terms of the eigenstructure of the connectivity matrix.

The interest in the study of synchronization in coupled systems grew considerably after the discovery of the small-world and scale-free properties of many natural and artificial complex networks. \citet{Barahona2002} analyze the implications of the small-world phenomenon on the synchronization of oscillator networks of arbitrary topology. Subsequent analysis provided evidence that well-connected networks are easier to globally synchronize than regular networks. For instance, \citet{Chen2012} found that synchronization processes typically start first from a small fraction of hub nodes and then spread to nodes with smaller degrees. This is replicated at the cluster level, where a typical synchronization process starts from partial synchronization within clusters to evolve to global complete synchronization.

\citet{Ren2008} analyzes coupled second-order linear harmonic oscillators, described as point masses, under rather mild network connectivity assumptions, showing that their positions and velocities can be synchronized in the presence or absence of a leader.
A similar approach is used in \citet{Zhan2013} where the vibration spectra of a classical mass spring model on complex networks are investigated; and, in \citet{Shang2012}, where synchronization of a leader-follower system of coupled harmonic oscillators connected by dampers is studied in presence of random noises and time delays with an interaction topology modeled by a weighted directed graph.

In general, the presence of dissipative couplings has the effect to equate the state of the nodes so that the system enters a positively invariant set, the absorbing domain, in finite time. It is well-known that if $\bf G$ is a constant matrix whose eigenvalues have only negative real parts, then $e^{t{\bf G}}$ acts as a uniform contraction; that is, $||e^{t{\bf G}}||\leq Ke^{-\eta t}$ for some positive real constants, $K$ and $\eta$. This fact guarantees that the solution of the equation is uniformly asymptotically stable and the synchronization manifold remains invariant under the flow of the equation that governs the system.

For numerical estimates and model prediction purposes, \citet{Pietras2019} stress the importance of reduction of complex oscillatory systems. A number of surveys are devoted to this topic \cite{Arenas2008}. \citet{Arenas2021} study the synchronized state in a population of network-coupled heterogeneous oscillators, showing that the steady-state solution of the linearized dynamics may be written as a geometric series whose terms represent the different spatial scales of the network. They prove that this expansion converges for arbitrary frequency distributions and for both undirected and directed networks, provided that the adjacency matrix is primitive.
Also \citet{Eroglu2017} propose a survey on the main theory behind complete, generalized and phase synchronization phenomena in complex networks with an approach close to the one followed in the present paper.

Investigations of the transient and steady-state dynamics of synchronous machines have significantly impacted the description of power grid networks. \citet{Sorrentino2021} study the swing equation and its linearization on transmission networks, highlighting the role of symmetries in reducing the computational complexity of modal decomposition\cite{SanchezGarcia2020}. \citet{Sorrentino2022} extend this approach to multilayer networks.

A crucial issue is the time and speed of synchronization. \citet{Grabow2011} focus on the speed of convergence as a collective time scale for synchronizing systems on different topologies, ranging from completely ordered and grid-like, to completely disordered and random topologies. They find, for instance, that, at fixed in-degree, topological randomness induces shorter synchronization times, whereas intermediate randomness in the small-world regime induces longer synchronization times.

Many other contributions have been given on specific issues, such as in \citet{Skardal2014} and \citet{Skardal2019}, where authors introduce a synchrony alignment function that encodes the interplay between the network structure and oscillator frequencies as an objective measure for the synchronization properties of a network of heterogeneous oscillators with natural frequencies or in \citet{Su2009}, where authors revisit the synchronization problem for coupled oscillators in a dynamic proximity network. \citet{Dorfler2013} prove a closed-form condition for synchronization that can be stated in terms of the Moore-Penrose pseudoinverse of the Laplacian matrix. \citet{Liu2016}, in the wider context of control theory, address the problem of how networks organize themselves to balance control with functionality. In particular, they review the process of local pinning synchronization, when controllers, applied to a subset of nodes, help synchronize the network. Also in \citet{Lu2014}, pinning synchronization in the new setting of networks of networks is investigated. \citet{Slotine2004} study the role, in both nature and system design, of co-existing power leaders, to which the networks synchronize, and knowledge leaders, to whose parameters the networks adapt.

One of the most important contributions to the study of synchronization in coupled systems is undoubtedly the Kuramoto model\cite{Dorfler2014}. In this model the underlying topology is fully connected, even if it was later adapted to the scenario of nearest neighbor interactions.
Despite the simplicity and effectiveness of the Kuramoto model and all its variants, several works have highlighted some of its shortcomings. For example, \citet{Shahal2020} address the issue of the importance of synchronization in human networks for our civilization. They study, in an accurate experimental setting, the synchronization between violin players, arranged in a network framework with full control in terms of network connectivity, coupling strength and delay. Their results show that players can adjust their playing period and cancel connections by ignoring frustrating signals, to find a stable solution. Their observations reveal that the usual models for coupled networks, such as Kuramoto model, cannot always be applied to human networks, as additional degrees of freedom enable new strategies and produce better solutions than are possible within such a model. Corroborating these observations, it is argued that, during social interactions, participants are continuously active, each modifying their actions in response to the changing actions of their partners. This continuous mutual adaptation results in an interactional synchrony to which all members contribute. \citet{Dumas2010} study the brain activity of the participants in this process, and found that states of interactional synchrony correlate with the emergence of an interbrain synchronization network between right centroparietal regions. These regions have been found to play a central role in social interaction within an interindividual brainweb.

Another limiting aspect of classical models is the presence of strictly local interactions.
For this reason, \citet{Estrada2017} study how long-range interactions can affect synchronization dynamics, proposing a model of coupled oscillators based on the d-path Laplacian matrices. Their analysis reveals that long-range interactions improve network synchronizability with an impact that depends on the original structure.

The dynamics of opinion formation within a social network is also a widely studied field. In this kind of network, each node can be assigned, for instance, a discrete or continuous score according to his/her agreement or disagreement on a given topic. By introducing a coupling between members in order to describe the dragging effect in their opinion formation process, each member in the network can change his/her opinion according to a synchronization mechanism or according to some kind of resonant phenomenon with some source/leader status.

Related to this issue, \citet{Jenkins2013} presents a very interesting analysis of the so-called self-oscillation phenomenon, which is distinct from the more familiar phenomenon of forced resonance, whereas \citet{Fiore2021} focus on networks with influencers, a group of hubs that couple strongly and connect different parts of the network. They show that, subjecting influencers to an optimal noise, can result in enhanced network synchronization. This effect, called coherence resonance, emerges from a synergy between network structure and stochasticity and is highly nonlinear, vanishing as the noise is too weak or too strong.

In the light of the previous literature review, this paper aims to analyze, in a sense from first principles, some aspects of the synchronization and resonance phenomena, and, in general, of the transient dynamics in networks of oscillators with a heterogeneous topological structure. A unified approach is maintained in the description of all phenomena. The main new contributions and results of the present work can be summarized as follows:

\begin{enumerate}
\item a synthetic framework in which to inscribe the description of synchronization and resonance phenomena on networks of coupled identical harmonic oscillators with arbitrary topology;
\item a closed expression of the average synchronization time that encodes information about the topological structure of the network;
\item new evidences that, for a weakly coupled system, complete and dense networks achieve the state of perfect synchronization more rapidly than, for instance, the path network and that the complete network shows a synchronization time which is shorter for weak coupling than for strong coupling;
\item closed expressions for the resonance frequencies of the network as a function of the eigenvalues of the Laplacian matrix;
\item a detailed account of how a single node in the network can serve as a leader (influencer) or a follower (receiver) depending on its topological location, when a resonance phenomenon arises inside the network; 
\item a proof that the presence of dissipative terms pulls down the possible resonance frequencies to a single one, corresponding to the synchronization frequency.
\item a new perspective on the linearization of the swing equation and its solution in the proposed setting;
\item some new mathematical properties of the polar decomposition of the symplectic matrix governing the synchronization phenomenon.
\end{enumerate}

The paper is organized as follows. Section \ref{section2} introduces the basic notations and the general model. Section \ref{section3} has a didactic purpose: it shortly illustrates the case of networked coupled identical linear oscillators without damping. In Section \ref{section4}, existing results are framed into the broader context described, and new ones about the mathematical properties of the involved operators and the rate of synchronization are presented. Section \ref{section5} is devoted to resonance phenomena and, after some recalls on the case of the single driven oscillator, provides new results related to resonant frequencies on networks when an external driving force is placed at a specific node. Section \ref{section6} deals with a special case of the general equation introduced, specifically the linear swing equation and its solution. Finally, section \ref{section7} applies the previous findings to two real world networks, the Zachary club social network and the Syrian power grid, by illustrating a number of interesting remarks. All the proofs of the propositions in this article are given in the Appendixes \ref{Appendix A} and \ref{Appendix B}.

\section{Preliminaries}
\label{section2}
Our purpose is basically to describe the behavior of a system of $n$ coupled harmonic oscillators within a network structure. To model the local interactions between the $n$ harmonic oscillators, we use an undirected graph $G=(V,E)$ without loops, where $V$ is the set of vertices or nodes and $E\subset V\times V$ is the set of edges or links. Throughout the paper we will assume that the network is connected i.e., that a path exists between each pair of distinct vertices.
The actual configuration is completely encoded in the adjacency matrix ${\bf A}=\{a_{ij}\},\ i,j=1,\dots, n$, where $a_{ij}=1$ when $(i,j)\in E$, that is when nodes $i$ and $j$ show an actual local interaction, and $0$ otherwise. We will assume $a_{ij}=a_{ji}$, which is a perfectly symmetric interaction within a pair of nodes. The degree of a node $i$ is defined as ${k_i}=\sum_{j=1}^{n}a_{ij}=({\bf A}^{2})_{ii}$ and ${\bf K}={\rm diag}\, \{ k_{i}\}$ is the diagonal matrix of the degrees of all the nodes. We will denote by $\lambda_{i}$ and $\psi_{i},\ i=1,\dots, n$, respectively the eigenvalues and eigenvectors of the adjacency matrix $\bf A$, where the eigenvalues are listed in decreasing order: $\lambda_{1}\geq \dots \geq \lambda_{n}$. The Laplacian matrix is then defined as ${\bf L}={\bf K}-{\bf A}$. We will denote by $\mu_{i}$ and $\phi_{i},\ i=1,\dots, n$, respectively the eigenvalues and eigenvectors of the Laplacian matrix $\bf L$, where again the eigenvalues are listed in decreasing order: $\mu_1> \mu_2 > \dots>\mu_n=0$. In particular, $\phi_n=\frac{1}{\sqrt n}[1,1,\dots ,1]^{T}=\frac{1}{\sqrt n}{\bf u}^{T}$, being ${\bf u}\in {\mathbb R}^{n}$ the vector of all $1$'s. Similarly, we will denote by ${\bf 1}={\bf u}{\bf u}^{T}$ the $n-$square matrix of all $1$'s and by $\bf I$ the identical matrix.

\begin{figure}[h]
	\includegraphics[width=0.50\linewidth]{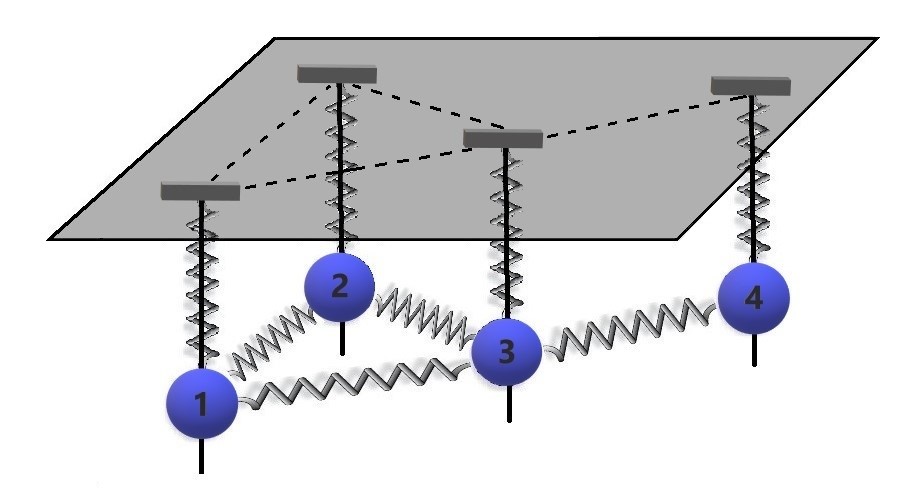}
	\centering
	\caption{A simple example of a network of coupled harmonic oscillators.}
	\label{fig1}      
\end{figure}

To have a concrete picture of the system of coupled oscillators, let us look at the network in Fig. \ref{fig1}. This example will be used throughout the paper as a toy model. We have $n=4$ objects of mass $m_i,\ i=1, \dots , 4$, attached to a rigid support by identical springs. These objects are also connected to each other according to a network scheme with adjacency matrix ${\bf A}$ by other identical springs.
All vertical springs have the same elastic constant $c_1$, and all the springs connecting different objects have the same elastic constant $c_2$. Each mass is allowed to oscillate only in the vertical direction around the equilibrium position. Let ${\bf x}=[x_1,\dots, x_n]^T$ be the vector of the displacements of each mass from its equilibrium position, so that $x_i\in {\mathbb R}$ represents the one-dimensional (vertical, in the figure) displacement of the $i$th object. Let $F_i$ be the force acting on mass $i$ due to the vertical spring and $F_{ij}$ be the force acting on mass $i$ due to its interaction with mass $j$. Then $F_i=-c_{1}x_{i}$ and $F_{ij}=-c_{2}(x_{i}-x_{j})$ and the total force on node $i$ is given by
\begin{equation}
\begin{split}
F_{i}^{\rm tot}=&-c_{1}x_{i}-c_{2}\sum_{j}a_{ij}(x_{i}-x_{j})\\=&-(c_{1}+c_{2}k_{i})x_{i}+c_{2}\sum_{j}a_{ij}x_{j}.
\end{split}
\end{equation}

The forces acting on all the nodes in the network are then collected in the force vector
\begin{equation}
{\bf F}^{\rm tot}_{\rm elastic}=-c_{1}{\bf x}-c_{2}({\bf K}-{\bf A}){\bf x}=-\left[ c_{1}{\bf I}+c_2{\bf L}\right]{\bf x}.
\label{elasticforce}
\end{equation}

Let us add now the damping forces, which are proportional to the velocities ${\bf \dot{x}}={d{\bf x}(t)}/{dt}$ with coefficients $c^{\prime}_{1}$ for the vertical damping and $c^{\prime}_{2}$ for the coupling between nodes. As above, in a similar way, they can be collected into the vector:
\begin{equation}
{\bf F}^{\rm tot}_{\rm damping}=-\left[ c^{\prime}_{1}{\bf I}+c^{\prime}_2{\bf L}\right]{\bf \dot{x}}.
\label{dampingforce}
\end{equation}

The most general equation of motion of the network of coupled and damped harmonic oscillators with masses ${\bf M}={\rm diag}\, \{ m_{i}\}$ is then given by
\begin{equation}
{\bf M}{\bf \ddot{x}}= -\left[ c_{1}{\bf I}+c_2{\bf L}\right]{\bf x}-\left[ c^{\prime}_{1}{\bf I}+c^{\prime}_2{\bf L}\right]{\bf \dot{x}}
\label{equationofmotion}
\end{equation}

which is equivalent to the system of $2n$ differential equations given by
\begin{equation}
\left\{ 
\begin{array}{l}
	{\bf \dot{x}}={\bf v} \\
	{\bf M}{\bf \dot{v}}= -\left[ c_{1}{\bf I}+c_2{\bf L}\right]{\bf x}-\left[ c^{\prime}_{1}{\bf I}+c^{\prime}_2{\bf L}\right]{\bf v} \\
\end{array}
\right.
\label{system1}
\end{equation}

with initial conditions ${\bf x}(0)={\bf x}_{0}$ and ${\bf v}(0)={\bf v}_{0}$. Let us define the vector ${\bf y}\in {\mathbb R}^{2n}$ as
\begin{equation}
{\bf y}\coloneqq
\left( \begin{array}{c}
{\bf x} \\
{\bf v} \\
\end{array} \right)
\label{definitionofy}
\end{equation}

so that, in matrix form,
\begin{equation}
\begin{split}
&\dot{\bf y}=
\left(\begin{array}{c}
{\bf \dot{x}} \\
{\bf \dot{v}} \\
\end{array} \right)=\\&
\left[ \begin{array}{cc} 
{\bf 0} & {\bf I} \\
 -{\bf M}^{-1}\left( c_{1}{\bf I}+c_2{\bf L}\right)  & -{\bf M}^{-1}\left( c^{\prime}_{1}{\bf I}+c^{\prime}_2{\bf L}\right) \\
\end{array} \right]
\left(\begin{array}{c}
{\bf x} \\
{\bf v} \\
\end{array} \right)
\label{system2}
\end{split}
\end{equation}

We will focus on the case of coupled \textit{identical} harmonic oscillators, so that ${\bf M}={\bf I}$. Hence, the system (\ref{system2}) becomes
\begin{equation}
\left\{ 
\begin{array}{l}
{\bf \dot{y}}={\bf G}{\bf y} \\
{\bf y}(0)={\bf y}_{0} \\
\end{array}
\right.
\label{system3}
\end{equation}

where the matrix ${\bf G}$ is
\begin{equation}
{\bf G}\coloneqq
\left[ \begin{array}{cc} 
{\bf 0} & {\bf I} \\
-\left( c_{1}{\bf I}+c_2{\bf L}\right)  & -\left( c^{\prime}_{1}{\bf I}+c^{\prime}_2{\bf L}\right) \\
\end{array} \right]
\label{matrixG}
\end{equation}
We stress that coefficients $c_{2}$ and $c^{\prime}_{2}$ determine the strength of the coupling between the oscillators and that in the following the coupling is meant to be weak when $c^{(\prime)}_{2}/c^{(\prime)}_{1}<1$.

Since we are interested in resonant phenomena, we add to the system an external periodic force applied to one node. Let us imagine that such a driving force $f(t)$, $f: [0, +\infty)\to {\mathbb R}$, acts on node $h,\ h=1,\dots, n$; node $h$ will be called the \textit{leader}. The total force on node $i$ is now
\begin{equation}
F_{i}^{\rm tot}=-c_{1}x_{i}-c_{2}\sum_{j}a_{ij}(x_{i}-x_{j})+f(t)\delta_{hi}
\end{equation}
where $\delta_{hi}=1$ only if $i=h$, $0$ otherwise. The forces acting on the nodes in the network are then described by the vector
\begin{equation}
\begin{split}
&{\bf F}^{\rm tot}=\\&-c_{1}{\bf x}-c_{2}({\bf K}-{\bf A}){\bf x}-c_{1}^{\prime}{\bf v}-c_{2}^{\prime}({\bf K}-{\bf A}){\bf v}+f(t){\bf e}_{h}
\label{totalforce}
\end{split}
\end{equation}

where ${\bf e}_{h}=[0,0,\dots,1,\dots, 0]^{T}$, with $1$ in position $h$, is in the standard basis of ${\mathbb R}^n$. System (\ref{system3}) is then generalized by
\begin{equation}
\left\{ 
\begin{array}{l}
{\bf \dot{y}}={\bf G}{\bf y}+{\bf b}(t) \\
{\bf y}(0)={\bf y}_{0} \\
\end{array}
\right.
\label{system4}
\end{equation}

where 
\begin{equation}
{\bf b}(t)=f(t)\left[ \begin{array}{c} 
{\bf 0} \\
{\bf e}_{h} \\
\end{array} \right]
\label{additiveterm}
\end{equation}

Finally, in addition to the cases described above, we will study a further variant of the system (\ref{system4}) in order to show how this formalism can accommodate the linearization of an important nonlinear problem. In particular, by setting $c_{1}=0$, $c_{2}=1$, $c^{\prime}_{1}=1$, $c^{\prime}_{2}=0$, and by giving the factor ${\bf b}(t)$ in Eq. (\ref{system4}) the form ${\bf b}(t)=[{\bf 0}, {\bf p} ]^{T}$, where ${\bf p}$ is a constant vector in ${\mathbb R}^{n}$, we will gain the \textit{linear swing equation} describing the dynamic response of a network of synchronous machines to small disturbances.

Shortly, we can summarize in the following Table \ref{summary} the cases we are interested in, where $c_{1}$, $c_{2}$, $c^{\prime}_{1}$, $c^{\prime}_{2}$ are the coupling terms in the matrix $\bf G$ in Eq. (\ref{matrixG}) and the last column specifies the form of the driving term ${\bf b}(t)$ in Eq. (\ref{system4}).
\begin{table}[H]
	\begin{center}
		\begin{tabular}{||r|c|c|c|c|c||}
			\hline
{\rm \bf Network System} & $\bf c_{1}$ & $\bf c_{2}$ & $\bf c^{\prime}_{1}$ & $\bf c^{\prime}_{2}$ & ${\bf b}(t)$ \tabularnewline \hline
Coupled  & 1 & 1 & 0 & 0 & ${\bf 0}$ \tabularnewline \hline
Coupled and damped & 1 & 0 & 0 & 1 & ${\bf 0}$ \tabularnewline \hline
Coupled and forced & 1 & 1 & 0 & 0 & $f(t)[{\bf 0}, {\bf e}_{h} ]^{T}$  \tabularnewline \hline
Coupled, damped and forced & 1 & 0 & 0 & 1 & $f(t)[{\bf 0}, {\bf e}_{h} ]^{T}$  \tabularnewline \hline
Linear swing equation & 0 & 1 & 1 & 0 & $[{\bf 0}, {\bf p} ]^{T}$ \tabularnewline \hline
\end{tabular}
\end{center}
\caption{The five cases discussed in the paper: coupled harmonic oscillators; coupled and damped harmonic oscillators; coupled and forced harmonic oscillators; coupled, damped and forced harmonic oscillators; the linear swing equation.}
\label{summary}
\end{table}
The first case in Section \ref{section3} is for illustrative purposes only and it is used to properly introduce notations. The second case in Section \ref{section4} will lead to synchronization phenomena on network. The third and fourth cases in Section \ref{section5}  are the core of the paper and lead to the resonance phenomena on networks. Finally, the fifth case in Section \ref{section6} will allow us to discuss some properties of the linear swing equation and its solution.

\section{Network of coupled identical harmonic oscillators}
\label{section3}
The analysis of this first case is worthwhile since it allows the exact solution of system (\ref{system3}) to be written in a closed form useful for later discussion. Therefore, let us go back to system (\ref{system3}) and set $c^{\prime}_{1}= 0$ and $c^{\prime}_{2}= 0$:
\begin{equation}
\left(\begin{array}{c}
{\bf \dot{x}} \\
{\bf \dot{v}} \\
\end{array} \right)=
\left[ \begin{array}{cc} 
{\bf 0} & {\bf I} \\
-\left( c_{1}{\bf I}+c_{2}{\bf L}\right)   & {\bf 0} \\
\end{array} \right]
\left(\begin{array}{c}
{\bf x} \\
{\bf v} \\
\end{array} \right)
\label{system5}
\end{equation}

Let us denote ${\bf H}=c_{1}{\bf I}+c_{2}{\bf L}$ and $\bf B$ the square root matrix of ${\bf H}$: ${\bf B}^2={\bf H}$.
Since, for any value of the positive coefficients $c_{1}$ and $c_{2}$, the matrix ${\bf H}$ is positive definite, $\bf B$ can be easily obtained as ${\bf B}={\bm \Phi} {\bm \Lambda}^{1/2}{\bm \Phi}^T$, where $\bm \Phi$ is the matrix whose columns are the eigenvectors of ${\bf H}$ and ${\bf L}$ and $\bm \Lambda$ is the diagonal matrix of the positive eigenvalues of ${\bf H}$. 

It is well known that the general solution of the equation  ${\bf \ddot{x}}+{\bf B}^{2}{\bf x}={\bf 0}$, equivalent to system (\ref{system5}), is given by
\begin{equation}
{\bf x}(t)=e^{i{\bf B}t}{\bf x}_{0}^{(1)}+ e^{-i{\bf B}t}{\bf x}_{0}^{(2)},
\label{generalsolution1}
\end{equation}

where ${\bf x}_{0}^{(1)}$ and ${\bf x}_{0}^{(2)}$ are complex conjugate constant vectors satisfying the following conditions:
\begin{equation}
\left\{ 
\begin{array}{l}
{\bf x}(0)=:{\bf x}_{0}={\bf x}_{0}^{(1)}+{\bf x}_{0}^{(2)} \\
{\bf \dot{x}}(0)=:{\bf v}_{0}= i{\bf B}\left( {\bf x}_{0}^{(1)}-{\bf x}_{0}^{(2)} \right). \\
\end{array}
\right.
\label{initialconditions1}
\end{equation}

Therefore
\begin{equation}
\left\{ 
\begin{array}{l}
{\bf x}_{0}^{(1)}=\frac{1}{2}\left[ {\bf x}_{0} -i{\bf B}^{-1}{\bf v}_{0}\right] \\
{\bf x}_{0}^{(2)}=\frac{1}{2}\left[ {\bf x}_{0} +i{\bf B}^{-1}{\bf v}_{0}\right] \\
\end{array}
\right.
\label{initialconditions3}
\end{equation}

and the general solution (\ref{generalsolution1}) takes the usual form
\begin{equation}
\begin{split}
{\bf x}(t)
&=\frac{e^{i{\bf B}t}+e^{-i{\bf B}t}}{2}\, {\bf x}_{0}-\frac{e^{i{\bf B}t}-e^{-i{\bf B}t}}{2}\, i{\bf B}^{-1}{\bf v}_{0}\\
&=\cos({\bf B}t)\, {\bf x}_{0}+ \sin({\bf B}t)\, {\bf B}^{-1}{\bf v}_{0}.\\
\end{split}
\label{generalsolution2}
\end{equation}

Let us notice that solution (\ref{generalsolution2}) of the system (\ref{system3}) can be deduced equivalently, by observing that 
\begin{equation}
{\bf y}(t)=e^{{\bf G}t}{\bf y}_{0}
\label{generalsolution}
\end{equation}

where now matrix ${\bf G}$ in Eq. (\ref{matrixG}) is
$
{\bf G}=
\left[ \begin{array}{cc} 
{\bf 0} & {\bf I} \\
- {\bf H}   & {\bf 0} \\
\end{array} \right].
$

A straightforward computation shows that, in this case,
\begin{equation}
\begin{split}
e^{{\bf G}t}=&
\left[ \begin{array}{cc} 
{\bf I} & {\bf 0} \\
{\bf 0} & {\bf I} \\
\end{array} \right]+
\left[ \begin{array}{cc} 
{\bf 0} & {\bf I} \\
-{\bf H} & {\bf 0} \\
\end{array} \right] t+
\frac{1}{2!}
\left[ \begin{array}{cc} 
-{\bf H} & {\bf 0} \\
{\bf 0} & -{\bf H} \\
\end{array} \right] t^2\\&+
\frac{1}{3!}
\left[ \begin{array}{cc} 
{\bf 0} & -{\bf H} \\
{\bf H}^2 & {\bf 0} \\
\end{array} \right] t^3+
\frac{1}{4!}
\left[ \begin{array}{cc} 
{\bf H}^2 & {\bf 0} \\
{\bf 0} & {\bf H}^2 \\
\end{array} \right] t^4+\dots\\
=&\left[ \begin{array}{cc} 
\cos{\bf B}t & {\bf B}^{-1}\sin{\bf B}t \\
-{\bf B}\sin{\bf B}t & \cos{\bf B}t \\
\end{array} \right]\\
\end{split}
\label{generalsolution3}
\end{equation}

where, as above, ${\bf B}^2={\bf H}$. Hence
\begin{equation}
{\bf y}(t)=\left[ \begin{array}{cc} 
\cos{\bf B}t & {\bf B}^{-1}\sin{\bf B}t \\
-{\bf B}\sin{\bf B}t & \cos{\bf B}t \\
\end{array} \right]
\left(\begin{array}{c}
{\bf x}_{0} \\
{\bf v}_{0} \\
\end{array} \right)
\end{equation}

which is equivalent to
\begin{equation}
\left\{ 
\begin{array}{l}
{\bf x}(t)=\cos({\bf B}t)\, {\bf x}_{0}+{\bf B}^{-1}\sin({\bf B}t){\bf v}_{0} \\
{\bf v}(t)=-{\bf B}\sin({\bf B}t)\, {\bf x}_{0}+ \cos({\bf B}t){\bf v}_{0}\\
\end{array}
\right.
\label{generalsolution4}
\end{equation}

\begin{figure}[H]
	\includegraphics[width=0.52\linewidth]{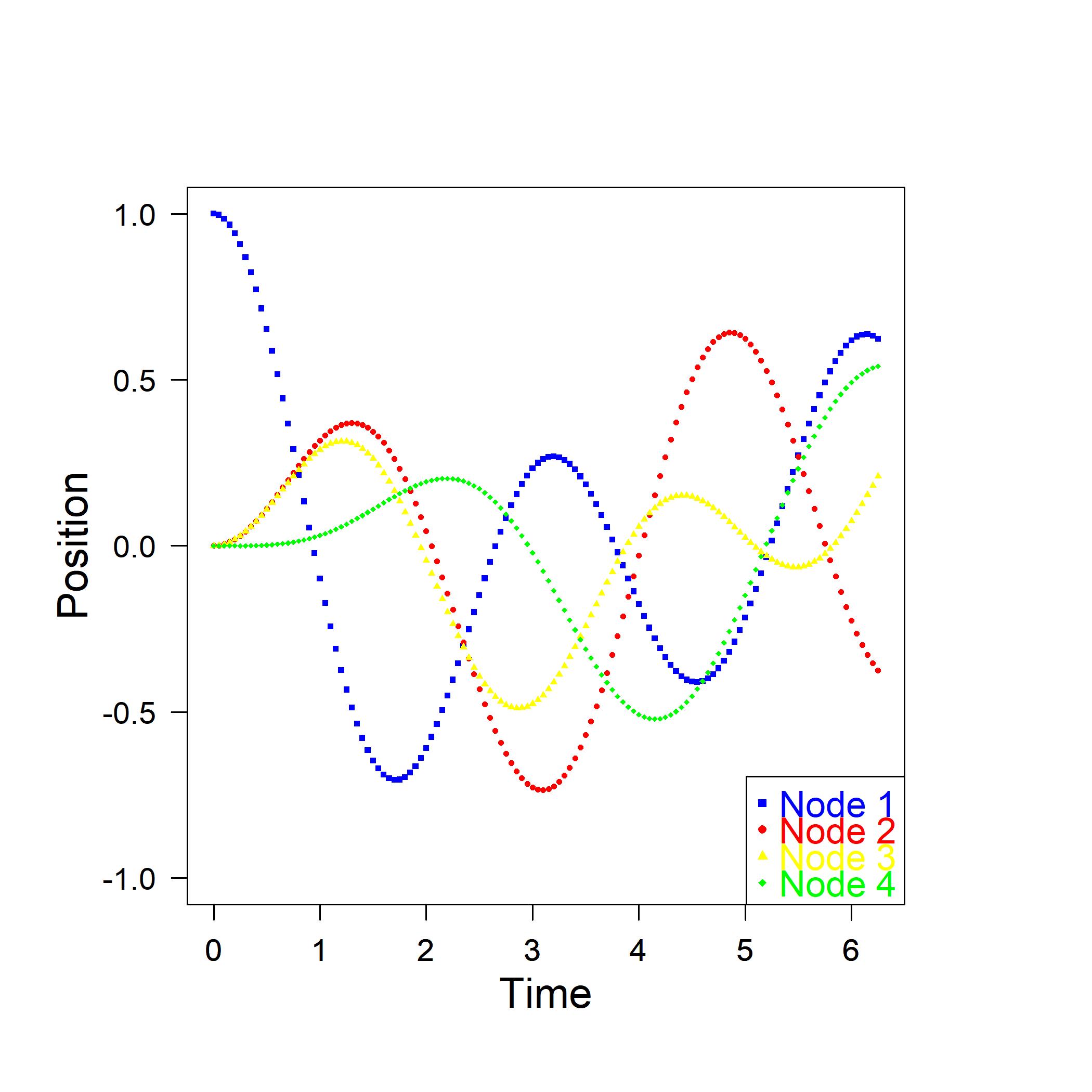}
	\centering
	\caption{Motion of the harmonic oscillators in the example network. Blue square dots refer to node 1, red circle dots to node 2, yellow triangular dots to node 3 and green kite dots to node 4.}
	\label{fig2}      
\end{figure}

Figure \ref{fig2} illustrates the harmonic motion of the four nodes of the network in figure \ref{fig1} given by solution $x_i(t)$ in Eq. (\ref{generalsolution4}). The values for the two coefficients are $c_{1}=c_{2}=1$ and the initial conditions are ${\bf x}_{0}=[1,0,0,0]^T$ and ${\bf v}_{0}=[0,0,0,0]^T$. On the horizontal axis, time ranges in $[0,2\pi]$.

Let us remark that, if we set $c_{1}=c^{\prime}_{1}=c^{\prime}_{2}=0$ and $c_{2}=1$ in system (\ref{system3}), we have a network of free harmonic oscillators connected by identical springs. 

In this case, ${\bf H}={\bf L}$, and the system is described by ${\bf \ddot{x}}+{\bf L}{\bf x}={\bf 0}$.
Since $\det {\bf B}=\det {\bf L}=0$, ${\bf B}^{-1}$ does not exist, and if ${\rm rank}({\bf B})={\rm rank}({\bf B}|{\bf v}_{0})=n-1$, system (\ref{initialconditions1}) admits $\infty^{1}$ solutions. In fact, if we set ${\bf x}_{0}^{(1)}={\bf a}-i{\bf b}$ and ${\bf x}_{0}^{(2)}={\bf a}+i{\bf b}$, system (\ref{initialconditions1}) is equivalent to
\begin{equation}
\left\{ 
\begin{array}{l}
{\bf x}_{0}=2{\bf a} \\
{\bf v}_{0}= 2{\bf B}{\bf b} \\
\end{array}
\right.
\label{initialconditions2}
\end{equation}

Therefore, there are $\infty^1$ values for the imaginary part $\bf b$ of the initial coefficients. Let us notice that ${\bf B}^{-1}$ doesn't exist because of the null eigenvalue $\mu_n=0$ of the Laplacian $\bf L$ and then of the matrix $\bf B$. This eigenvalue is associated with the global translation of all the nodes with the same velocity, i.e. to the stationary state of the whole network, which, in this case, is free and not bounded to any fixed support. Let us also observe that, since $\sum_{j}B_{ij}=\sum_{i}B_{ij}=0$, i.e. the sum of the elements along a row or a column in the matrix $\bf B$ is null, the initial velocity vector ${\bf v}_{0}$ has to satisfy $\sum_{i}v_{0i}=2\sum_{i}\sum_{j}B_{ij}b_{j}=2\sum_{j}\sum_{i}B_{ij}b_{j}=0$. Then, for the system (\ref{initialconditions1}) to be consistent, the initial velocity vector must satisfy the condition $\sum_{i}v_{0i}=0$, equivalent to the classical linear momentum conservation.

\section{Synchronization of coupled identical harmonic oscillators}
\label{section4}

\subsection{Mathematical properties}

Let us set in system (\ref{system3}): $c_{1}=c^{\prime}_{2}=1$ and $c_{2}=c^{\prime}_{1}=0$. This choice for the values of the parameters corresponds to a system of harmonic oscillators attached to a fixed support, free to oscillate around their rest point with constant $c_{1}=1$ but coupled in a damping network with coupling coefficient $c^{\prime}_{2}=1$. Basically it represents a system of $n$ networked point-mass agents connected and interacting only by virtual dampers and can be viewed as a general distributed algorithm for the synchronization of their positions and velocities. Matrix ${\bf G}$ in Eq. (\ref{matrixG}) now becomes
\begin{equation}
{\bf G}=
\left[ \begin{array}{cc} 
{\bf 0} & {\bf I} \\
- {\bf I}   & -{\bf L} \\
\end{array} \right]
\label{matrixG1}
\end{equation}

Let us first observe that $\bf G$ is a symplectic matrix, that is ${\bf G}^{T}{\bf J}{\bf G}={\bf J}$, where $\bf J$ is the $2n-$square antisymmetric matrix (symplectic identity)
\begin{equation}
{\bf J}=
\left[ \begin{array}{cc} 
{\bf 0} & {\bf I} \\
- {\bf I}   & {\bf 0} \\
\end{array} \right]
\end{equation}

As a consequence, we have that if $\lambda$ is an eigenvalue of $\bf G$, so are its complex conjugate $\lambda^{\star}$, $\frac{1}{\lambda}$ and, hence, $\frac{1}{\lambda^{\star}}$. Moreover, if $\lambda$ has multiplicity $m$, then $\frac{1}{\lambda}$ has multiplicity $m$. Finally, both ${\bf G}$ and ${\bf G}^{-1}=-{\bf J}{\bf G}^{T}{\bf J}$ have the same set of eigenvalues. In particular, denoting by $\lambda^{\rm G \pm}_i,\ i= 1,\dots ,n$, the eigenvalues of $\bf G$ and $\psi^{\rm G \pm}_i$ the corresponding eigenvectors, \citet{Ren2008} proves the following proposition

\begin{proposition}
\label{proposition1}
The eigenvalues of $\bf G$ are the $n$ couples of reciprocal values given by
\begin{equation}
\lambda^{\rm G\pm}_i=\frac{1}{2}\left[-\mu_i\pm \sqrt{\mu_i^2-4}\right]
\label{eigenvalues_G}
\end{equation}
and the corresponding normalized eigenvectors are given by
\begin{equation}
\psi^{\rm G \pm}_i=
\frac{1}{\sqrt{1+|{\lambda^{\rm G\pm}_i}|^{2}}}
\left[ \begin{array}{c} 
\phi_i \\
\lambda^{\rm G\pm}_i\phi_i \\
\end{array} \right]
\label{eigenvectors_G}
\end{equation}
\end{proposition}

where $\mu_i$ and $\phi_i$ are the eigenvalues and the corresponding eigenvectors of the Laplacian matrix $\bf L$. In Appendix \ref{Appendix B}, we propose a straightforward alternative proof for the undirected case in line with the one proposed by \citet{Ren2008} for the directed case. Here, we focus on the following remarks.

\begin{remark} 
Proposition \ref{proposition1} has the following straightforward consequence: $\bf G$ is a Hurwitz matrix, i.e. all its eigenvalues are negative or have a negative real part: ${\rm Re}({\lambda^{\rm G \pm}_i})<0, \ \forall i=1,\dots, n$. Indeed, we may have the following four possibilities as in Table \ref{table2}:
\begin{table}[H]
\begin{center}
\begin{tabular}{||l|ll||}
\hline
{\rm if} $\mu_i>2$ & ${\lambda^{\rm G}_i}\in {\mathbb{R}}$, & ${\lambda^{\rm G\pm }_i}<0$ \tabularnewline
\hline
{\rm if} $\mu_i=2$ & ${\lambda^{\rm G}_i}\in {\mathbb R}$, & ${\lambda^{\rm G\pm }_i}=-1$ \tabularnewline
\hline
{\rm if} $0<\mu_i<2$	& ${\lambda^{\rm G}_i}\in {\mathbb C}$, & $-1<{\rm Re}({\lambda^{\rm G\pm }_i})=-\frac{\mu_i}{2}<0$  \tabularnewline
\hline
{\rm for} $\mu_n=0$	& ${\lambda^{\rm G}_n}\in {\mathbb I}$, & ${\lambda^{\rm G \pm}_n}=\pm i$  \tabularnewline
\hline
\end{tabular}
\end{center}
\caption{Eigenvalues of the matrix $\bf G$}
\label{table2}
\end{table}
Let us also notice that, since $\mu_1 >k_{\rm max}\geq 2$ for any non trivial network, there is always a real negative eigenvalue. Moreover, as expected, the trace of the matrix $\bf G$ is equal to
\begin{equation*}
{\rm tr}\, {\bf G}=\sum_{j=1}^{n}\left( \lambda^{\rm G +}_j +\lambda^{\rm G -}_j \right)=-\sum_{j=1}^{n}\mu_{j}=-\sum_{j=1}^{n}k_j=-k_{\rm tot}
\end{equation*}
where $k_{\rm tot}$ is the total degree of the network.
\end{remark}

A general solution to system (\ref{system3}) is again of the usual form
\begin{equation}
{\bf y}(t)=e^{{\bf G}t}{\bf y}_{0}
\label{generalsolution5}
\end{equation}
with matrix ${\bf G}$ as in Eq. (\ref{matrixG1}).

\begin{remark}
Let us immediately notice that, if we apply the general time evolution in Eq. (\ref{generalsolution5}) to an initial state equal to an eigenvector $\psi^{\rm G \pm}_i$ as in Eq. (\ref{eigenvectors_G}), such that $\mu_i>0$, we get that, for $t\to +\infty$:
\begin{equation*}
{\bf y}(t)=e^{{\bf G}t}\psi^{\rm G \pm}_i=e^{{\lambda^{\rm G\pm}_i}t}\psi^{\rm G \pm}_i\to 0
\end{equation*}
since ${\rm Re}({\lambda^{\rm G\pm }_i})<0$. We can then identify the eigenstates of the matrix $\bf G$ as the initial states such that the synchronization leads to the extinction of the oscillation for all the nodes of the network.
\end{remark}

In Appendix \ref{Appendix A}, we propose some theoretical results about the polar decomposition of the matrix $\bf G$ in Eq. (\ref{matrixG1}), which have relevance in the synchronization process. In particular, we will show that the asymptotic behavior is entirely encoded in the unitary component $\bf U$ of the matrix ${\bf G}$. Since $\bf G$ is a real non-singular matrix, its polar decomposition can be directly obtained as ${\bf G}={\bf G}({\bf G}^{T}{\bf G})^{-1/2}({\bf G}^{T}{\bf G})^{1/2}$. Therefore, we set
\begin{equation}
{\bf G}={\bf U}{\bf P}
\label{polardecomposition_text}
\end{equation}
where
\begin{equation}
\label{definitions_UP_text}
{\bf U}={\bf G}({\bf G}^{T}{\bf G})^{-1/2} \qquad {\rm and}\qquad  {\bf P}=({\bf G}^{T}{\bf G})^{1/2}
\end{equation}

The next proposition states that, in terms of asymptotic behaviors, the evolution operators $e^{{\bf G}t}$ and $e^{{\bf U}t}$ are equivalent.

\begin{proposition}
\label{proposition2}
For any initial condition ${\bf y}(0)={\bf y}_{0}$, both $e^{{\bf G}t}{\bf y}_{0}$ and $e^{{\bf U}t}{\bf y}_{0}$ show the same asymptotic behavior ${\bf {\tilde{y}}}(t)$, for $t \to +\infty$, where
\begin{equation}
\begin{split}
\label{asymptoticsolution}
{\bf {\tilde{y}}}(t) &\coloneqq
\frac{1}{n}
\left[ \begin{array}{cc} 
\cos t\, {\bf 1} &
\sin t\, {\bf 1} \\
-\sin t\, {\bf 1} &
\cos t\, {\bf 1} \\
\end{array} \right]
\left(\begin{array}{c}
{\bf x}_{0} \\
{\bf v}_{0} \\
\end{array} \right)\\
&=
\frac{1}{n}
\left( \begin{array}{c} 
\cos t\, {\bf 1}{\bf x}_{0}+ \sin t\, {\bf 1}{\bf v}_{0} \\
-\sin t\, {\bf 1}{\bf x}_{0}+ \cos t\, {\bf 1}{\bf v}_{0} \\
\end{array} \right)
\end{split}
\end{equation}
and $\bf 1$ is the $n$-square all $1$'s matrix.
\end{proposition}

\subsection{Example}

\begin{figure}[H]
\centering
	\subfloat[]{\includegraphics[width=0.30\textwidth]{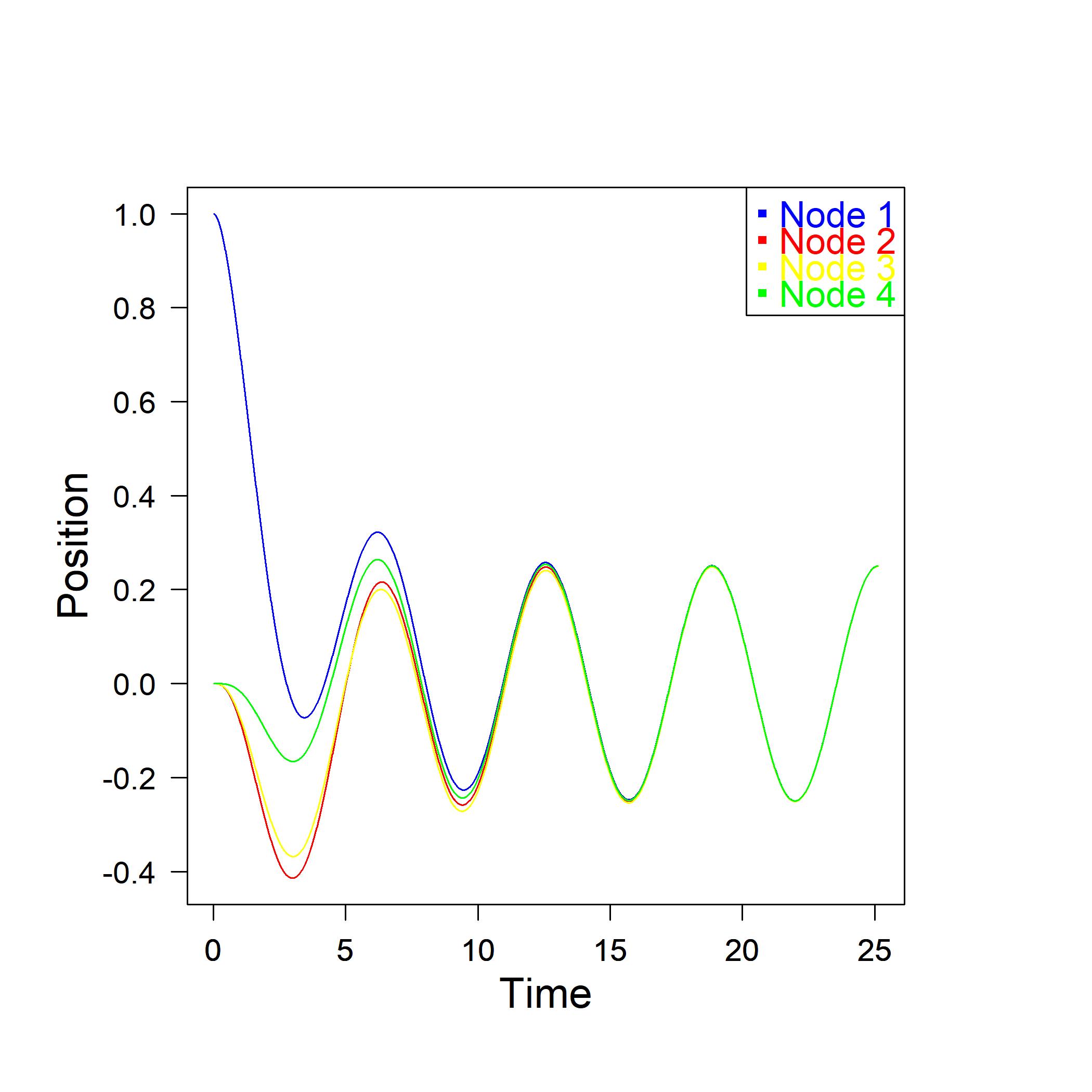}}
	\\
	\subfloat[]{\includegraphics[width=0.30\textwidth]{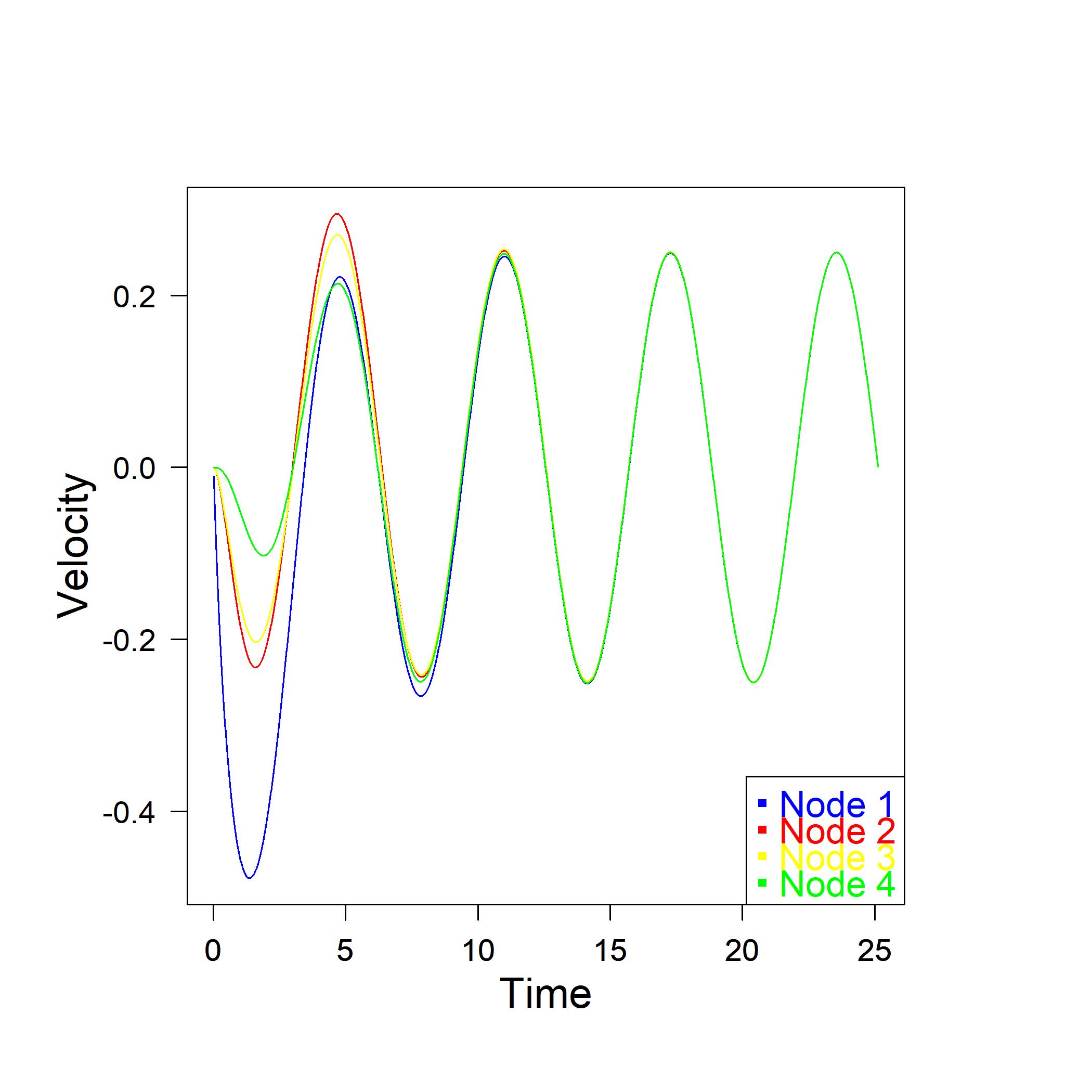}}
	\\
	\subfloat[]{\includegraphics[width=0.30\textwidth]{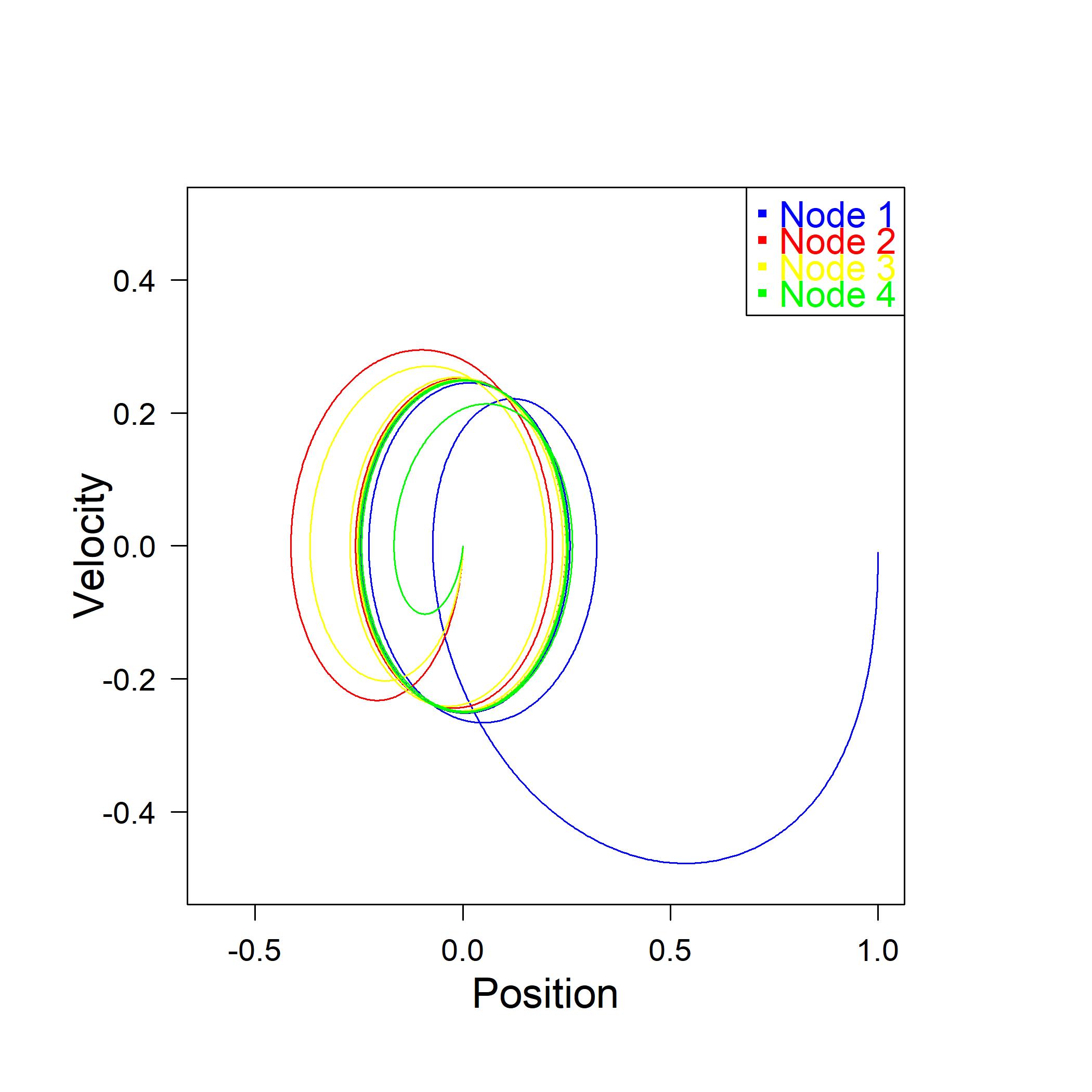}}
\caption{Synchronization of the harmonic oscillators in the example network with initial conditions ${\bf x}_{0}=[1,0,0,0]^T$ and ${\bf v}_{0}=[0,0,0,0]^T$: panel (a) position, panel (b) velocity, and panel (c) phase portrait.}
\label{fig3} 
\end{figure}

\begin{figure}[H]
\centering
	\subfloat[]{\includegraphics[width=0.30\textwidth]{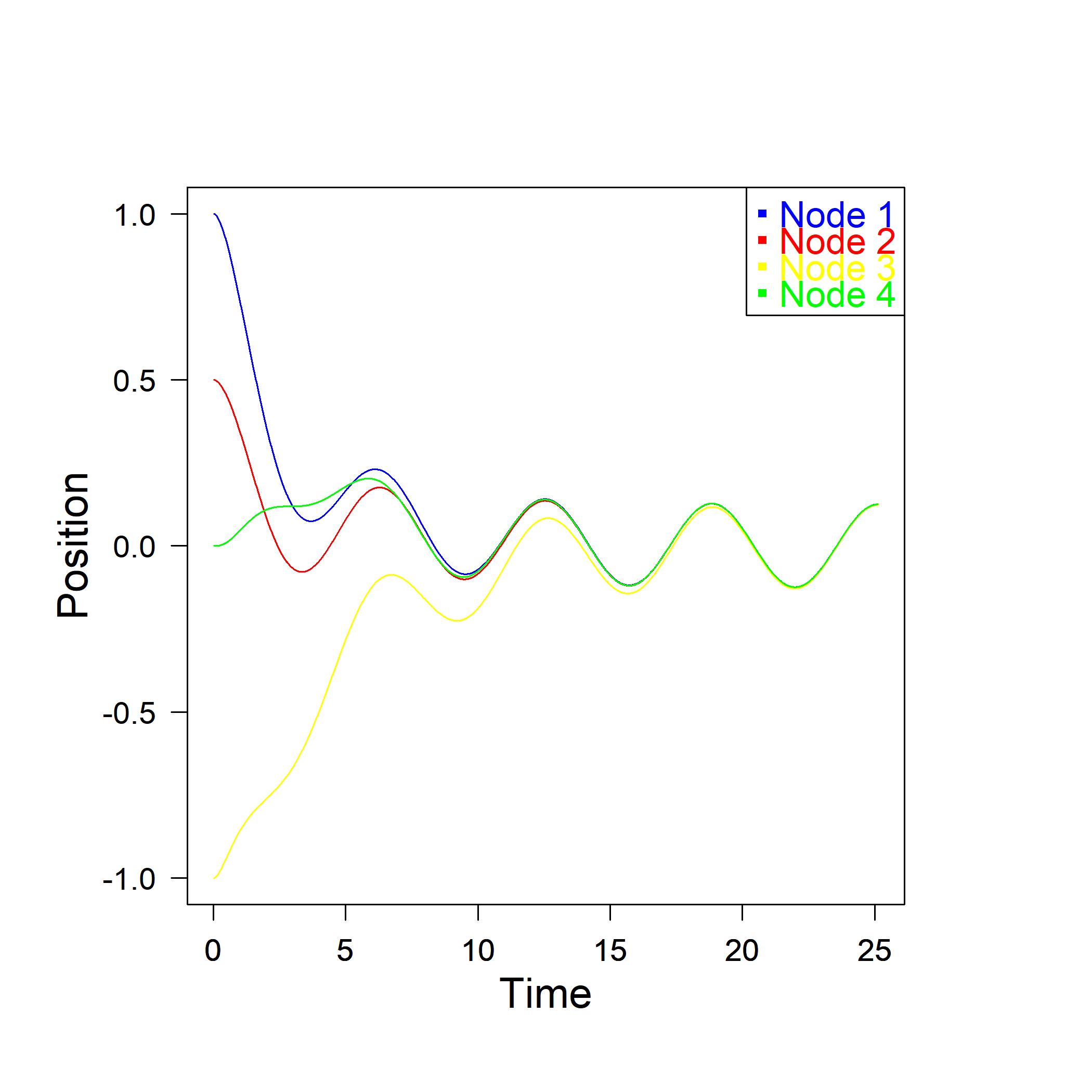}}
	\\
	\subfloat[]{\includegraphics[width=0.30\textwidth]{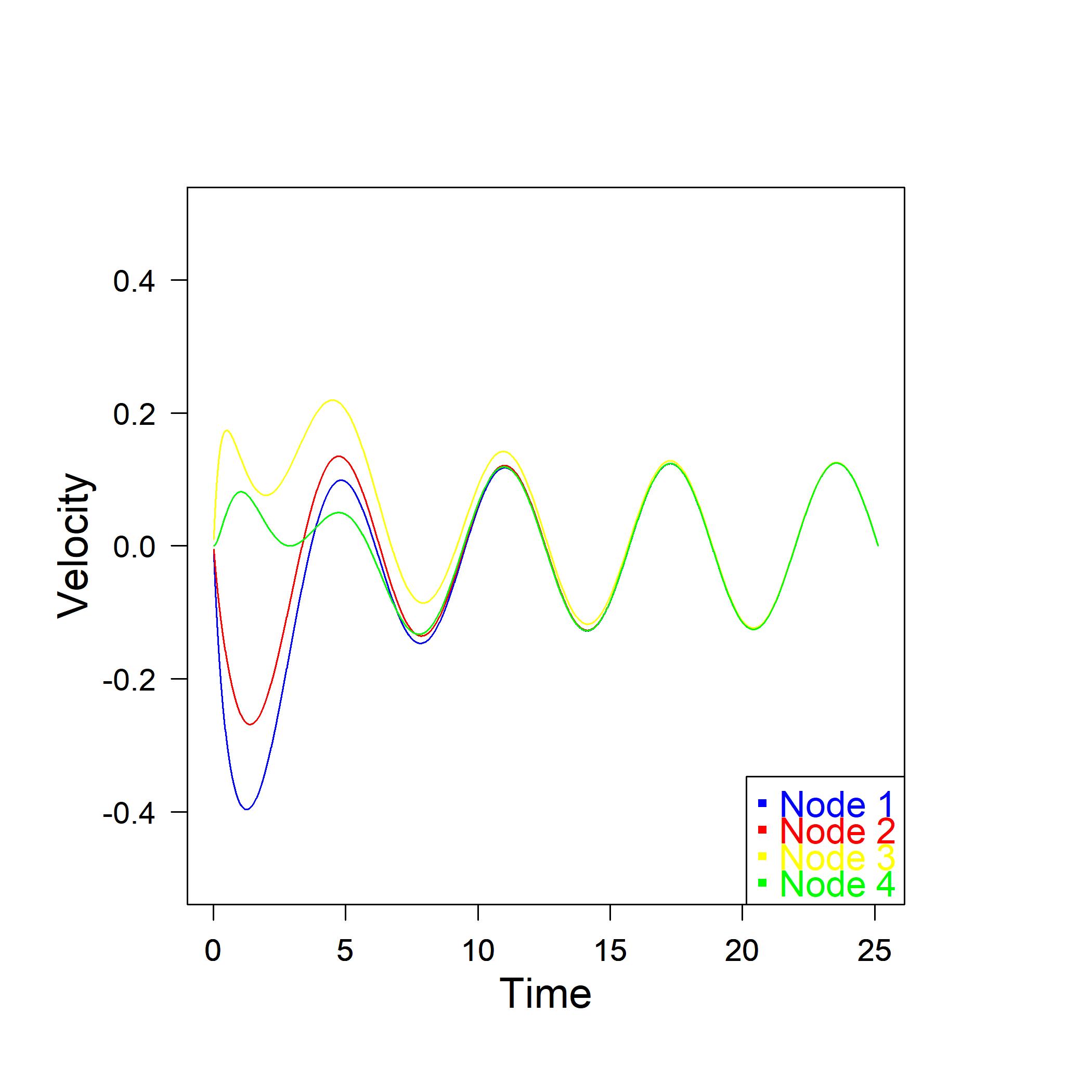}}
	\\
	\subfloat[]{\includegraphics[width=0.30\textwidth]{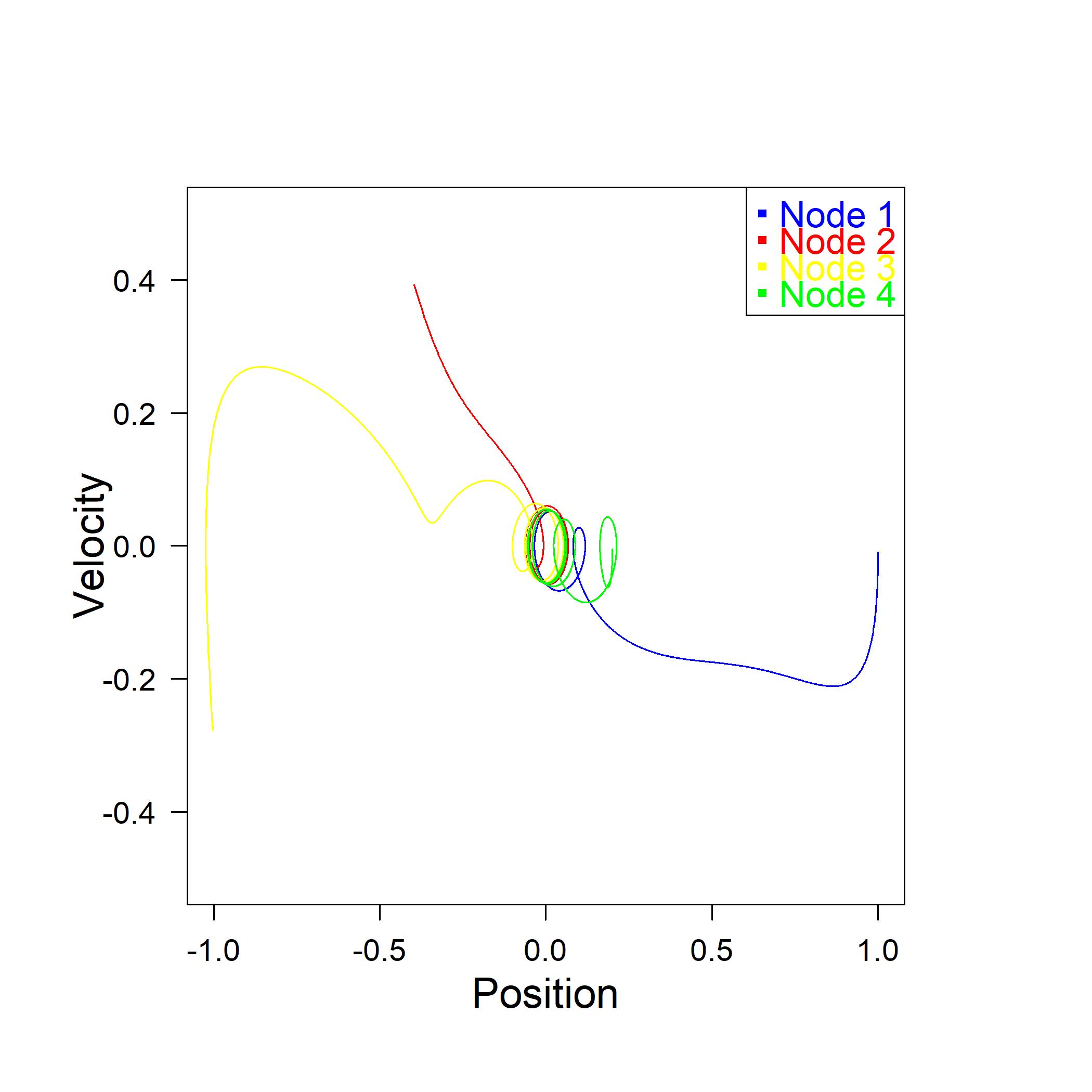}}
\caption{Synchronization of the harmonic oscillators in the example network with initial conditions ${\bf x}_{0}=[1,0.5,-1,0]^T$ and ${\bf v}_{0}=[0,0,0,0]^T$: panel (a) position, panel (b) velocity, and panel (c) phase portrait.}
	\label{fig4} 
\end{figure}

\begin{figure}[H]
\centering
	\subfloat[]{\includegraphics[width=0.30\textwidth]{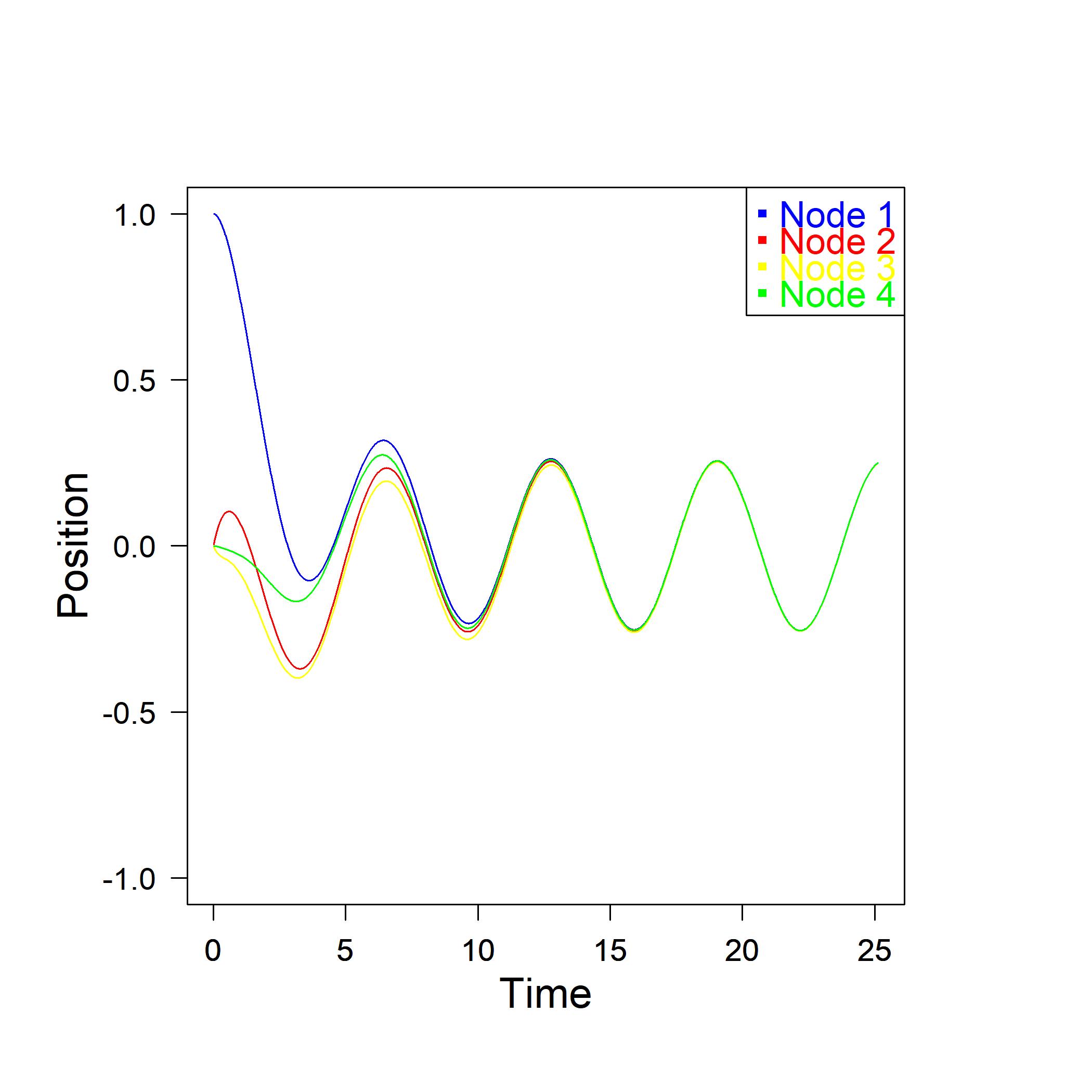}}
	\\
	\subfloat[]{\includegraphics[width=0.30\textwidth]{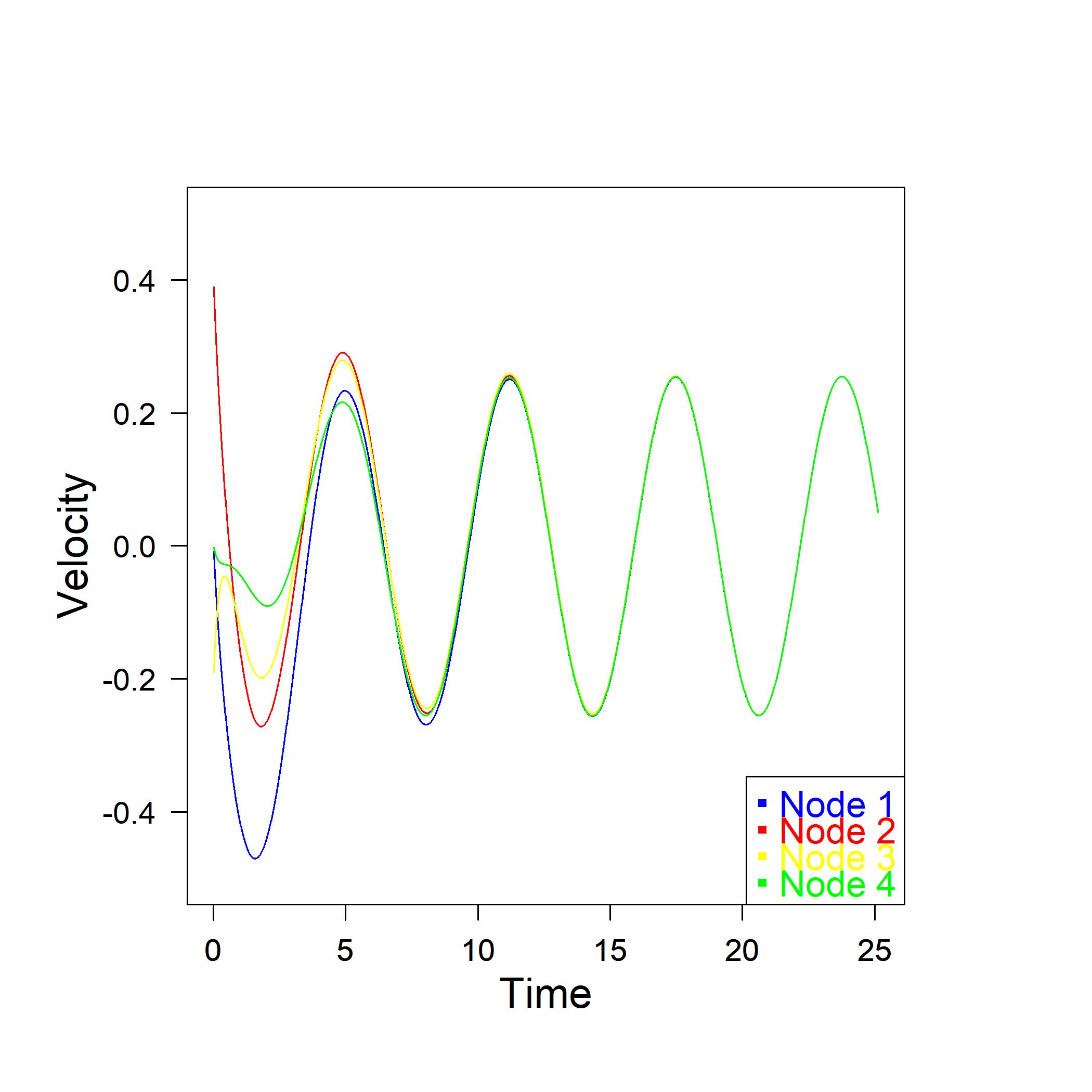}}
	\\
	\subfloat[]{\includegraphics[width=0.30\textwidth]{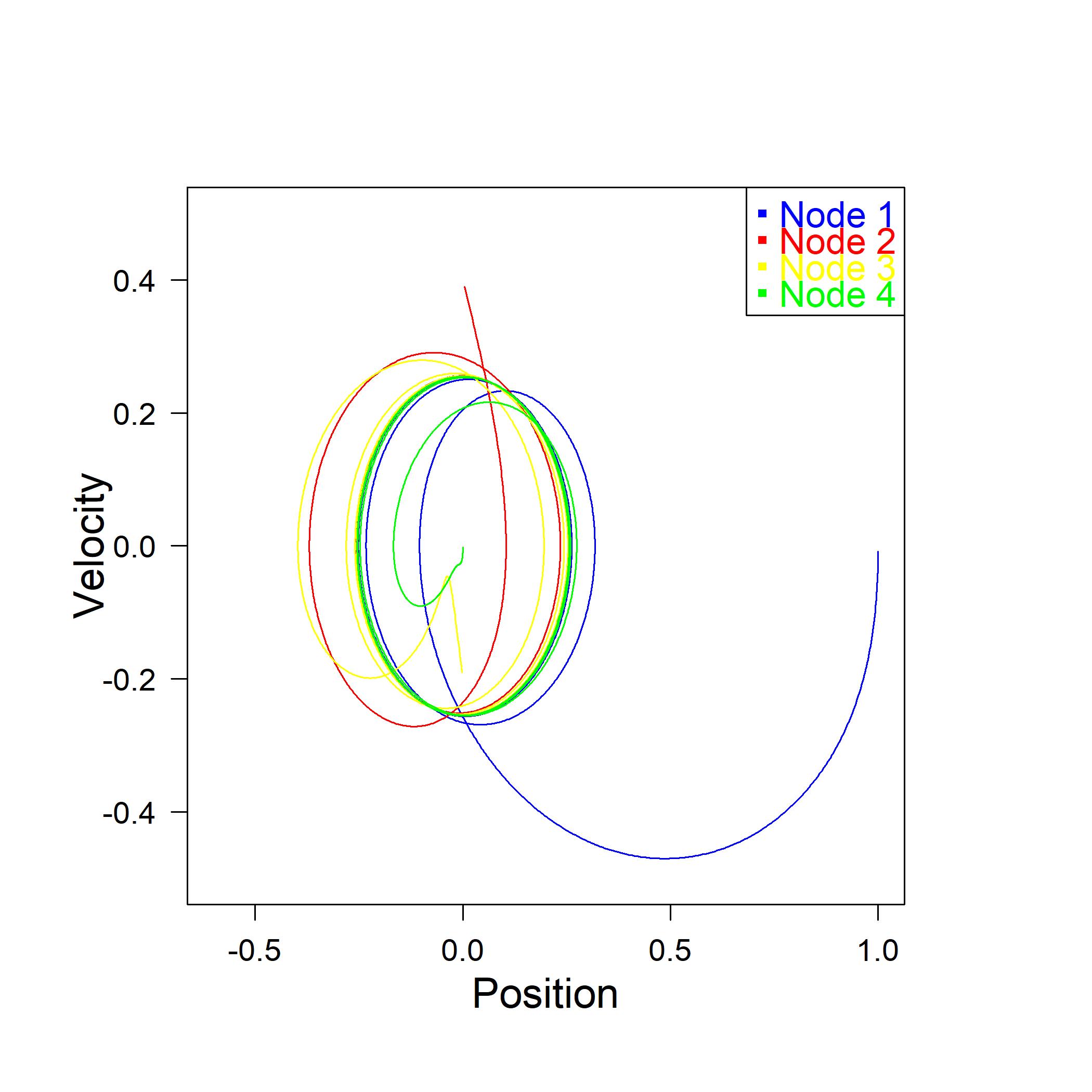}}
\caption{Synchronization of the harmonic oscillators in the example network with initial conditions ${\bf x}_{0}=[1,0,0,0]^T$ and ${\bf v}_{0}=[0,0.4,-0.2,0]^T$: panel (a) position, panel (b) velocity, and panel (c) phase portrait.}
	\label{fig5} 
\end{figure}

\begin{figure}[H]
\centering
	\subfloat[]{\includegraphics[width=0.30\textwidth]{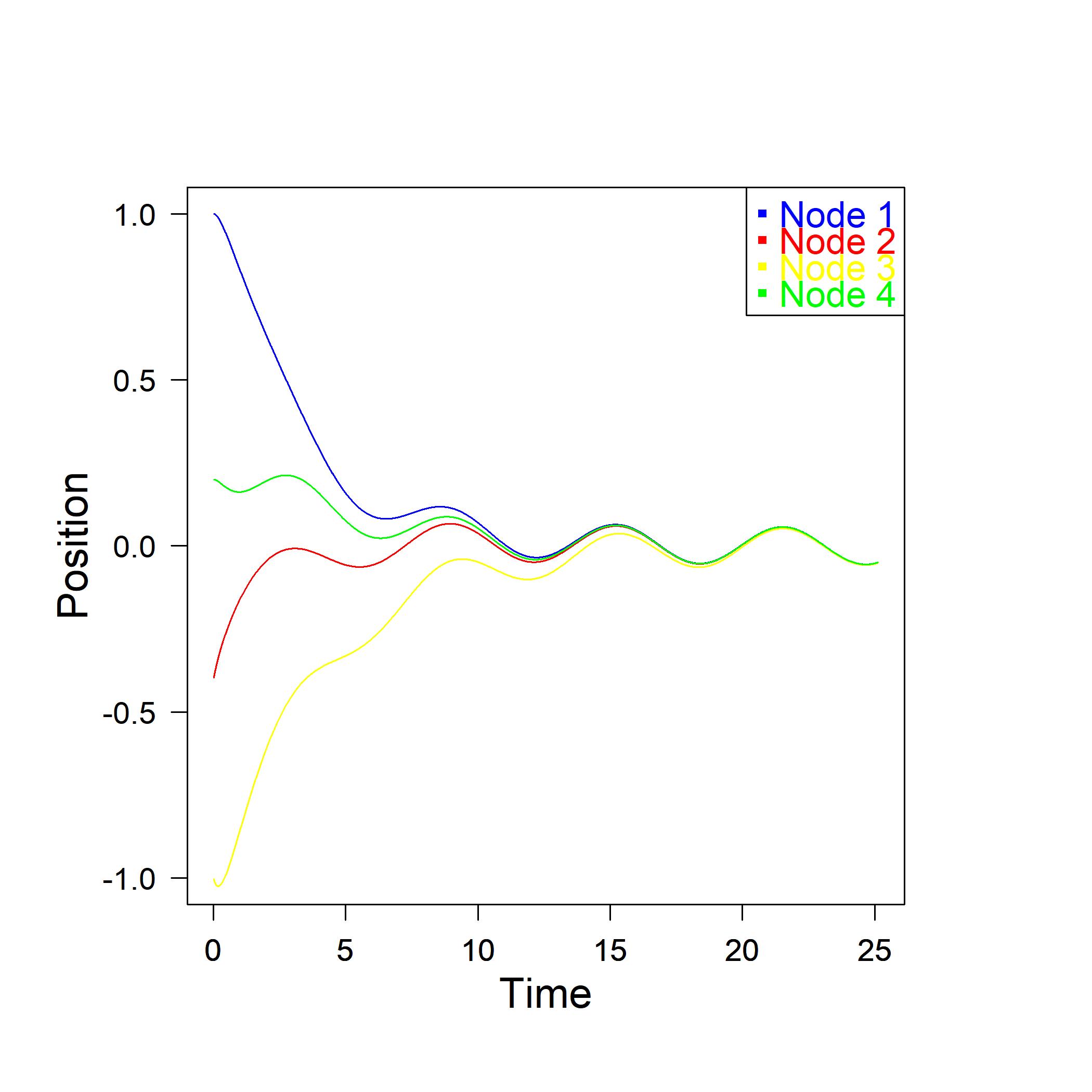}}
	\\
	\subfloat[]{\includegraphics[width=0.30\textwidth]{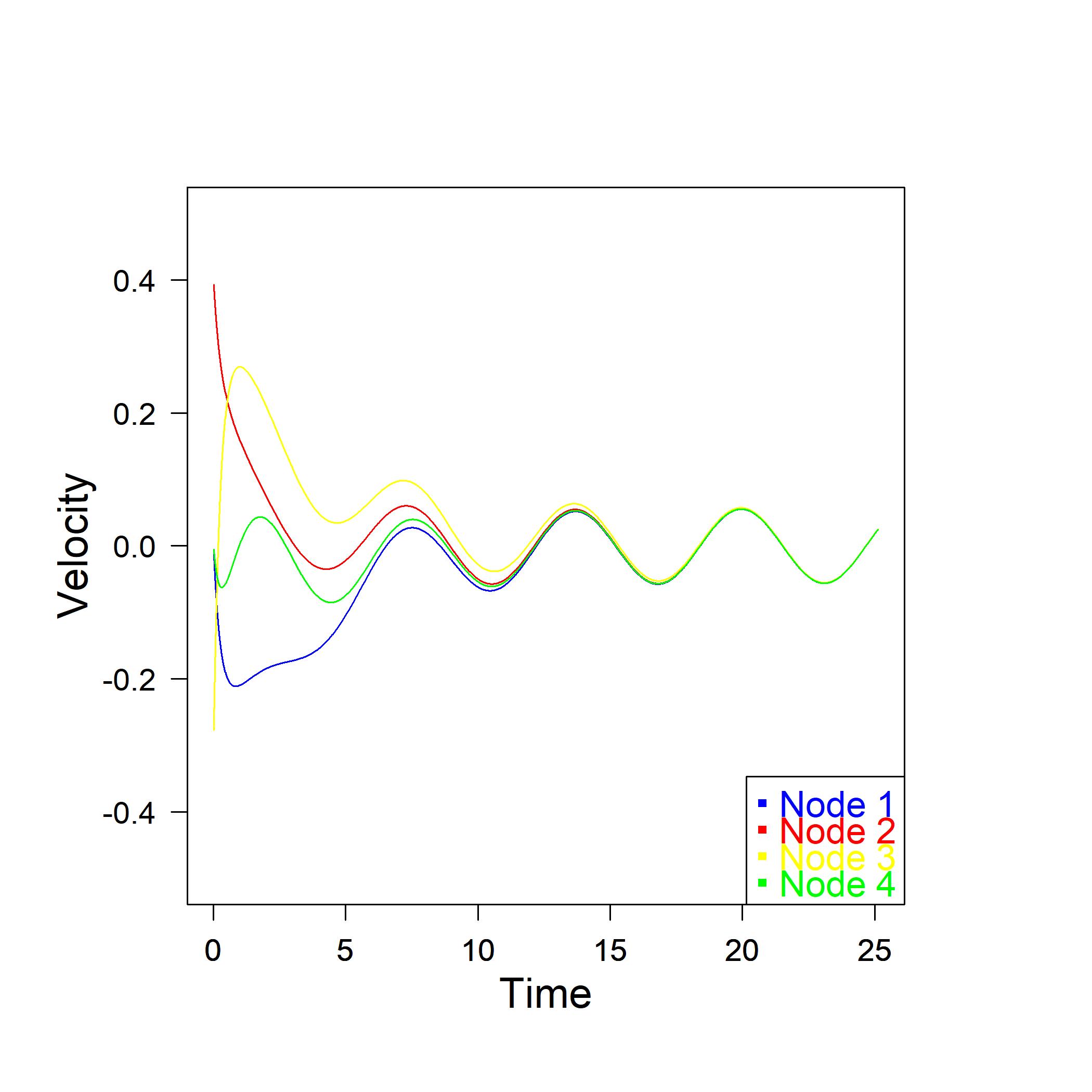}}
	\\
	\subfloat[]{\includegraphics[width=0.30\textwidth]{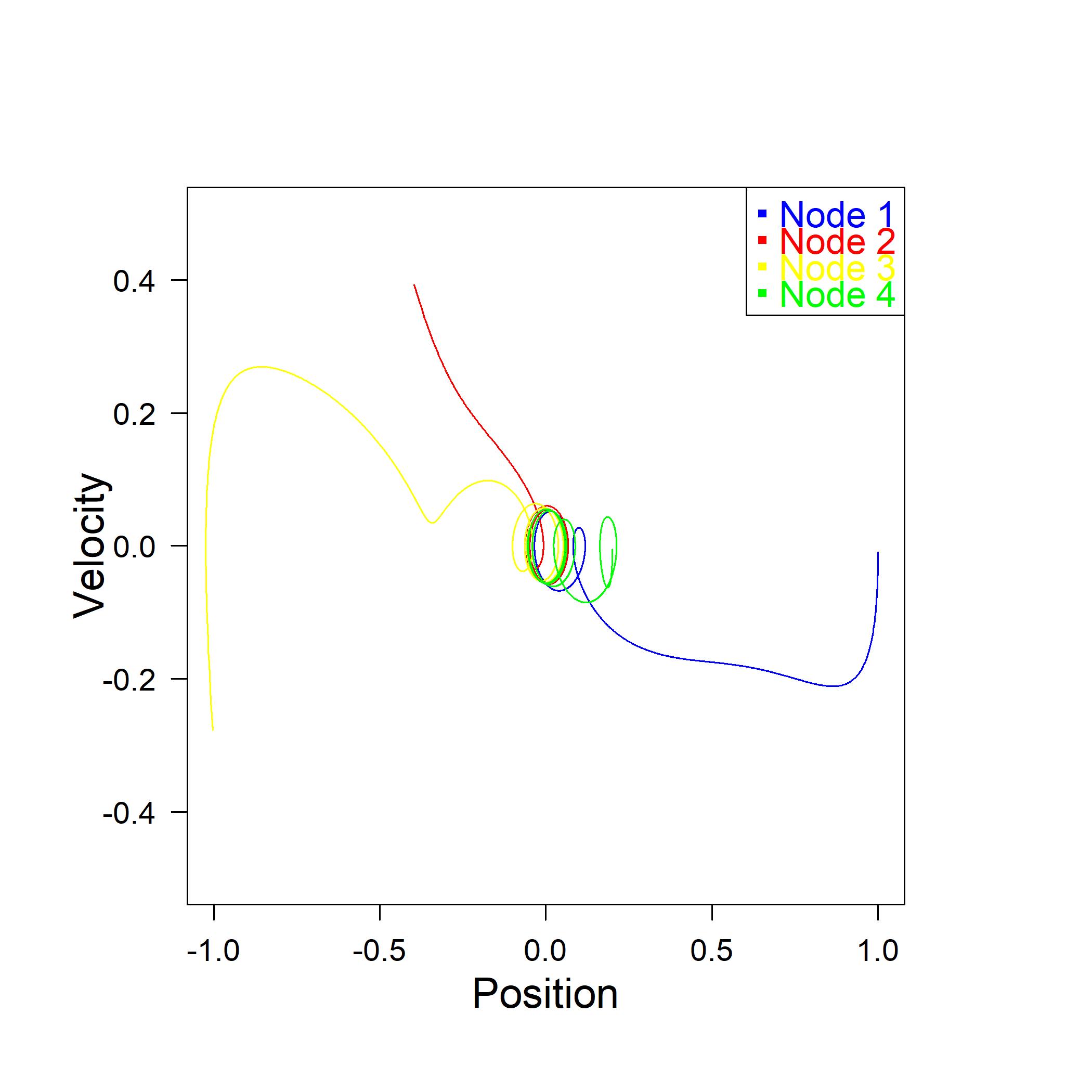}}
\caption{Synchronization of the harmonic oscillators in the example network with initial conditions ${\bf x}_{0}=[1,-0.4,-1,0.2]^T$ and ${\bf v}_{0}=[0,0.4,-0.3,0]^T$: panel (a) position, panel (b) velocity, and panel (c) phase portrait.}
	\label{fig6} 
\end{figure}

Let us consider again the toy model in Fig. \ref{fig1}. The four eigenvalues of $\bf L$ are: $\mu_1=4,\ \mu_2=3,\ \mu_3=1,\ \mu_4=0$. The eight eigenvalues of operators $\bf G$ , $\bf P$ , $\bf U$ and the four angles defined in Eq. (\ref{angles}) in Appendix \ref{Appendix A} are listed in Table \ref{table3}
\begin{table}[H]
	\begin{center}
		\begin{tabular}{||l|l|l|l||l||}
			\hline
			$\lambda$ & $\bf G$ & $\bf P$ & $\bf U$ & {\rm Angles}  \tabularnewline \hline
			${\lambda^{\rm + }_1}$ & $-0.268$ & $4.236$  & $-0.894+0.447i$ & ${\theta_1}=153^{\circ}$   \tabularnewline \hline
			${\lambda^{\rm + }_2}$ & $-0.382$ & $3.303$  & $-0.832+0.555i$ & ${\theta_2}=146^{\circ}$   \tabularnewline \hline
			${\lambda^{\rm + }_3}$ & $-0.500+0.866i$ & $1.618$  & $-0.447+0.894i$ & ${\theta_3}=117^{\circ}$ \tabularnewline \hline
			${\lambda^{\rm + }_4}$ & $+i$ & $1$  & $+i$ & ${\theta_4}=90^{\circ}$ \tabularnewline \hline
			${\lambda^{\rm - }_4}$ & $-i$ & $1$  & $-i$ &  $\quad \ \ -$ \tabularnewline \hline
			${\lambda^{\rm - }_3}$ & $-0.500-0.866i$ & $0.618$  & $-0.447-0.894i$ &  $\quad \ \ -$ \tabularnewline \hline
			${\lambda^{\rm - }_2}$ & $-2.618$ & $0.303$  & $-0.832-0.555i$ &   $\quad \ \ -$ \tabularnewline \hline
			${\lambda^{\rm - }_1}$ & $-3.732$ & $0.236$  & $-0.894-0.447i$ &  $\quad \ \ -$ \tabularnewline \hline
		\end{tabular}
	\end{center}
\caption{Eigenvalues of the matrices $\bf G$, $\bf P$ and $\bf U$ and angles ${\theta_i}=\arg(\lambda^{\rm U +}_i)= 2\arctan {\lambda^{\rm P +}_i}$ for the toy model.}
\label{table3}
\end{table}
In the Figures \ref{fig3}, \ref{fig4}, \ref{fig5}, and \ref{fig6} we illustrate some different synchronization processes on the same network, with different initial conditions, as specified in the caption of the figures. For each initial condition, plots represent the position and velocity of the four nodes and the phase portrait with the corresponding limiting cycle.

\subsection{Synchronization time}
\label{Synchronization_time}
The preceding analysis gives us the opportunity to add some insights into the widely discussed topic of synchronization time and speed. What follows contributes to this topic by highlighting the role of network topology in the synchronization process in order to compare the behavior of different systems and to discriminate which systems converge more or less rapidly to the synchronized state.

In the literature, a measure of synchronization is defined as ${\cal C}(t)=\frac{1}{n}\sum_{i=1}^{n}{\Theta}\left(\varepsilon - d_{i}(t)\right)$ where $d_{i}(t)=||{y_i}(t)-{\tilde y}_{i}(t)||$, ${\bf {\tilde y}}(t)$ is the synchronization solution in Eq. (\ref{asymptoticsolution}), $\varepsilon$ is a threshold and ${ \Theta}$ is the step function (see, for instance, \citet{Grabow2011}). This index provides the percentage of nodes whose distance from the synchronization solution is less than a given threshold. As the network becomes completely synchronized on this solution, ${\cal C}(t)\to 1$.

There is actually no univocal definition of an index that measures the speed of synchronization. For example, it could also be taken as an indicator of the time at which the maximum of the absolute difference between couples of the oscillators positions, that is $\max_{ij}\left| x_{i}(t)-x_{j}(t) \right|$, or alternatively, the mean of all the differences between couples of positions, drops below a given threshold.

In light of the closed expression of the asymptotic solution in Eq. (\ref{asymptoticsolution}), we will refer directly to the mean absolute difference between positions of nodes and such a limiting solution.

We aim here at defining a proper synchronization time as a function of the topological properties of the underlying network.
Since, in the singular value decomposition, the more negative are the eigenvalues, the faster the corresponding exponentials go to zero, we have to select exponentials with less negative exponents if we are to identify the dominant behavior that slows synchronization and to set a bound on the synchronizability of the global process. Let us consider an initial state given by
${\bf y}_{0}=\left(
{\bf x}_{0}, 
{\bf v}_{0} 
\right)^{T}$
and let us write the expansion of the solution to problem (\ref{system3}) on the basis of the eigenvectors of matrix $\bf G$:

\begin{equation*}
\begin{split}
\label{expansion}
& \quad {\bf y}(t)\\
=& \quad e^{{\bf G}t}{\bf y}_{0}\\
=& \quad e^{t\lambda_{n}^{G+}}\left(\psi_{n}^{G+}\cdot {\bf y}_{0} \right)\psi_{n}^{G+}
 +e^{t\lambda_{n}^{G-}}\left(\psi_{n}^{G-}\cdot {\bf y}_{0} \right)\psi_{n}^{G-}+\\
&\quad\ e^{t\lambda_{n-1}^{G+}}\left(\psi_{n-1}^{G+}\cdot {\bf y}_{0} \right)\psi_{n-1}^{G+}
 +e^{t\lambda_{n-1}^{G-}}\left(\psi_{n-1}^{G-}\cdot {\bf y}_{0} \right)\psi_{n-1}^{G-}
 +\dots\\
=& \quad e^{it}\frac{1}{\sqrt{2n}}
\left(\begin{array}{c}
{\bf u} \\
i{\bf u} \\
\end{array} \right)^{\dagger}
\cdot
\left(\begin{array}{c}
{\bf x}_{0} \\
{\bf v}_{0} \\
\end{array} \right)
\cdot
\frac{1}{\sqrt{2n}}
\left(\begin{array}{c}
{\bf u} \\
i{\bf u} \\
\end{array} \right)+\\& \quad
e^{-it}\frac{1}{\sqrt{2n}}
\left(\begin{array}{c}
{\bf u} \\
-i{\bf u} \\
\end{array} \right)^{\dagger}
\cdot
\left(\begin{array}{c}
{\bf x}_{0} \\
{\bf v}_{0} \\
\end{array} \right)
\cdot
\frac{1}{\sqrt{2n}}
\left(\begin{array}{c}
{\bf u} \\
-i{\bf u} \\
\end{array} \right)+\\
&\quad\ e^{t\lambda_{n-1}^{G+}}\left(\psi_{n-1}^{G+}\cdot {\bf y}_{0} \right)\psi_{n-1}^{G+}
+e^{t\lambda_{n-1}^{G-}}\left(\psi_{n-1}^{G-}\cdot {\bf y}_{0} \right)\psi_{n-1}^{G-}
+\dots\\
=& \quad \frac{1}{2n}e^{it}\left[{\bf u}^{T} {\bf x}_{0}-i{\bf u}^{T}{\bf v}_{0} \right]
\left(\begin{array}{c}
{\bf u} \\
i{\bf u} \\
\end{array} \right)+\\& \quad
\frac{1}{2n}e^{-it}\left[{\bf u}^{T} {\bf x}_{0}+i{\bf u}^{T}{\bf v}_{0} \right]
\left(\begin{array}{c}
{\bf u} \\
-i{\bf u} \\
\end{array} \right)+\\
&\quad\ e^{t\lambda_{n-1}^{G+}}\left(\psi_{n-1}^{G+}\cdot {\bf y}_{0} \right)\psi_{n-1}^{G+}
+e^{t\lambda_{n-1}^{G-}}\left(\psi_{n-1}^{G-}\cdot {\bf y}_{0} \right)\psi_{n-1}^{G-}
+\dots\\
=& \quad \frac{1}{n}\left[\begin{array}{c}
(+\cos t\, {\bf u}^{T} {\bf x}_{0}+\sin t\, {\bf u}^{T} {\bf v}_{0}){\bf u} \\
(-\sin t\, {\bf u}^{T} {\bf x}_{0}+\cos t\, {\bf u}^{T} {\bf v}_{0}){\bf u} \\
\end{array} \right]+\\
&\quad\ e^{t\lambda_{n-1}^{G+}}\left(\psi_{n-1}^{G+}\cdot {\bf y}_{0} \right)\psi_{n-1}^{G+}
+e^{t\lambda_{n-1}^{G-}}\left(\psi_{n-1}^{G-}\cdot {\bf y}_{0} \right)\psi_{n-1}^{G-}
+\dots
\end{split}
\end{equation*}

The first term in the last step is equal to ${\bf {\tilde y}}(t)$ in Eq. (\ref{asymptoticsolution}) and it is the synchronization solution, so that
\begin{equation}
\begin{split}
\label{gap}
&{\bf y}(t)-{\bf {\tilde y}}(t)=\\
&e^{t\lambda_{n-1}^{G+}}\left(\psi_{n-1}^{G+}\cdot {\bf y}_{0} \right)\psi_{n-1}^{G+}
+e^{t\lambda_{n-1}^{G-}}\left(\psi_{n-1}^{G-}\cdot {\bf y}_{0} \right)\psi_{n-1}^{G-}
+\dots
\end{split}
\end{equation}

All the $2n-2$ terms on the right hand side contain an exponential whose coefficient is $\lambda_{i}^{G\pm},\ i=n-1, \dots,2,1$. These coefficients can be complex or real but in both cases their real part is always negative, according to Table \ref{table1}. Therefore the process is governed by the term that contains the least negative eigenvalue real part. The eigenvalues of $\bf G$ are expressed by Eq. (\ref{eigenvalues_G}), $\lambda^{\rm G\pm}_i=\frac{1}{2}\left[-\mu_i\pm \sqrt{\mu_i^2-4}\right]$, as an increasing function of $\mu_i$. Since $\mu_1 >k_{\rm max}\geq 2$ for any nontrivial network, there is always a real negative eigenvalue and the maximum real value is attained for $i=1$, so that the less negative real eigenvalue is $\lambda^{\rm G +}_1=\frac{1}{2}[-\mu_1+ \sqrt{\mu_1^2-4}]$. On the other hand, when $\lambda^{\rm G\pm}_i\in {\mathbb C}$, the less negative real part is gained for the minimum non zero $\mu_i$, i.e. for $i=n-1$, corresponding to the algebraic connectivity, and it is given by ${\rm Re}(\lambda^{\rm G\pm}_{n-1})=-\frac{\mu_{n-1}}{2}$. Therefore, the less negative exponential coefficient is given by
\begin{equation}
\label{max_eigenvalue}
\lambda_{S}=\max \left(\frac{1}{2}\left[-\mu_1+ \sqrt{\mu_1^2-4}\right]; -\frac{\mu_{n-1}}{2} \right)
\end{equation}

If we look at Eq. (\ref{gap}) by components
\begin{equation}
\begin{split}
\label{gap_components}
&{y_i}(t)-{{\tilde y}_i}(t)=\\
&e^{t\lambda_{n-1}^{G+}}\left(\psi_{n-1}^{G+}\cdot {\bf y}_{0} \right)\psi_{n-1}^{G+}(i)
+e^{t\lambda_{n-1}^{G-}}\left(\psi_{n-1}^{G-}\cdot {\bf y}_{0} \right)\psi_{n-1}^{G-}(i)
+\dots
\end{split}
\end{equation}

then we can say that, as $t\to +\infty$,
\begin{equation}
\label{main_term}
\left|{y_i}(t)-{{\tilde y}_i}(t)\right|\sim e^{t\lambda_{S}}\left(\psi_{S}\cdot {\bf y}_{0} \right)\psi_{S}(i)
\end{equation}

where $\lambda_{S}$ is the eigenvalue in Eq. (\ref{max_eigenvalue}) and $\psi_{S}(i)$ is the component $i$ of the corresponding eigenvector (as far as this component is not zero). Now, we want to find out the time after which, for any node $i$, the difference between the actual solution ${y_i}(t)$ and the synchronization solution ${{\tilde y}_i}(t)$ becomes less than a given threshold $\varepsilon>0$: $e^{t\lambda_{S}}\left(\psi_{S}\cdot {\bf y}_{0} \right)\psi_{S}(i)< \varepsilon$. By solving for $t$, we get:
\begin{equation}
\label{time}
t_{i}>\frac{1}{\left|\lambda_{S}\right|}\log \left| \frac{\left(\psi_{S}\cdot {\bf y}_{0} \right)\psi_{S}(i)}{\varepsilon} \right|
\end{equation}

Finally, the average synchronization time can be bounded by
\begin{equation}
\label{mean_time}
t_{\rm mean}>\frac{1}{n\left|\lambda_{S}\right|}\sum_{i=1}^{n}\log \left| \frac{\left(\psi_{S}\cdot {\bf y}_{0} \right)\psi_{S}(i)}{\varepsilon} \right|
\end{equation}

\begin{remark}
When $\left|\lambda_{S}\right|=\frac{\mu_{n-1}}{2}$, the greater the algebraic connectivity, the smaller is the synchronization time. Moreover, the algebraic connectivity is bounded by the inequality $\mu_{n-1}\leq \frac{n}{n-1}k_{\rm min}$ with $k_{\rm min}$ minimum degree in the network, and since the density of the network is given by $\delta =\frac{k_{\rm mean}}{n-1}$, we have that
\begin{equation}
\mu_{n-1}\leq \frac{n}{n-1}k_{\rm min} \leq \frac{n}{n-1}k_{\rm mean} = n \delta
\end{equation}

so that, in this case, we have
\begin{equation}
\begin{split}
\label{mean_time_vs_density}
t>&\frac{2}{n\left|\mu_{n-1}\right|}\sum_{i=1}^{n}\log \left| \frac{\left(\psi_{S}\cdot {\bf y}_{0} \right)\psi_{S}(i)}{\varepsilon} \right|\\ \geq& 
\frac{2}{n^2\delta}\sum_{i=1}^{n}\log \left| \frac{\left(\psi_{S}\cdot {\bf y}_{0} \right)\psi_{S}(i)}{\varepsilon} \right|
\end{split}
\end{equation}

For $n$ fixed then, as $\delta$ grows then the synchronization time decreases.
\end{remark}

\begin{figure}[H]
	\includegraphics[width=1\linewidth]{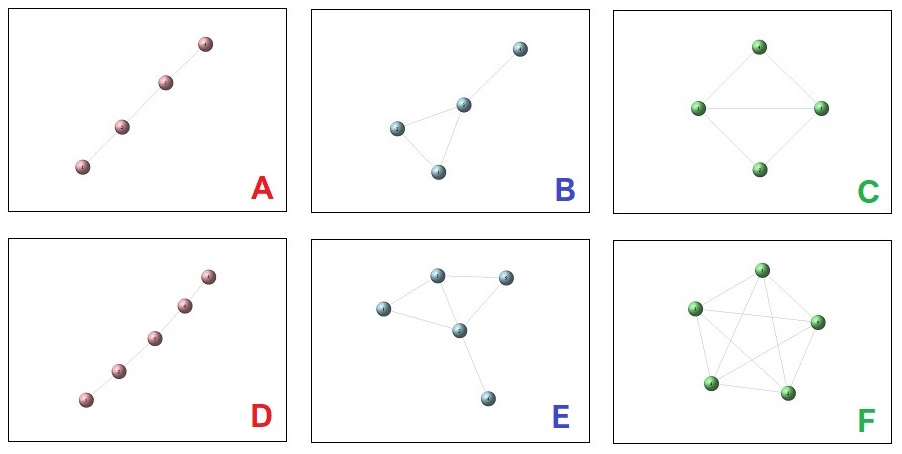}
	\centering
	\caption{Small graphs with different synchronization times analyzed in the text.}
	\label{fig7}      
\end{figure}

Let us test the synchronization time in Eq. (\ref{mean_time}) on the toy networks with $n=4$ and $n=5$ nodes in figure \ref{fig7}.
For the networks A, B and C with $n=4$, the values of the parameter $\lambda_{S}$ are: $\lambda_{S}^{\rm A}=-0.2928932$, $\lambda_{S}^{\rm B}=-0.2679492\ (-0.3819660)$ and $\lambda_{S}^{\rm C}=-0.2679492\ (-0.2679492)$.\footnote{Note that, since the two less negative eigenvalues for network B and C are equal, we have to look at the second less negative eigenvalues in round brackets.} Then we expect that network A would synchronize faster than B, which in turn would synchronize faster than C.

\begin{figure}[H]
\centering
	\subfloat[]{\includegraphics[width=0.45\textwidth]{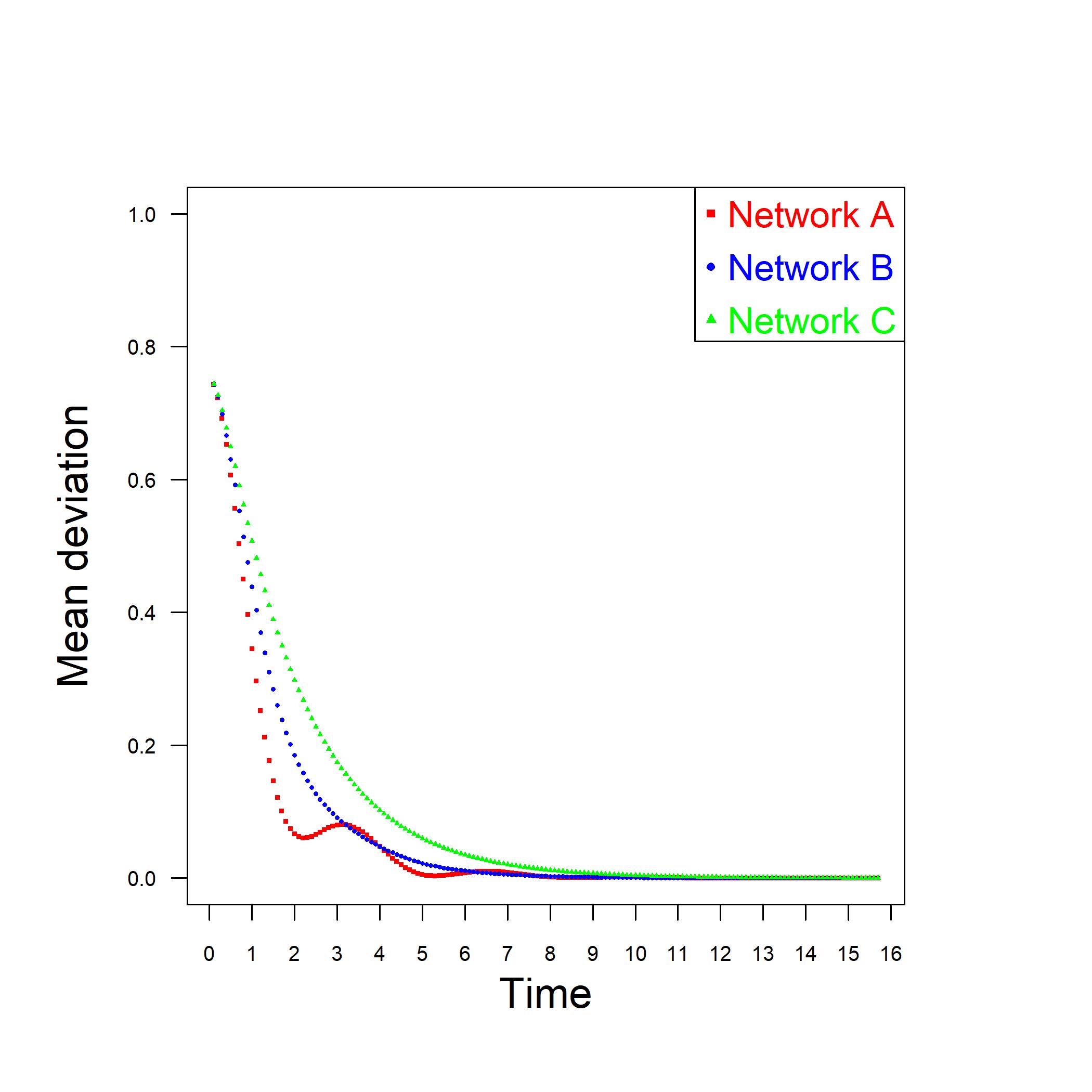}}
	\\
	\subfloat[]{\includegraphics[width=0.45\textwidth]{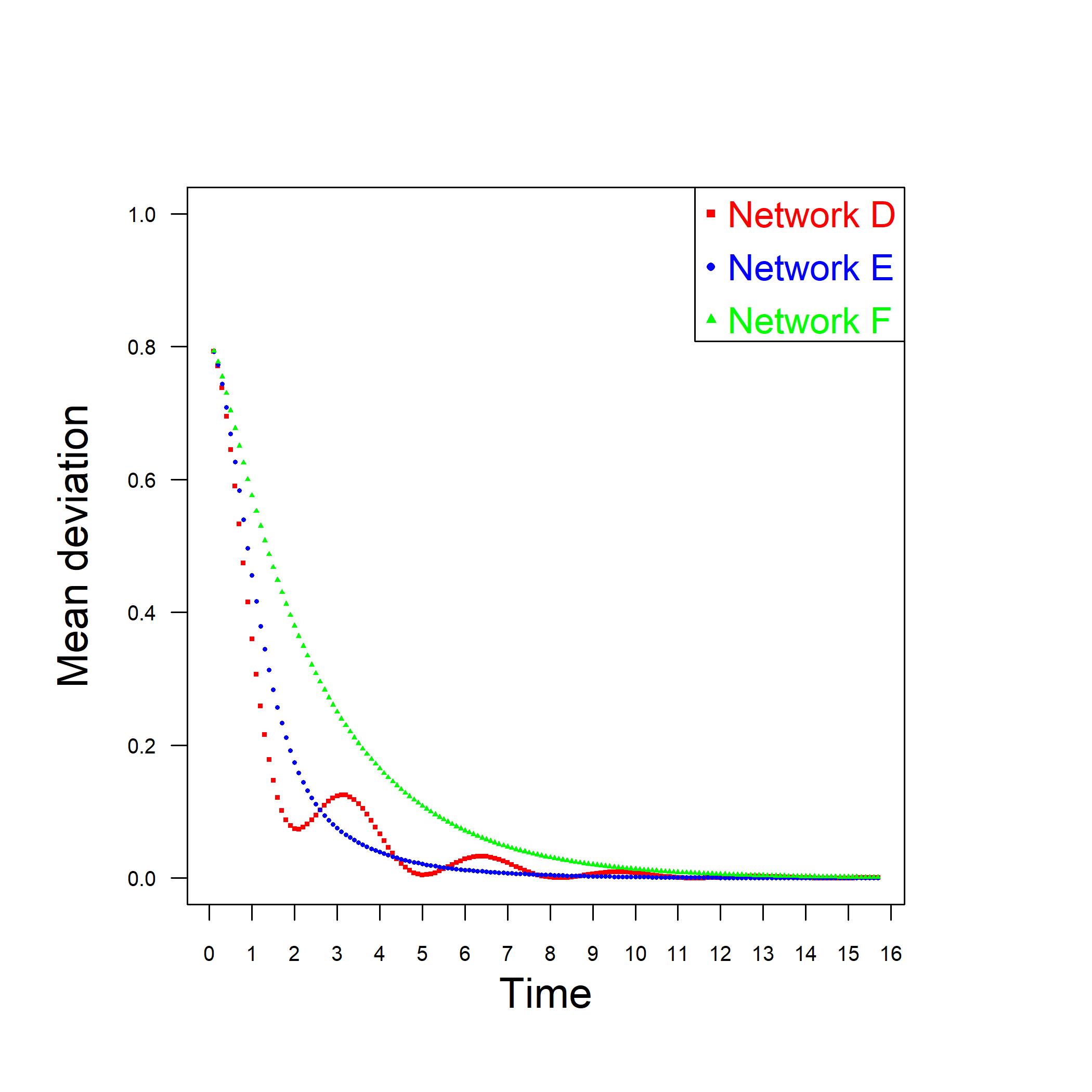}}
	\caption{Different speed of synchronization for the two sets of toy networks: panel (a) A, B and C; panel (b) D, E and F.}
	\label{fig8} 
\end{figure}

In figure \ref{fig8}, panel (a), we plot the mean of the absolute values of differences between node position at time $t$ and the synchronization solution, as in Eq. (\ref{gap}), as a function of time $t$: we observe that the red line (network A) decays on average more rapidly than the blue line (network B), which decays more rapidly than the green line (network C). If we compute the mean synchronization times given in formula (\ref{mean_time}) with $\varepsilon=10^{-3}$, we get: $t_{\rm mean}^{\rm A}=8.421$, $t_{\rm mean}^{\rm B}=17.273$ and $t_{\rm mean}^{\rm C}=18.264$. The empirical values, computed as the mean of the minimum times $t$ such that the absolute value of the difference between the solution for a given node and the asymptotic solution goes below $\epsilon$, are $t_{\rm emp}^{\rm A}=8.275$, $t_{\rm emp}^{\rm B}=17.525$ and $t_{\rm emp}^{\rm C}=21.925$.

Similarly, for networks D, E and F with $n=5$, the values of the parameter $\lambda_{S}$ are: $\lambda_{S}^{\rm A}=-0.1909830$, $\lambda_{S}^{\rm B}=-0.2087122\ (-0.2679492)$ and $\lambda_{S}^{\rm C}=-0.2087122\ (-0.2087122)$. Then we expect that network D would synchronize faster than E, which in turn would synchronize faster than F.
In figure \ref{fig8}, panel (b), we plot again the mean deviation of the node positions at time $t$ and the synchronization solution, as in Eq. (\ref{gap}), as a function of time $t$: we observe that red line (network D) decays on average more rapidly than blue line (network E), which decays more rapidly than green line (network F). Again, if we compute the mean synchronization times in formula (\ref{mean_time}) with $\varepsilon=10^{-3}$, we get: $t_{\rm mean}^{\rm D}=11.849$, $t_{\rm mean}^{\rm E}=19.868$ and $t_{\rm mean}^{\rm F}=26.510$,
whereas the empirical values are $t_{\rm emp}^{\rm D}=11.140$, $t_{\rm emp}^{\rm E}=21.760$ and $t_{\rm emp}^{\rm F}=26.940$.

This observation could be taken as a hint that a path will synchronize faster than a more integrated or even complete network.
Things are actually more elaborate as argued in the following important remark.

\begin{figure}[H]
	\includegraphics[width=1\linewidth]{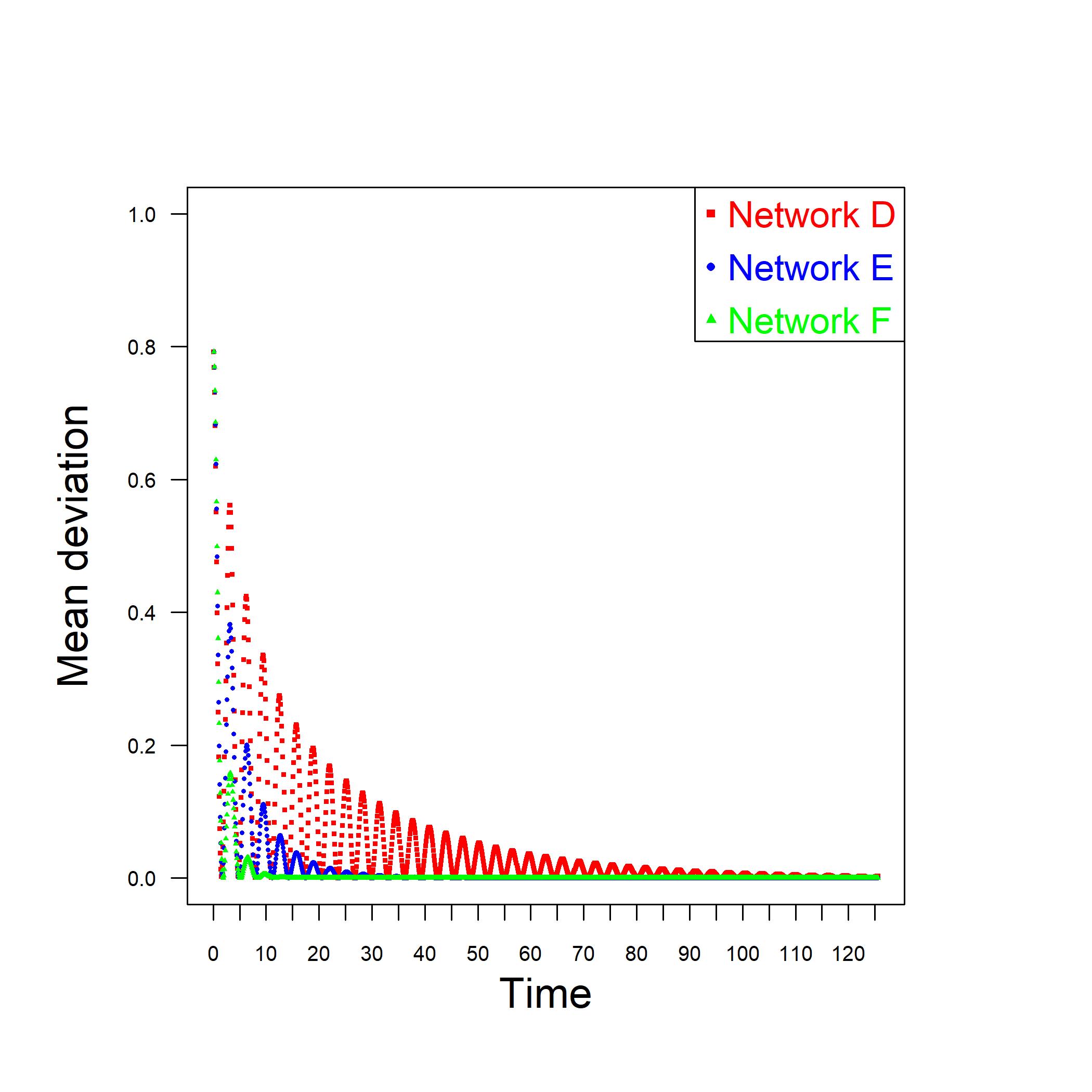}
	\centering
	\caption{Synchronization for weak coupled networks}
	\label{fig9}      
\end{figure}

\begin{remark}
So far we have assumed a constant value equal to $1$ for the coupling coefficient $c_{2}^{\prime}$ between nodes. In this remark we want to show that different networks may exhibit different rankings in the speed and time of synchronization according to the values of the coupling coefficient. Let us consider again the triplet of networks D, E and F in figure \ref{fig7} and let us assume as synchronization time the lower bound in formula (\ref{mean_time}) with a threshold equal to $\varepsilon=10^{-6}$.

When coefficients in system (\ref{system2}) are $c_{1}=c^{\prime}_{2}=1$ and $c_{2}=c^{\prime}_{1}=0$, formula (\ref{mean_time}) yields
\begin{equation}
t_{\rm mean}^{D}=48.02 \quad t_{\rm mean}^{E}=52.97 \quad t_{\rm mean}^{F}=59.61
\end{equation}

and the empirical values are
\begin{equation}
t_{\rm emp}^{D}=38.10 \quad t_{\rm emp}^{E}=49.74 \quad t_{\rm emp}^{F}=60.04
\end{equation}

These values show again that the path network gains synchronization more rapidly than the complete network. However, if we lower the coupling coefficient to $c^{\prime}_{2}=0.1$, the same times become:
\begin{equation}
t_{\rm mean}^{A}=278.57 \quad t_{\rm mean}^{B}=82.20 \quad t_{\rm mean}^{C}=43.07
\end{equation}

with empirical values
\begin{equation}
t_{\rm emp}^{A}=294.50 \quad t_{\rm emp}^{B}=105.70 \quad t_{\rm emp}^{C}=38.24
\end{equation}

The different behavior is illustrated in figure \ref{fig9}, which is the analogous of figure \ref{fig8}, panel (b), for the weak coupled networks D, E and F. This example shows that, for a weakly coupled system, the complete network achieves the state of perfect synchronization far more rapidly than the path network. This is an outcome opposite to the previous one. More surprisingly, the complete network shows a synchronization time which is shorter for weak coupling than for strong coupling. This fact makes it interesting to analyze how this time can depend on the coupling coefficients through the eigenvalues of the involved matrices.

If we set $c^{\prime}_{2}=\alpha$, with $\alpha>0$, the eigenvalues of $\bf G$ become
\begin{equation}
\lambda^{\rm G\pm}_i=\frac{1}{2}\left[-\alpha\mu_i\pm \sqrt{\alpha^2\mu_i^2-4}\right]
\end{equation}

As $\alpha$ decreases more eigenvalues become complex and if $\alpha<\frac{2}{\mu_{1}}$ they are all complex. Therefore, what matters is the real part which is just $-\frac{\alpha\mu_i}{2}$. The maximum non-null value is then for all $-\frac{\alpha\mu_{n-1}}{2}$. For the path graph $\mu_{n-1}=2\left[ 1-\cos(\pi/n)\right]$, for the complete graph $\mu_{n-1}=n$. Then we have to compare $-\alpha\left( 1-\cos(\pi/n)\right)$ and $-\frac{\alpha n}{2}$. Since $n>2\left( 1-\cos(\pi/n)\right)$ for every $n>2$, we have that $-\frac{\alpha n}{2}<-\alpha\left( 1-\cos(\pi/n)\right)$, which shows that, for $\alpha<\frac{2}{\mu_{1}}$, the high connectivity of the complete network enables it to synchronize faster than other topologies. Similar results for network A, B and C in figure \ref{fig7} confirm these conclusions.
\end{remark}

\section{Resonant states in complex networks}
\label{section5}

We move now to the third case in Table \ref{summary} by introducing a driving periodic force acting on or even produced by a single node. We will analyze the effect of a standard periodic sinusoidal force with frequency $\omega$. Our purpose is to illustrate some new results about resonant phenomena on networks. In order to introduce method and notations, we quickly revise some ideas about the single driven resonant oscillator.

\subsection{The driven harmonic oscillator}

Let us consider a harmonic oscillator subjected to an external periodic force ${\bf F}(t)$: $m\ddot{{\bf x}}={\bf F}(t)-k{\bf x}$. Let us set ${\bf y}$ as in Eq. (\ref{definitionofy}). Then $\dot {\bf y}={\bf A}{\bf y}+{\bf b}$ where ${\bf A}=\left( \begin{array}{cc} 
0 & 1 \\
 -\omega_{0}^{2} & 0\\
\end{array} \right)$, with
$\omega_0=\sqrt{\frac{k}{m}}$ and ${\bf b}=\frac{{\bf F}(t)}{m}\, {\bf e}\coloneqq f(t)\, {\bf e}$, with ${\bf e}=[0,1]^{T}$. The general solution, with the initial condition ${\bf y}_{0}=[{\bf x}_{0}, {\bf v}_{0}]^{T}$, is
\begin{equation}
{\bf y}(t)=e^{{\bf A}t}\left[ \int_{0}^{t}e^{-{\bf A}u}{\bf b}(u)du +{\bf c} \right]
\label{generalsolution6}
\end{equation}

where
\begin{equation}
e^{{\bf A}t}=\left[ \begin{array}{cc} 
\cos \omega_{0}t & \frac{1}{\omega_{0}}\sin \omega_{0}t \\
 -\omega_{0}\sin \omega_{0}t & \cos \omega_{0}t\\
\end{array} \right]
\label{expAt}
\end{equation}

Let us consider a driving force $F(t)=F_{0}\cdot \sin(\omega t)$ so that ${\bf b}(t)=f(t)\, {\bf e}$ with $f(t)=\frac{F_{0}}{m}\sin(\omega t)$. From Eq. (\ref{expAt}), we have
\begin{equation}
e^{-{\bf A}u}\,{\bf b}(u)=\frac{F_{0}}{m}\left[ \begin{array}{c} -\frac{1}{\omega_0}\sin \omega_{0}u \cdot \sin \omega u \\
 \cos \omega_{0}u \cdot \sin \omega u\\
\end{array} \right]
\end{equation}

and, by a straightforward computation of the integral in Eq. (\ref{generalsolution6}), we get
\begin{equation}
{\bf y}(t)=\frac{F_{0}}{m}\frac{\omega}{\omega_{0}^{2}-\omega^{2}}\left[ \begin{array}{c} 
\frac{\sin \omega t}{\omega} -\frac{\sin \omega_{0} t}{\omega_{0}}  \\
 \cos \omega t - \cos \omega_{0}t \\
\end{array} \right]
\label{synchronized_solution}
\end{equation}

Let us note that, for $\omega \to \omega_{0}$, the asymptotic behaviors
\begin{equation}
\begin{split}
x(t)&\sim \frac{F_{0}t}{2m\omega_{0}}\cos \omega_{0}t\\
v(t)&\sim \frac{F_{0}t}{2m}\sin \omega_{0}t
\end{split}
\end{equation}

show that the amplitude of both position and velocity grows linearly in time. 

Incidentally, it might be useful to mention that in the case of a periodical delta-type force $F(t)=F_{0}\sum_{k=1}^{n}\delta(t-kT)$, where $\delta$ is now a Dirac delta function and $T=\frac{t}{n}=\frac{2\pi}{\omega}$, the previous solution turns out to be very different.
Indeed, from Eq. (\ref{expAt}) we have

\begin{equation}
e^{-{\bf A}u}\,{\bf b}(u)=\frac{F_{0}}{m}
\sum_{k=1}^{n}\delta\left( u-\frac{2k\pi}{\omega}\right)
\left[ \begin{array}{c} 
-\frac{1}{\omega_0}\sin \omega_{0}u \\
 \cos \omega_{0}u \\
\end{array} \right]
\end{equation}

and, by computing the integral in Eq. (\ref{generalsolution6}), we get
\begin{equation}
\int_{0}^{t}e^{-{\bf A}u}{\bf b}(u)du=\frac{F_{0}}{m}
\left[ \begin{array}{c} 
-\frac{1}{\omega_0}\sum_{k=1}^{n}\sin\left(2\pi\frac{\omega_{0}}{\omega}k\right) \\
 \sum_{k=1}^{n}\cos\left(2\pi\frac{\omega_{0}}{\omega}k\right)\\
\end{array} \right]
\end{equation}

By the Dirichlet kernels
%
and replacing $2n+1=\frac{\pi+\omega t}{\pi}$, we finally obtain

%
\begin{equation}
{\bf y}(t)=\frac{F_{0}}{m \omega_{0}}
\frac{\sin \frac{\omega_{0}}{2}t}{\sin \frac{\omega_{0}}{\omega}\pi}
\left[ \begin{array}{c}
\sin \left(\frac{\omega_{0}}{2}t-\frac{\omega_{0}}{\omega}\pi\right)\\
\omega_{0}\cos \left(\frac{\omega_{0}}{2}t-\frac{\omega_{0}}{\omega}\pi\right)\\
\end{array} \right]
\label{eq69}
\end{equation}

Let us observe that the limit for $\omega\to \omega_{0}$ in Eq. (\ref{eq69}) shows a different behavior than in Eq. (\ref{synchronized_solution}): when $\omega$ approaches $\omega_{0}$, only the term $\sin \frac{\omega_{0}}{\omega}\pi$ in the denominator approaches $0$, while other terms represent an oscillation of frequency $\omega_{0}$. Therefore, the resonance produced by this kind of force preserves the proper frequency $\omega_{0}$ while the amplitude goes to infinity uniformly for all $t$.

\subsection{Network case}

Let us now move to the network case. The general solution of system (\ref{system4}) is given by
\begin{equation}
{\bf y}(t)=e^{{\bf G}t}\left[ \int_{0}^{t}e^{-{\bf G}u}{\bf b}(u)du +{\bf c} \right].
\label{generalsolution7}
\end{equation}

In the absence of dissipative terms, by setting ${c_1}={c_2}=1$, matrix $\bf G$ in Eq. (\ref{matrixG}) becomes
\begin{equation}
{\bf G}=
\left[ \begin{array}{cc} 
{\bf 0} & {\bf I} \\
-{\bf H}  & {\bf 0} \\
\end{array} \right]
\label{matrixG2}
\end{equation}

where now ${\bf H}={\bf I}+{\bf L}$. In Eq. (\ref{generalsolution3}), we proved that
\begin{equation}
e^{{\bf G}t}=\left[ \begin{array}{cc} 
\cos{\sqrt{\bf H}}t & (\sqrt{\bf H})^{-1}\sin\sqrt{\bf H}t \\
-\sqrt{\bf H}\sin{\sqrt{\bf H}}t & \cos\sqrt{\bf H}t \\
\end{array} \right]\\
\label{expGt}
\end{equation}

Let us assume a driving force $f(t)=F_{0}\sin(\omega t)$ applied on node $h$, $h=1,\dots, n$. The general solution and the corresponding resonant states are described by the next proposition.

\begin{proposition}
\label{proposition3}
The solution in Eq. (\ref{generalsolution7}) to system (\ref{system4}), with ${\bf y}(0)={\bf y}_{0}={\bf 0}$, for a sinusoidal driving force $f(t)=F_{0}\sin \omega t$ acting on node $h$, $h=1,\dots, n$, is
\begin{equation}
{\bf y}(t)=F_{0} \sum_{i=1}^{n}
\frac{\omega}{\omega_{i}^{2}-\omega^2}\phi_{i}(h)
\left[\begin{array}{c} 
\left(\frac{\sin \omega t}{\omega}- \frac{\sin \omega_{i} t}{\omega_{i}} \right) \phi_{i} \\
\left(\cos \omega t- \cos \omega_{i} t \right) \phi_{i}\\
\end{array} \right]\\
\label{synchronized_solution1}
\end{equation}
where $\omega_{i}\coloneqq\sqrt{1+\mu_{i}}$ and $\phi_{i}$ are the eigenvectors of the Laplacian matrix.
\end{proposition}

The reader can find the proof of proposition (\ref{proposition3}) in the Appendix \ref{Appendix B}. Here we want to focus on its consequences.  Eq. (\ref{synchronized_solution1}) can be expressed by components as
\begin{equation}
x_{k}(t)=F_{0}\omega \sum_{i=1}^{n} \frac{\phi_{i}(h)\phi_{i}(k)}{\omega^{2}-\omega_{i}^2}\left(\frac{\sin \omega_{i} t}{\omega_{i}}- \frac{\sin \omega t}{\omega} \right)
\label{position}
\end{equation}
\begin{equation}
v_{k}(t)=F_{0}\omega \sum_{i=1}^{n} \frac{\phi_{i}(h)\phi_{i}(k)}{\omega^{2}-\omega_{i}^2}\left(\cos \omega_{i} t - \cos \omega t \right)
\label{velocity}
\end{equation}

In particular,
\begin{equation}
\begin{split}
x_{k}(t)=&
F_{0}\omega \sum_{i=1}^{n-1} \frac{\phi_{i}(h)\phi_{i}(k)}{\omega^{2}-\omega_{i}^2}\left(\frac{\sin \omega_{i} t}{\omega_{i}}- \frac{\sin \omega t}{\omega} \right)\\
&+ \frac{F_{0}}{n}\frac{\omega}{\omega^{2}-1}\left( \sin t - \frac{\sin \omega t}{\omega}\right)
\label{separatesolution1}
\end{split}
\end{equation}

since $\mu_{1}>\mu_{2}>\dots>\mu_{n}=0$, $\omega_{1}>\omega_{2}>\dots>\omega_{n}=1$ and $\phi_{n}(k)=\frac{1}{\sqrt{n}}$, $\forall k$.

Let us analyze first what happens when $\omega \to \omega_{\tilde{\imath}}$ for some specific $\tilde{\imath}=1,\dots, n$, that is when the frequency of the external force approaches one of the \textit{proper} frequencies of the network. The diverging term in the position components in Eq. (\ref{position}) is
\begin{equation}
\begin{split}
x_{k}(t)&\sim F_{0}\omega \frac{\phi_{\tilde{\imath}}(h)\phi_{\tilde{\imath}}(k)}{\omega^{2}-\omega_{\tilde{\imath}}^2}\left(\frac{\sin \omega_{\tilde{\imath}} t}{\omega_{\tilde{\imath}}}- \frac{\sin \omega t}{\omega} \right)\\
&= F_{0} \frac{\phi_{\tilde{\imath}}(h)\phi_{\tilde{\imath}}(k){}}{\omega_{\tilde{\imath}}(\omega + \omega_{\tilde{\imath}})(\omega - \omega_{\tilde{\imath}})}\left[\omega \sin \omega_{\tilde{\imath}} t - \omega_{\tilde{\imath}} \sin \omega t \right]\\
\end{split}
\end{equation}

Since, in general, $\frac{x\sin tx_{0}-x_{0}\sin t x}{x-x_{0}}\sim \sin tx_{0}-tx_{0}\cos tx_{0}$ for $x\to x_{0}$, we have
\begin{equation}
\begin{split}
x_{k}(t)&\sim F_{0} \frac{\phi_{\tilde{\imath}}(h)\phi_{\tilde{\imath}}(k){}}{\omega_{\tilde{\imath}}(\omega + \omega_{\tilde{\imath}})}\left[ \sin \omega_{\tilde{\imath}} t - \omega_{\tilde{\imath}}t \cos \omega_{\tilde{\imath}} t \right]\\
&= F_{0} \frac{\phi_{\tilde{\imath}}(h)\phi_{\tilde{\imath}}(k){}}{2\omega_{\tilde{\imath}}^{2}}\left[ \sin \omega_{\tilde{\imath}} t - \omega_{\tilde{\imath}}t \cos \omega_{\tilde{\imath}} t \right]
\end{split}
\end{equation}

In a similar vein, the diverging term in the velocity components in Eq. (\ref{velocity}), for $\omega \to \omega_{\tilde{\imath}}$, is
\begin{equation}
\begin{split}
v_{k}(t)&\sim F_{0}\omega \frac{\phi_{\tilde{\imath}}(h)\phi_{\tilde{\imath}}(k)}{\omega^{2}-\omega_{\tilde{\imath}}^2}\left(\cos \omega_{\tilde{\imath}} t - \cos \omega t \right)\\
&= F_{0} \omega \phi_{\tilde{\imath}}(h)\phi_{\tilde{\imath}}(k) \frac{\sin \frac{\omega+\omega_{\tilde{\imath}}}{2}t}{\omega + \omega_{\tilde{\imath}}}\cdot \frac{\sin \frac{\omega-\omega_{\tilde{\imath}}}{2}t}{\frac{\omega-\omega_{\tilde{\imath}}}{2}} \\
&\sim \frac{F_{0}}{2} \phi_{\tilde{\imath}}(h)\phi_{\tilde{\imath}}(k)\, t\sin  \omega_{\tilde{\imath}} t
\end{split}
\end{equation}

The first conclusion is that the amplitudes of both $x_{k}(t)$ and $v_{k}(t)$, \textit{in general}, grow linearly with $t$ as $\omega$ approaches one of the proper frequencies of the network. However, things are slightly more involved than this, as discussed in the next remark. 

\begin{remark}
Let us imagine that, for $i=\tilde{\imath}$, the product $\phi_{\tilde{\imath}}(h)\phi_{\tilde{\imath}}(k)$ is equal to $0$.

Let us begin by examining the case in which $\phi_{\tilde{\imath}}(k)=0$. Then, even if $\omega \to \omega_{\tilde{\imath}}$, the node $k$ cannot resonate with the source node $h$ at such a frequency $\omega_{\tilde{\imath}}$. In other terms, if we make the source node $h$ oscillate with frequency $\omega=\omega_{\tilde{\imath}}$ and $\phi_{\tilde{\imath}}(k)=0$, then node $k$ is not affected by the oscillation in $h$. It is as if node $k$ is transparent to the propagation of the oscillation of frequency $\omega_{\tilde{\imath}}$ from the source node $h$. The network as a whole cannot transmit such a frequency to node $k$ and this node does not participate to the global resonant process.

On the other hand, if it is $\phi_{\tilde{\imath}}(h)=0$, then an oscillation of frequency $\omega = \omega_{\tilde{\imath}}$ in the source node $h$ cannot be spread throughout the network to any node and no resonance phenomenon can be induced from such a node with that frequency.

The only frequency that is able to always put in resonance the entire network is $\omega_{\tilde{\imath}}=\omega_{n}=1$. In fact, if we look at Eq. (\ref{separatesolution1}), being $\phi_{n}(h)\neq 0$ and $\phi_{n}(k)\neq 0$, $\forall h,k$, when $\omega \to 1$, any receiver node $k$ resonates with any source node $h$ oscillating at such a ground frequency.

More, specifically, let us observe that, for $\omega\to 1$, we have
\begin{equation}
\begin{split}
&x_{k}(t)\sim \\
&F_{0} \sum_{i=1}^{n-1} \frac{\phi_{i}(h)\phi_{i}(k)}{\mu_{i}}\left(\sin t - \frac{\sin \omega_{i} t}{\omega_{i}} \right)+
\frac{F_{0}}{2n}\left( \sin t - t\cos t\right)
\label{frequencyone}
\end{split}
\end{equation}

since $1-\omega_{i}^{2}=-\mu_{i}$. The quantity $
L^{+}_{hk}\coloneqq \sum_{i=1}^{n-1}\frac{\phi_{i}(h)\phi_{i}(k)}{\mu_{i}}$ represents the $hk$-component of the Moore-Penrose pseudo-inverse of the Laplacian operator. Sometimes it is also called vibrational communicability and denoted by $G^{\rm (v)}_{hk}\coloneqq L^{+}_{hk}$. It can be interpreted as the correlation function between the displacements $x_{h}(t)$ and $x_{k}(t)$. Formula (\ref{frequencyone}) contains a term that is exactly the product of the vibrational communicability between source and receiver nodes and the driving force $F_{0}\sin t$. This term tunes the extent of the interaction between source and receiver node when $\omega$ tends to the ground resonance frequency. More generally, the ratio
\begin{equation}
\frac{\phi_{i}(h)\phi_{i}(k)}{\omega^{2}-\omega_{i}^2}
\end{equation}

in Eqs. (\ref{position}) and (\ref{velocity}) contains the information about the local communicability between nodes and it conveys the role of the network topology in the resonant process from the source node $h$ to the receiver node $k$ at the ground frequency. 
\end{remark}

To exemplify, let us return to the initial toy network in the figure \ref{fig1}. Let us take a look at the explicit matrix of eigenvectors
\begin{equation}
{\bf M}=\left[
\begin{array}{cccc} 
0.289 & 0.707 & -0.408 & -0.500 \\
0.289 & -0.707 & -0.408 & -0.500 \\
-0.866 & 0 & 0 & -0.500 \\
0.289 & 0 & 0.816 & -0.500 \\
\end{array} \right]
\end{equation}

with eigenvalues $\mu_{1}=4$, $\mu_{2}=3$, $\mu_{3}=1$, $\mu_{4}=0$.
The eigen-frequencies are then $\omega_{1}=\sqrt{5}$, $\omega_{2}=2$, $\omega_{3}=\sqrt{2}$, $\omega_{4}=1$. Let us suppose that the driving force is applied in the first node; that is, $h=1$. Then the motion of the four nodes is described by the following functions:
\begin{equation*}
\begin{split}
x_{1}(t)
=&F_{0}\omega\Bigg[
\frac{0.084}{\omega^{2}-5}\left(\frac{\sin \sqrt{5} t}{\sqrt{5}}- \frac{\sin \omega t}{\omega} \right)\\
&+\frac{0.500}{\omega^{2}-4}\left(\frac{\sin 2t}{2}- \frac{\sin \omega t}{\omega} \right) \\
&+\frac{0.166}{\omega^{2}-2}\left(\frac{\sin \sqrt{2} t}{\sqrt{2}}- \frac{\sin \omega t}{\omega} \right)\\
&+\frac{0.250}{\omega^{2}-1}\left(\sin t- \frac{\sin \omega t}{\omega} \right)\Bigg]
\end{split}
\end{equation*}

\begin{equation*}
\begin{split}
x_{2}(t)
=&F_{0}\omega\Bigg[
\frac{0.084}{\omega^{2}-5}\left(\frac{\sin \sqrt{5} t}{\sqrt{5}}- \frac{\sin \omega t}{\omega} \right)\\
&-\frac{0.500}{\omega^{2}-4}\left(\frac{\sin 2t}{2}- \frac{\sin \omega t}{\omega} \right) \\
&+\frac{0.166}{\omega^{2}-2}\left(\frac{\sin \sqrt{2} t}{\sqrt{2}}- \frac{\sin \omega t}{\omega} \right)\\
&+\frac{0.250}{\omega^{2}-1}\left(\sin t- \frac{\sin \omega t}{\omega} \right)\Bigg]
\end{split}
\end{equation*}

\begin{equation*}
\begin{split}
x_{3}(t)
=&F_{0}\omega\Bigg[
-\frac{0.250}{\omega^{2}-5}\left(\frac{\sin \sqrt{5} t}{\sqrt{5}}- \frac{\sin \omega t}{\omega} \right)\\
&+\frac{0.250}{\omega^{2}-1}\left(\sin t- \frac{\sin \omega t}{\omega} \right)\Bigg]
\end{split}
\end{equation*}

\begin{equation*}
\begin{split}
x_{4}(t)
=&F_{0}\omega\Bigg[
\frac{0.084}{\omega^{2}-5}\left(\frac{\sin \sqrt{5} t}{\sqrt{5}}- \frac{\sin \omega t}{\omega} \right) \\
&-\frac{0.333}{\omega^{2}-2}\left(\frac{\sin \sqrt{2} t}{\sqrt{2}}- \frac{\sin \omega t}{\omega} \right)\\
&+\frac{0.250}{\omega^{2}-1}\left(\sin t- \frac{\sin \omega t}{\omega} \right)\Bigg]
\end{split}
\end{equation*}

\begin{figure}[H]
	\includegraphics[width=1\linewidth]{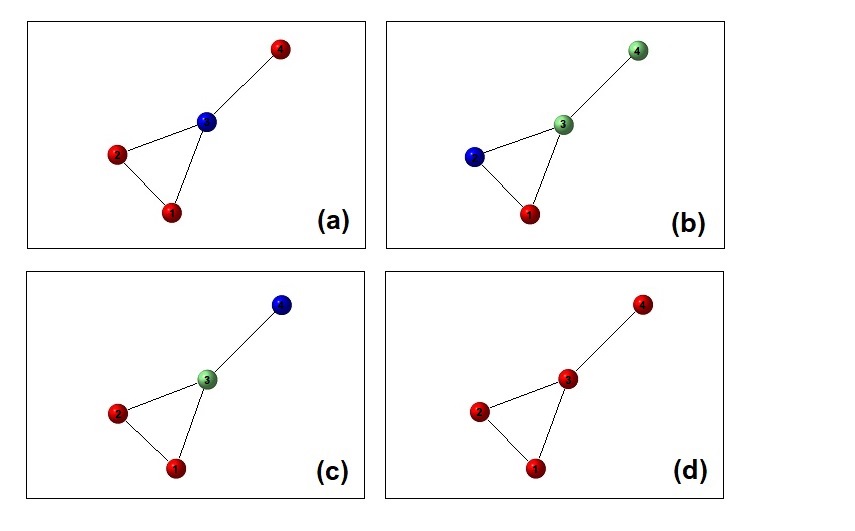}
	\centering
	\caption{Phase configuration of the whole network for different resonant driving frequencies: (a) $\omega_{1}=\sqrt{5}$; (b) $\omega_{2}=2$; (c) $\omega_{3}=\sqrt{2}$; (d) $\omega_{4}=1$. Red nodes have phase $+$, green nodes phase $0$ and blue nodes phase $-$.}
	\label{fig10} 
\end{figure}

The four possible configurations of the phases of nodes, for the different resonant driving frequencies applied in node $1$, are depicted in figure \ref{fig10}.

In particular, if the driving force in node number $1$ has frequency $\omega_{2}=2$ nodes number $3$ and $4$ do not resonate since $\phi_{2}(3)=\phi_{2}(4)=0$; similarly, when node number $1$ oscillates with frequency $\omega_{3}=\sqrt{2}$, only node number $3$ is prevented from entering the resonant state.

We also said that, if the driving force is placed in a node $h$ such that $\phi_{i}(h)=0$, then if $\omega \to \omega_{i}$, we have $x_{k}(t)=0,\ \forall k=1,\dots, n$. For instance, if the force is placed in the node $3$, being $\phi_{2}(3)=0$, this implies that for $\omega\to \omega_{2}=2$, all the nodes stay at rest. The network is not able to transmit such a frequency from such a node. Other frequencies may have effect but there will be no resonance with frequency $\omega_{2}=2$ placed in node number $3$.

\begin{figure}[H]
\centering
	\subfloat[]{\includegraphics[width=0.30\textwidth]{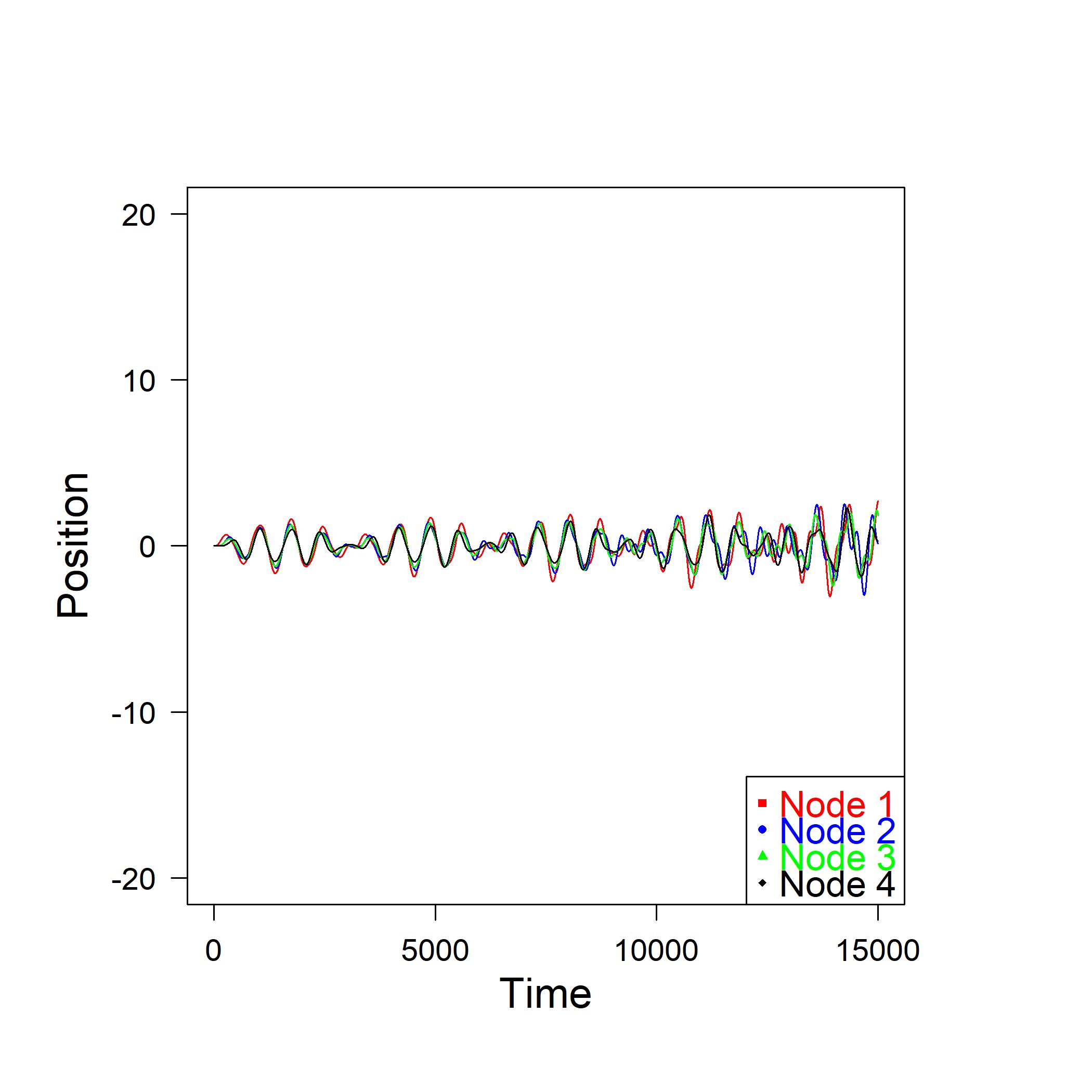}}
	\\
	\subfloat[]{\includegraphics[width=0.30\textwidth]{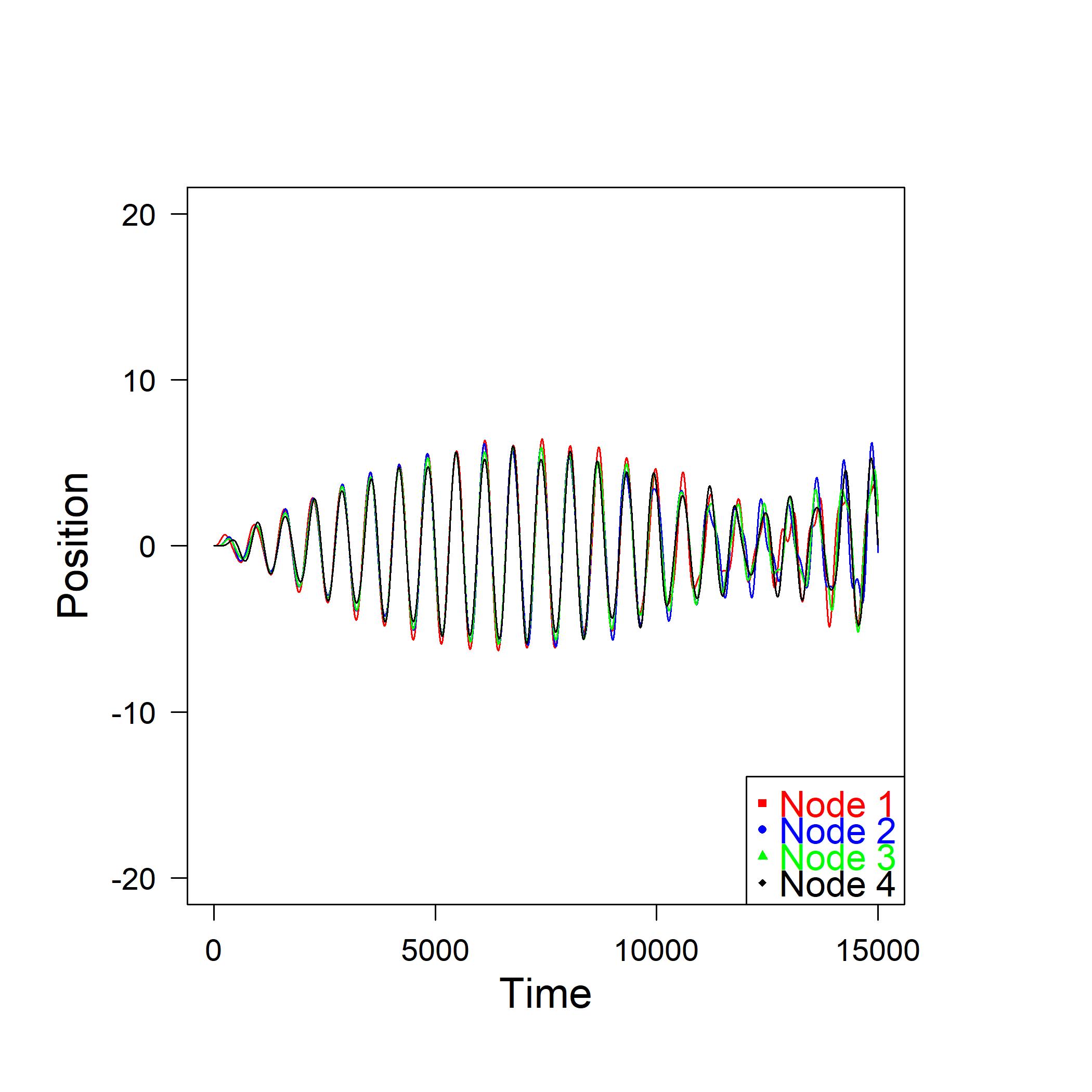}}
	\\
	\subfloat[]{\includegraphics[width=0.30\textwidth]{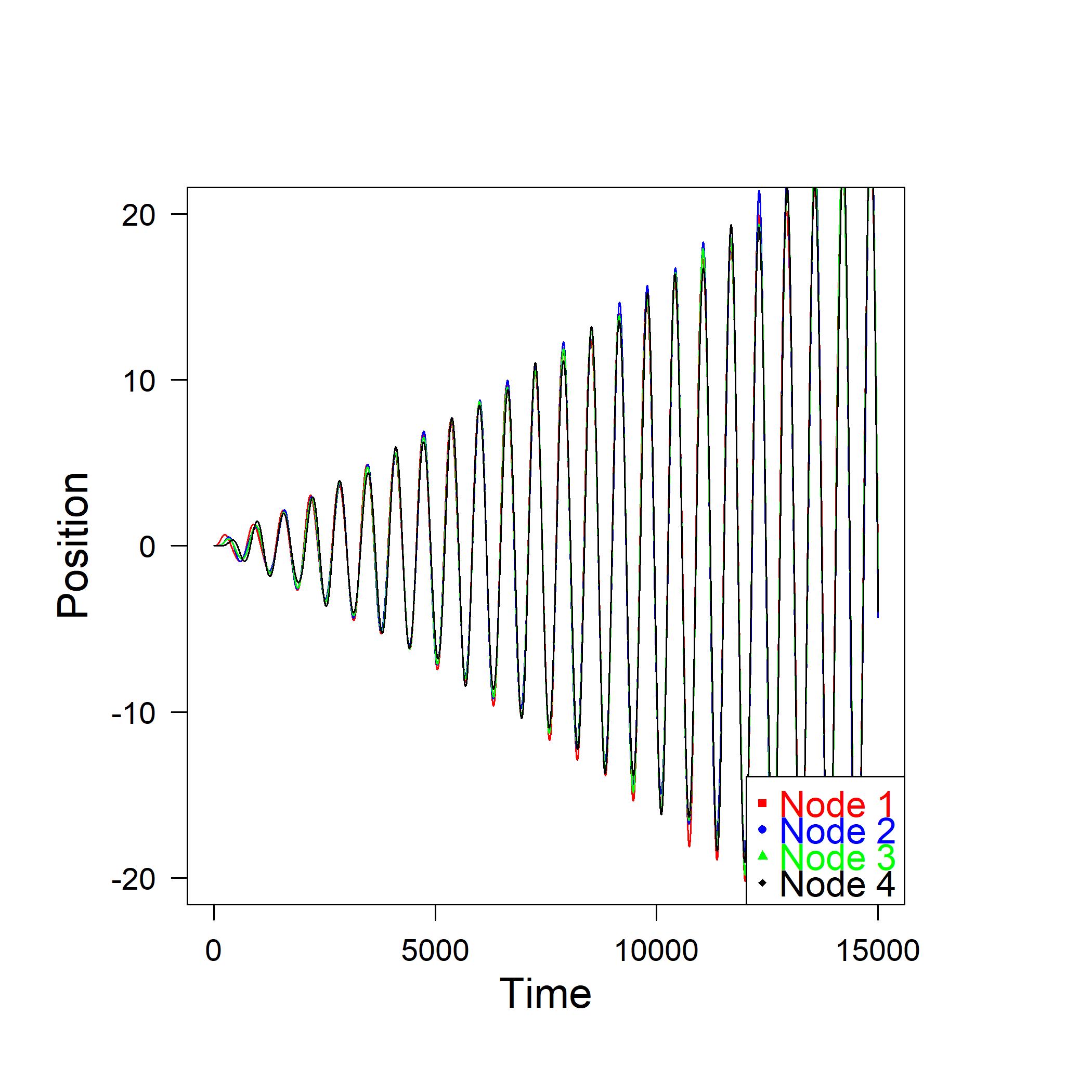}}
\caption{Approaching the resonant state with a driving force in node number $1$ of the toy network and ${\bf x}_{0}=[0,0,0,0]^T$ and ${\bf v}_{0}=[0,0,0,0]^T$; the frequencies are: panel (a) $\omega=0.80$; panel (b) $\omega=0.95$; and panel (c) $\omega=0.99$.}
	\label{fig11} 
\end{figure}

\begin{figure}[H]
\centering
	\subfloat[]{\includegraphics[width=0.30\textwidth]{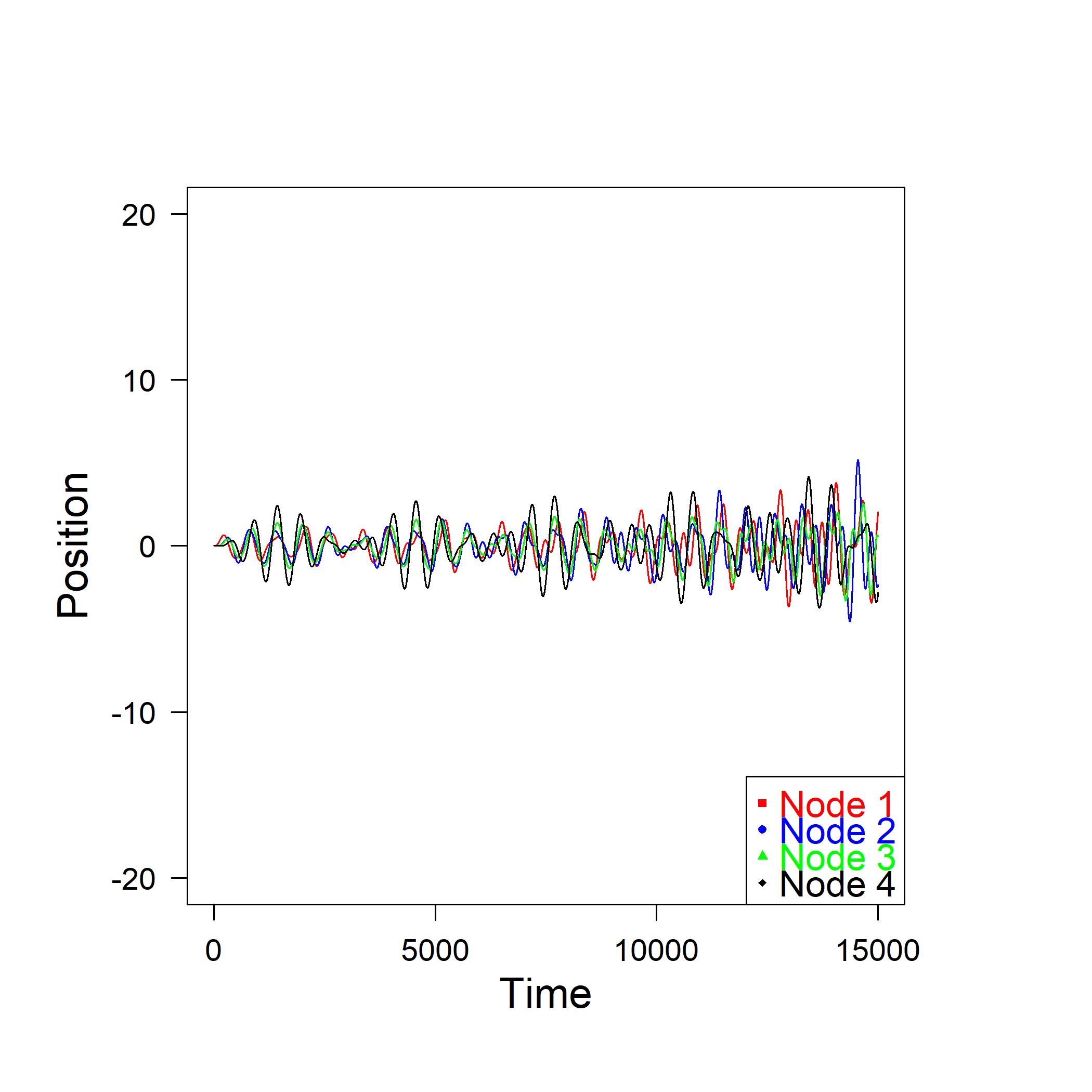}}
	\\
	\subfloat[]{\includegraphics[width=0.30\textwidth]{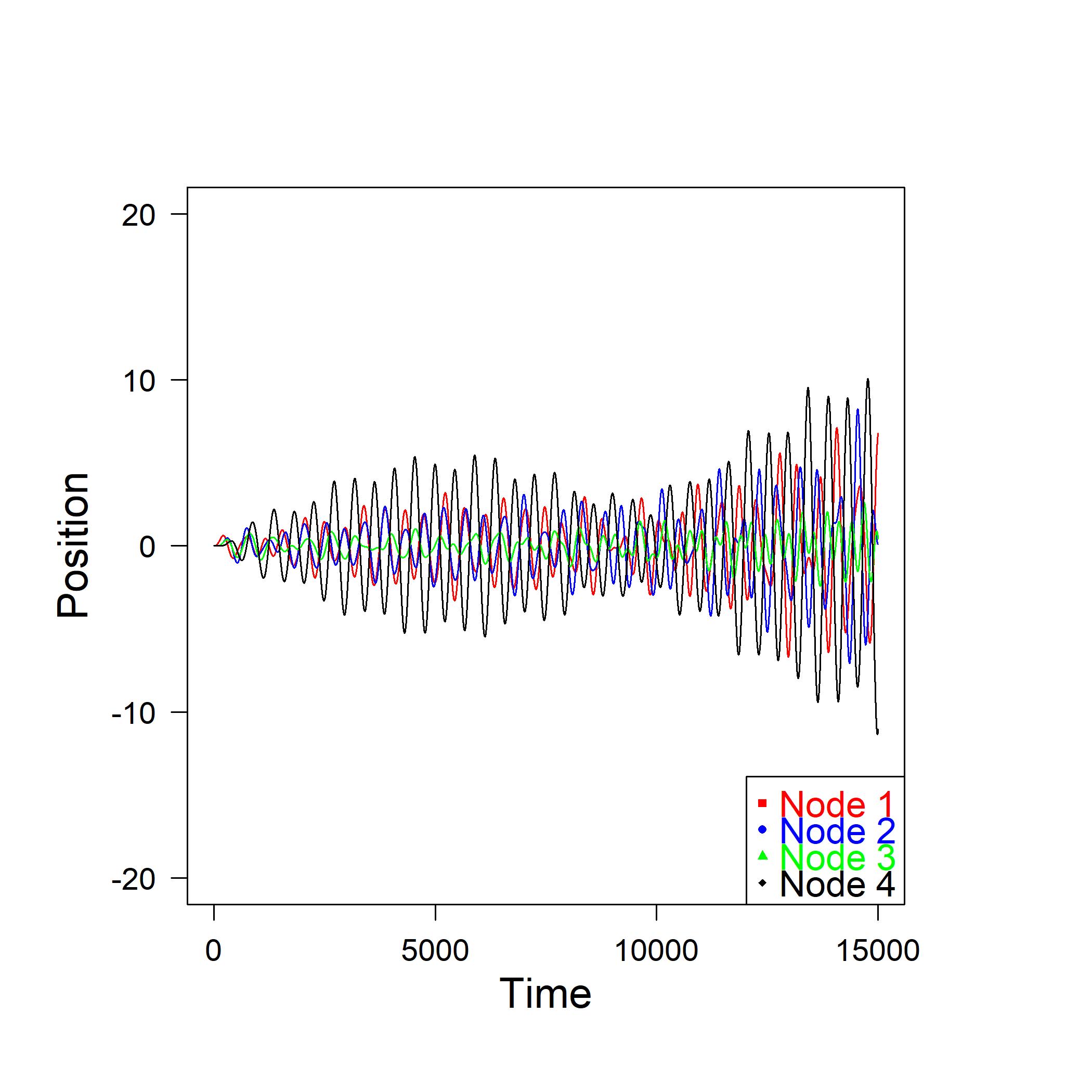}}
	\\
	\subfloat[]{\includegraphics[width=0.30\textwidth]{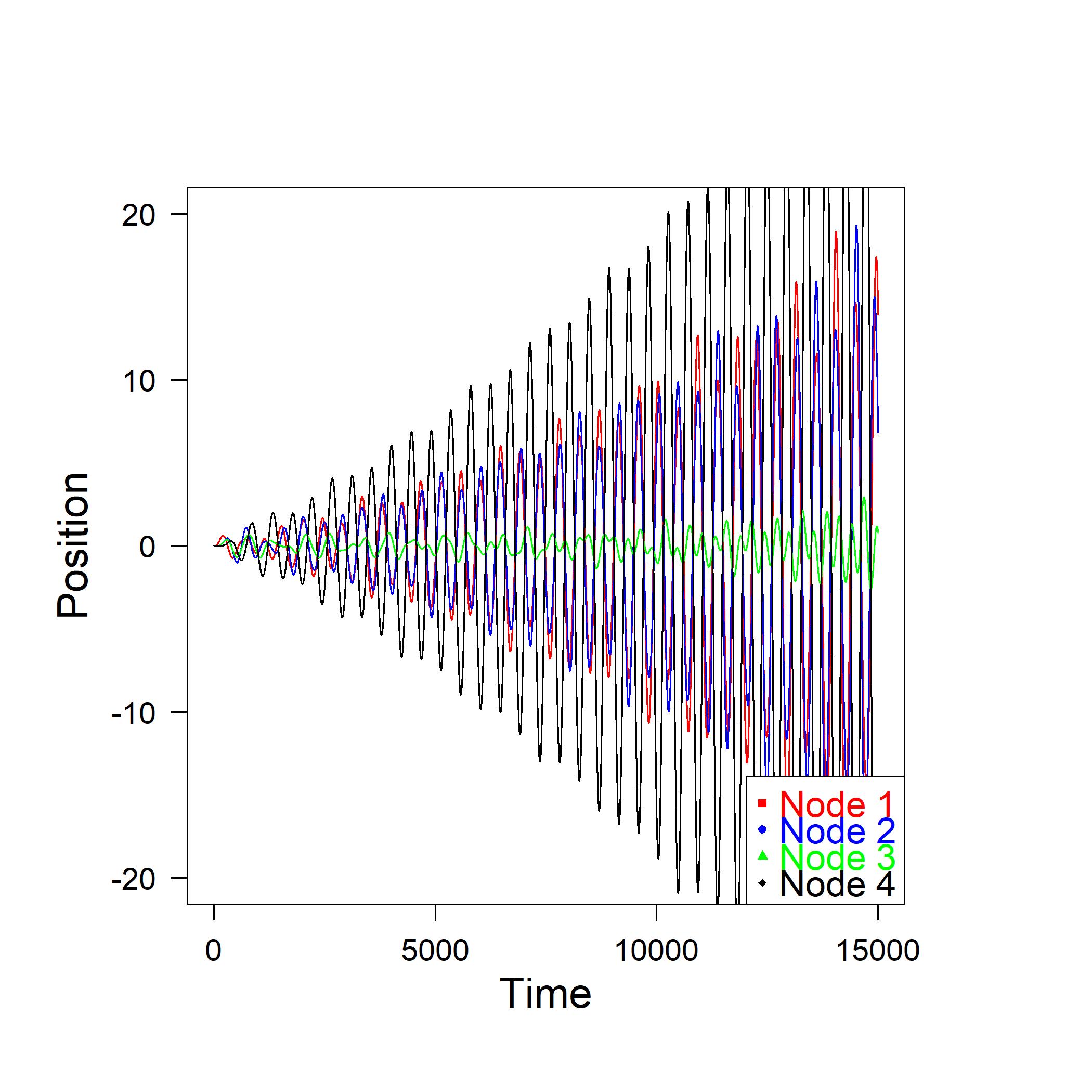}}
\caption{Approaching the resonant state with a driving force in node number $1$ of the toy network and ${\bf x}_{0}=[0,0,0,0]^T$ and ${\bf v}_{0}=[0,0,0,0]^T$; the frequencies are: panel (a) $\omega=1.20$; panel (b) $\omega=1.35$; and panel (c) $\omega=1.40$.}
	\label{fig12} 
\end{figure}

In figure \ref{fig11}, we illustrate the position plots of the four nodes in the toy network when an external force $f(t)=\sin (\omega t)$ is applied to node number $1$ starting from initial conditions ${\bf x}_{0}=[0,0,0,0]^T$ and ${\bf v}_{0}=[0,0,0,0]^T$. As $\omega$ approaches $\omega_{4}=1$ from below, all the nodes are activated and enter a global resonant state, as shown in figure \ref{fig9} panel (d).
In figure \ref{fig12}, we illustrate the position plots of the four nodes under the same conditions as before but as $\omega$ approaches $\omega_{3}=\sqrt{2}$ from below; let us notice that all nodes are activated, except node number $3$ (green line in the plot), which does not respond to node $1$ stimuli at this frequency, as shown in figure \ref{fig9} panel (c).

\subsection{Synchronization with resonance}

We conclude with some remarks on the case in which the coupling between nodes is due to dissipative terms in the presence of a periodic external force of the same type discussed in the previous subsection.
In this case, the fourth in Table \ref{summary}, matrix $\bf G$ in Eq. (\ref{matrixG}) takes the form
\begin{equation}
{\bf G}=
\left[ \begin{array}{cc} 
{\bf 0} & {\bf I} \\
-{\bf I}  & -{\bf L} \\
\end{array} \right]
\label{matrixG3}
\end{equation}

A straightforward computation shows that now
\begin{widetext}
\begin{equation}
\begin{split}
&\quad e^{{\bf G}t}=\\
&
\left[ \begin{array}{cc} 
\frac{1}{2}{\bf B}^{-1}\left[
({\bf B}+{\bf L})e^{\frac{{\bf B}-{\bf L}}{2}t}+
({\bf B}-{\bf L})e^{-\frac{{\bf B}+{\bf L}}{2}t}
\right]
&
{\bf B}^{-1}\left[e^{\frac{{\bf B}-{\bf L}}{2}t}-
e^{-\frac{{\bf B}+{\bf L}}{2}t}
\right]
\\
{\bf B}^{-1}\left[e^{-\frac{{\bf B}+{\bf L}}{2}t}-
e^{\frac{{\bf B}-{\bf L}}{2}t}
\right]
&
\frac{1}{2}{\bf B}^{-1}\left[
({\bf B}-{\bf L})e^{\frac{{\bf B}-{\bf L}}{2}t}+
({\bf B}+{\bf L})e^{-\frac{{\bf B}+{\bf L}}{2}t}
\right] \\
\end{array} \right]
\end{split}
\label{expGts}
\end{equation}
\end{widetext}

where ${\bf B}=\sqrt{{\bf L}^2-4{\bf I}}$.\footnote{Let us notice that, in order ${\bf B}^{-1}$ exists, it is necessary that ${\bf L}^2-4{\bf I}$ is non singular, that is $\mu_{i}\neq 2, \forall i=1,\dots, n$. However, from a computational point of view, it can be proven that this matrix exists for $\mu_{i}\to 2$.}
Let us now suppose again a driving force acting on node $h$ of the form $f(t)=F_{0}\sin \omega t$. The general solution and the corresponding resonant states are described by the next proposition.

\begin{proposition}
\label{proposition4}
The solution ${\bf y}(t)=[{\bf x}(t),{\bf v}(t)]^{T}=e^{{\bf G}t}\left[ \int_{0}^{t}e^{-{\bf G}u}{\bf b}(u)du +{\bf c} \right]$ with ${\bf b}(u)=f(u)\left[ {\bf 0}, {\bf e}_{h}  \right]^{T}$, ${\bf y}(0)={\bf y}_{0}={\bf 0}$ and $\bf G$ as in Eq. (\ref{matrixG3}) for a sinusoidal driving force $f(t)=F_{0}\sin \omega t$ acting on node $h$, $h=1,\dots, n$, is given by the following components:
\begin{equation*}
\begin{split}
&{\bf x}(t)=\\
&\sum_{i=1}^{n} \frac{F_{0}\phi_{i}^{\star}(h)}{\xi_{i}} 
\Bigg[ \left( 
\frac{1-(\lambda_{i}^{-})^2}{\xi_{i}(\omega^2+(\lambda_{i}^{-})^2)}+
\frac{1-(\lambda_{i}^{+})^2}{\xi_{i}(\omega^2+(\lambda_{i}^{+})^2}
\right) \sin \omega t \\
&+ \left( 
\frac{1}{\omega^2+(\lambda_{i}^{-})^2}-
\frac{1}{\omega^2+(\lambda_{i}^{+})^2}
\right) \omega \cos \omega t \\
&+\left( 
\frac{\omega}{\omega^2+(\lambda_{i}^{+})^2}e^{\lambda_{i}^{+} t}-\frac{\omega}{\omega^2+(\lambda_{i}^{-})^2}e^{\lambda_{i}^{-} t}
\right)\Bigg]\phi_{i}\\
\end{split}
\end{equation*}

and 

\begin{equation*}
\begin{split}
&{\bf v}(t)=\\
&\sum_{i=1}^{n} \frac{F_{0}\phi_{i}^{\star}(h)}{\xi_{i}} 
\Bigg[ \left( 
\frac{\lambda_{i}^{-}(1-(\lambda_{i}^{-})^2)}{\xi_{i}(\omega^2+(\lambda_{i}^{-})^2)}+
\frac{\lambda_{i}^{+}(1-(\lambda_{i}^{+})^2)}{\xi_{i}(\omega^2+(\lambda_{i}^{+})^2}
\right) \sin \omega t \\
&+ \left( 
\frac{\lambda_{i}^{-}}{\omega^2+(\lambda_{i}^{-})^2}-
\frac{\lambda_{i}^{+}}{\omega^2+(\lambda_{i}^{+})^2}
\right) \omega \cos \omega t \\
&+
\left( 
\frac{\omega\lambda_{i}^{-}}{\omega^2+(\lambda_{i}^{+})^2}e^{\lambda_{i}^{+} t}-\frac{\omega\lambda_{i}^{+}}{\omega^2+(\lambda_{i}^{-})^2}e^{\lambda_{i}^{-} t}
\right)\Bigg]\phi_{i}
\end{split}
\end{equation*}

where $\xi_{i}=\sqrt{\mu_{i}^{2}-4}$, $\lambda_{i}^{+}=\frac{1}{2}(\xi_{i}-\mu_{i})$, and $\lambda_{i}^{-}=-\frac{1}{2}(\xi_{i}+\mu_{i})$ as in Eq. (\ref{eigenvalues_G}).
\end{proposition}

Notice that, with respect to Eq. (\ref{eigenvalues_G}), we dropped the index G in the eigenvalues to make the expressions lighter.

Leaving the proof of this proposition to Appendix \ref{Appendix B}, let us now focus on the $x$ component and its behavior for $t$ large. After some further steps, it can be written as

\begin{equation*}
\begin{split}
&{\bf x}(t)=\\
&\sum_{i=1}^{n} \frac{F_{0}\phi_{i}^{\star}(h)}{(\omega^2+(\lambda_{i}^{-})^2)\cdot (\omega^2+(\lambda_{i}^{+})^2)} 
\Bigg[ \left(1-\omega^2
\right) \sin \omega t \\
& -\mu_{i} \omega \cos \omega t +
\frac{\omega}{\xi_{i}}
\left( 
{(\omega^2+(\lambda_{i}^{-})^2)}e^{\lambda_{i}^{+} t}-
{(\omega^2+(\lambda_{i}^{+})^2)}e^{\lambda_{i}^{-} t}
\right)\Bigg]\phi_{i}\\
\end{split}
\end{equation*}

and, by components,
\begin{equation}
\begin{split}
&{x}_{k}(t)=\\
&\sum_{i=1}^{n-1} \frac{F_{0}\phi_{i}(h)\phi_{i}(k)}{(\omega^2+(\lambda_{i}^{-})^2)\cdot (\omega^2+(\lambda_{i}^{+})^2)} 
\Bigg[ \left(1-\omega^2
\right) \sin \omega t -\\
& \mu_{i} \omega \cos \omega t +
\frac{\omega}{\xi_{i}}
\left( 
{(\omega^2+(\lambda_{i}^{-})^2)}e^{\lambda_{i}^{+} t}-
{(\omega^2+(\lambda_{i}^{+})^2)}e^{\lambda_{i}^{-} t}
\right)\Bigg]\\
&+ \frac{F_{0}}{n}\frac{\omega \sin t-\sin \omega t}{\omega^2-1}\\
\end{split}
\label{separatesolution2}
\end{equation}
formula in which we made the final term for $i=n$ explicit.

The question now is what frequency or frequencies are able to induce resonance phenomena in the network as $t$ grows to infinity, bearing in mind the fact that the presence of dissipative terms still drives synchronization between nodes. This amounts to asking what are the differences between formula (\ref{separatesolution1}) and formula (\ref{separatesolution2}).

First, let us note that the ratio before the square bracket in formula (\ref{separatesolution2}) is always a real value. This is due to the fact that $\lambda_{i}^{-}$ and $\lambda_{i}^{+}$ are reciprocal and to the equality $(\lambda_{i}^{-})^2+(\lambda_{i}^{+})^2=\mu_{i}^2-2$. Moreover, as previously discussed, the exponential terms $e^{\lambda_{i}^{\pm} t}$ inside the brackets go to zero, as $t$ grows. Therefore, in general, for $t\to +\infty$ we have coefficients in the sum depending on $k$ which prevents nodes from synchronizing. 

However, when $\omega\to 1$, the last term in the sum becomes dominant and this term is common to all the nodes. Indeed, the only frequency $\omega$ that can prompt resonance phenomena is the one such that $\omega^2+(\lambda_{i}^{-})^2=0$ and $\omega^2+(\lambda_{i}^{+})^2=0$ and the only value that can make these quantities vanishing is the pure imaginary value $\lambda_{n}^{+}=\lambda_{n}^{-}=i$. 

We conclude that, under these particular conditions, the only resonant state is exactly the one involving complete synchronization between the nodes. Therefore, in general, a driving force does not favor the synchronization process and, vice versa, a resonant state is, in general, prevented by the presence of this kind of coupling for frequencies other than $\omega=1$. There is, thus, only one fully synchronized and resonant state for $\omega \to 1$ and $t \to +\infty$. This state is described by the following asymptotic behavior, still linearly growing in time:
\begin{equation*}
{x}_{k}(t)\sim
\frac{F_{0}}{n}\frac{\omega \sin t-\sin \omega t}{\omega^2-1}
\sim \frac{F_{0}}{2n}(\sin t -t\cos t)
\end{equation*}

\begin{figure}[H]
\centering
	\subfloat[]{\includegraphics[width=0.30\textwidth]{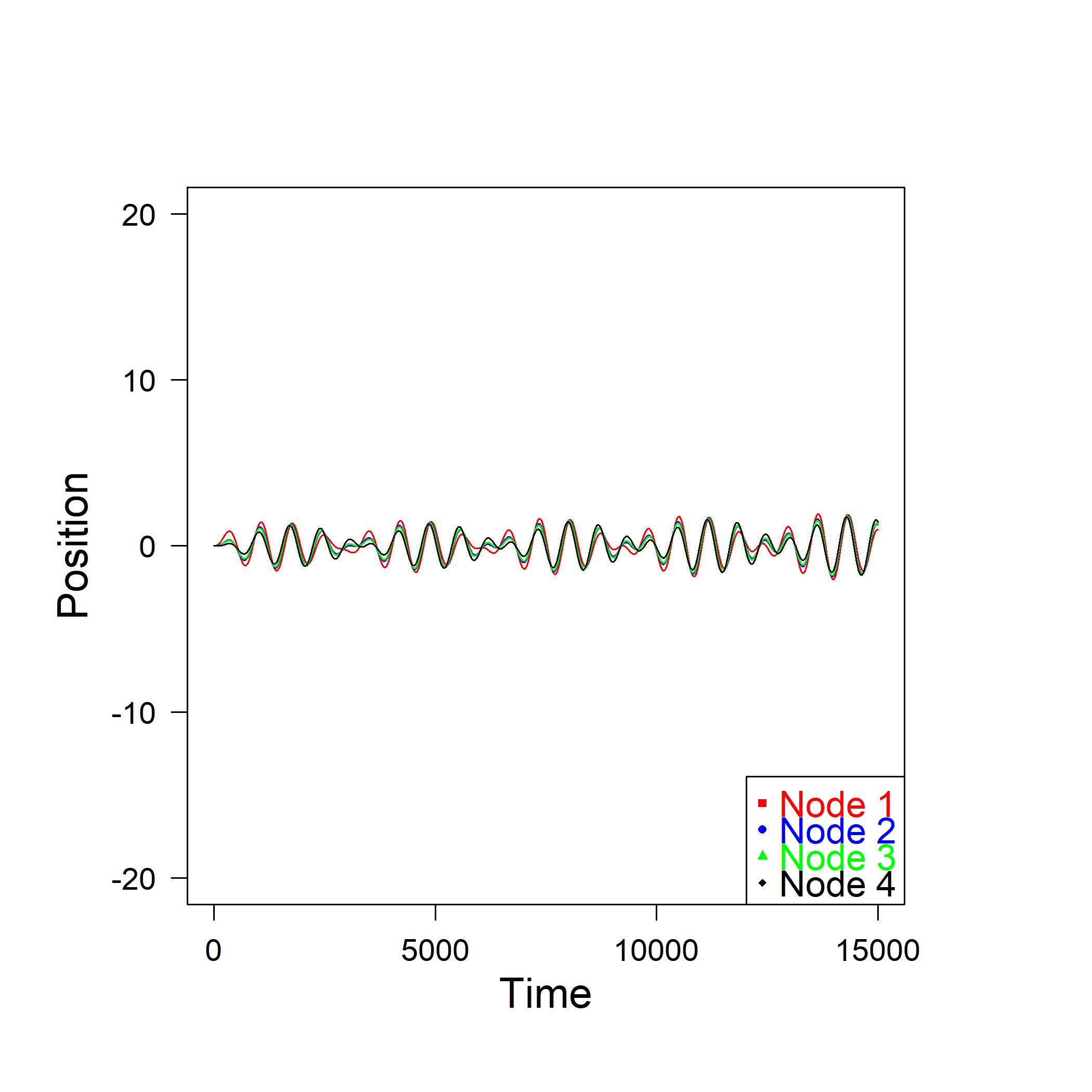}}
	\\
	\subfloat[]{\includegraphics[width=0.30\textwidth]{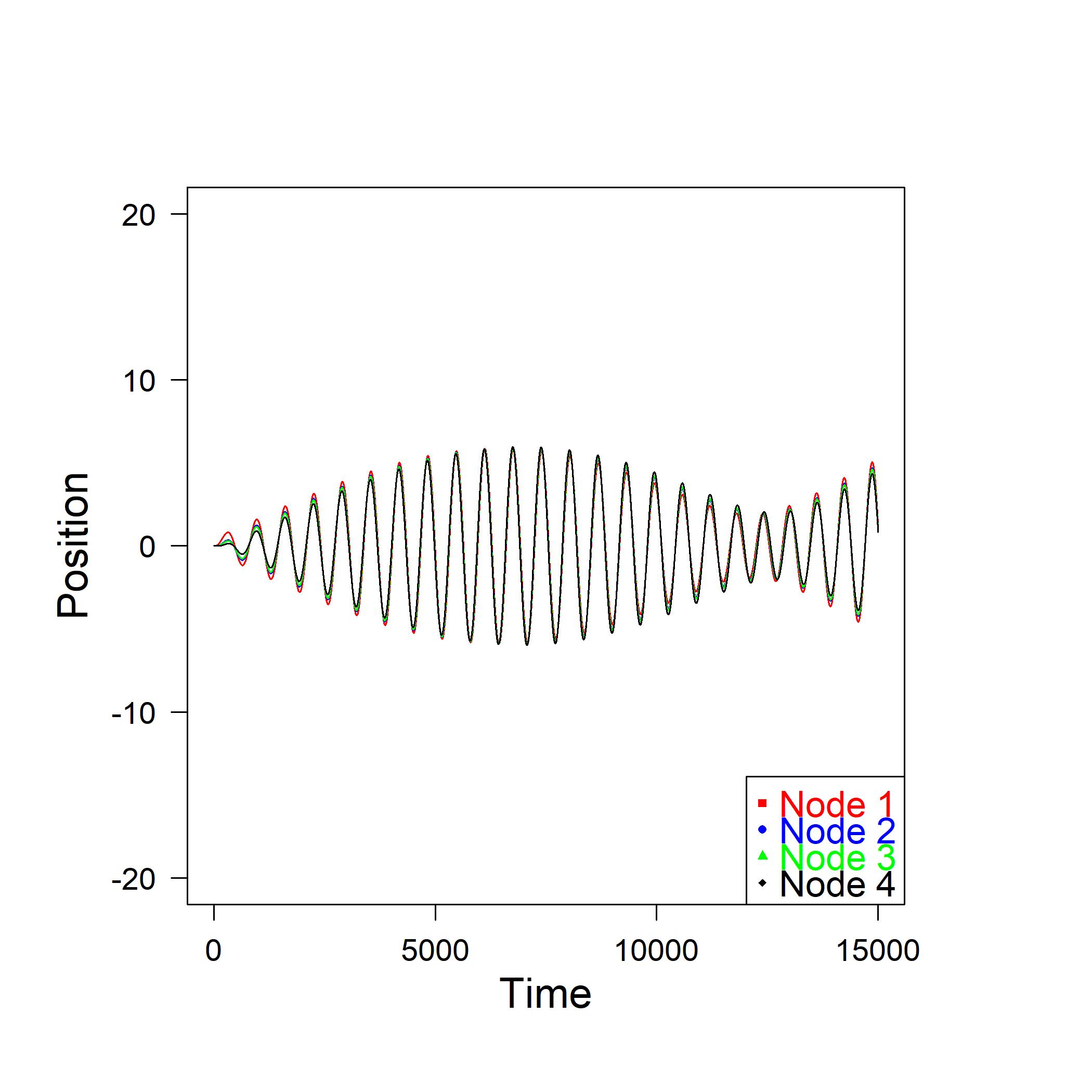}}
	\\
	\subfloat[]{\includegraphics[width=0.30\textwidth]{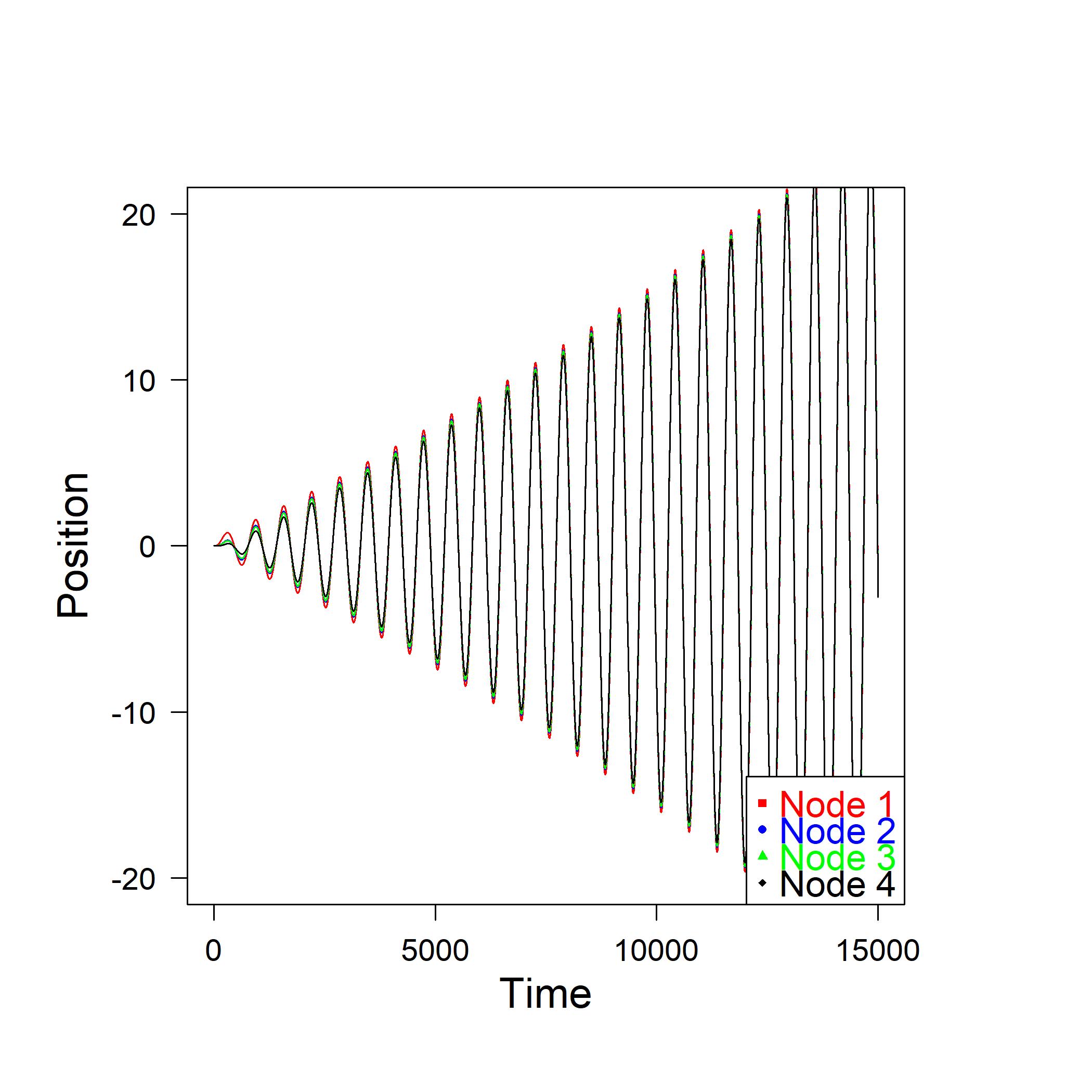}}
\caption{Approaching the resonant state with a driving force in node number $1$ of the toy network with dissipative couplings and ${\bf x}_{0}=[0,0,0,0]^T$ and ${\bf v}_{0}=[0,0,0,0]^T$; the frequencies are: panel (a) $\omega=0.80$; panel (b) $\omega=0.95$; and, panel (c) $\omega=0.99$.}
	\label{fig13}  
\end{figure}

\begin{figure}[H]
\centering
	\subfloat[]{\includegraphics[width=0.30\textwidth]{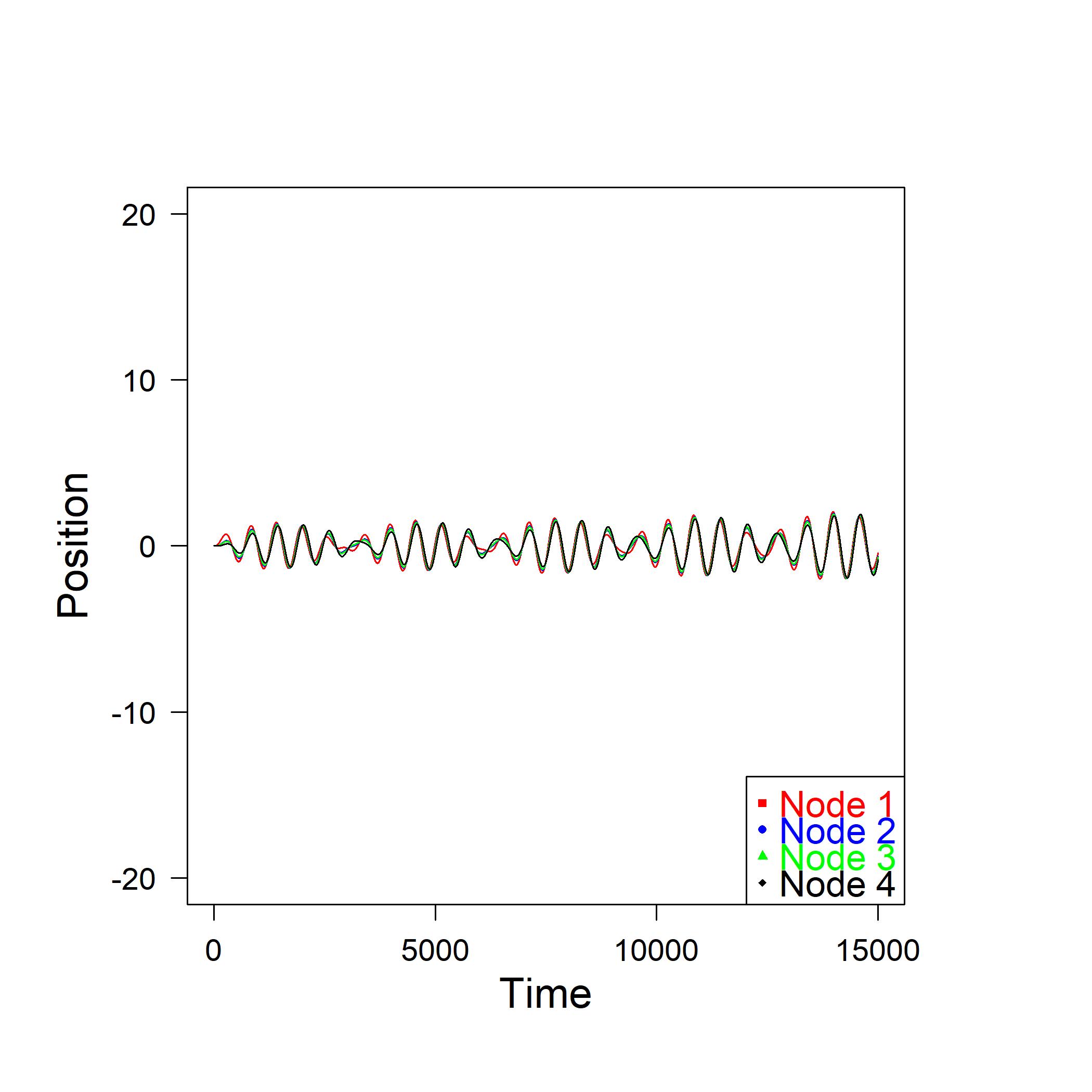}}
	\\
	\subfloat[]{\includegraphics[width=0.30\textwidth]{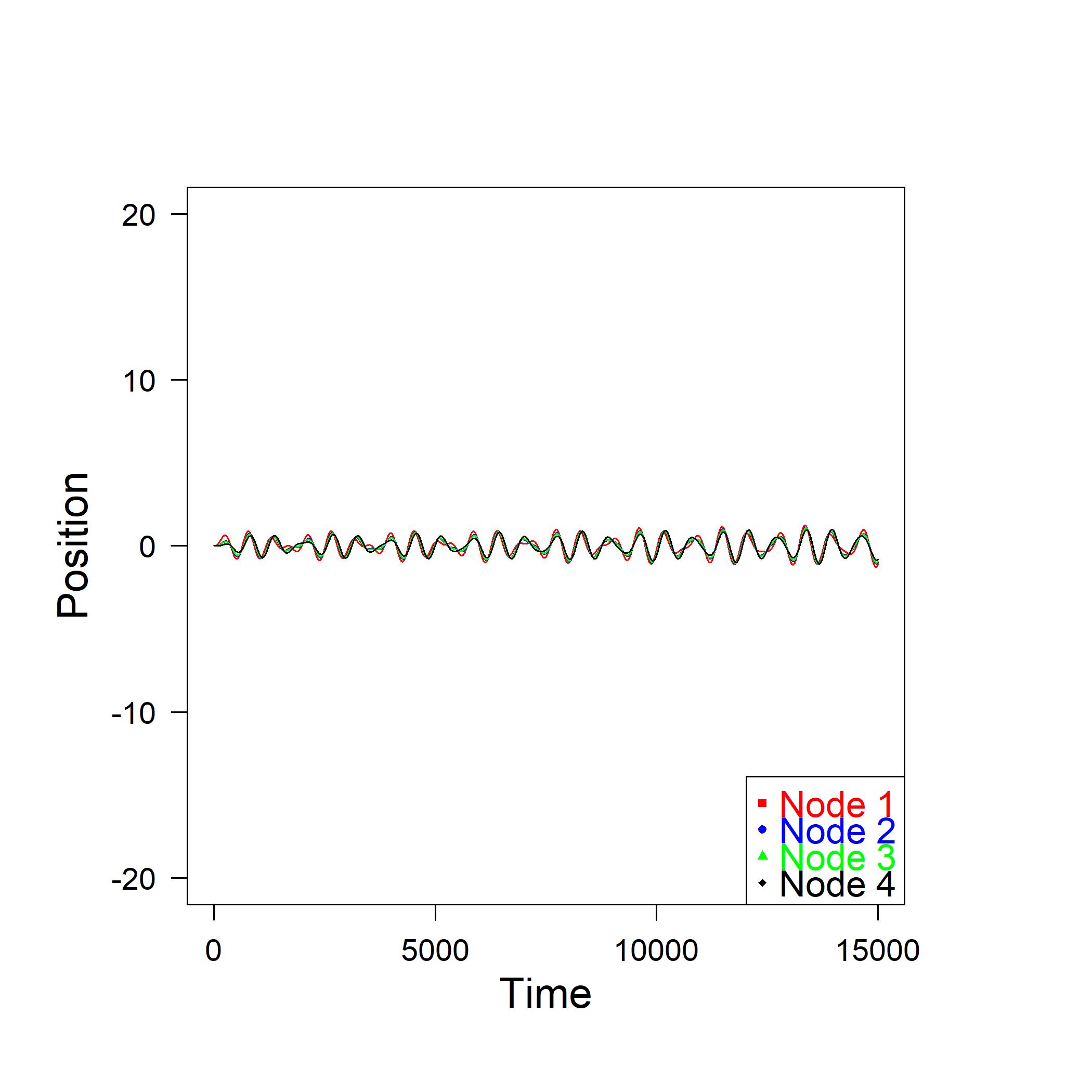}}
	\\
	\subfloat[]{\includegraphics[width=0.30\textwidth]{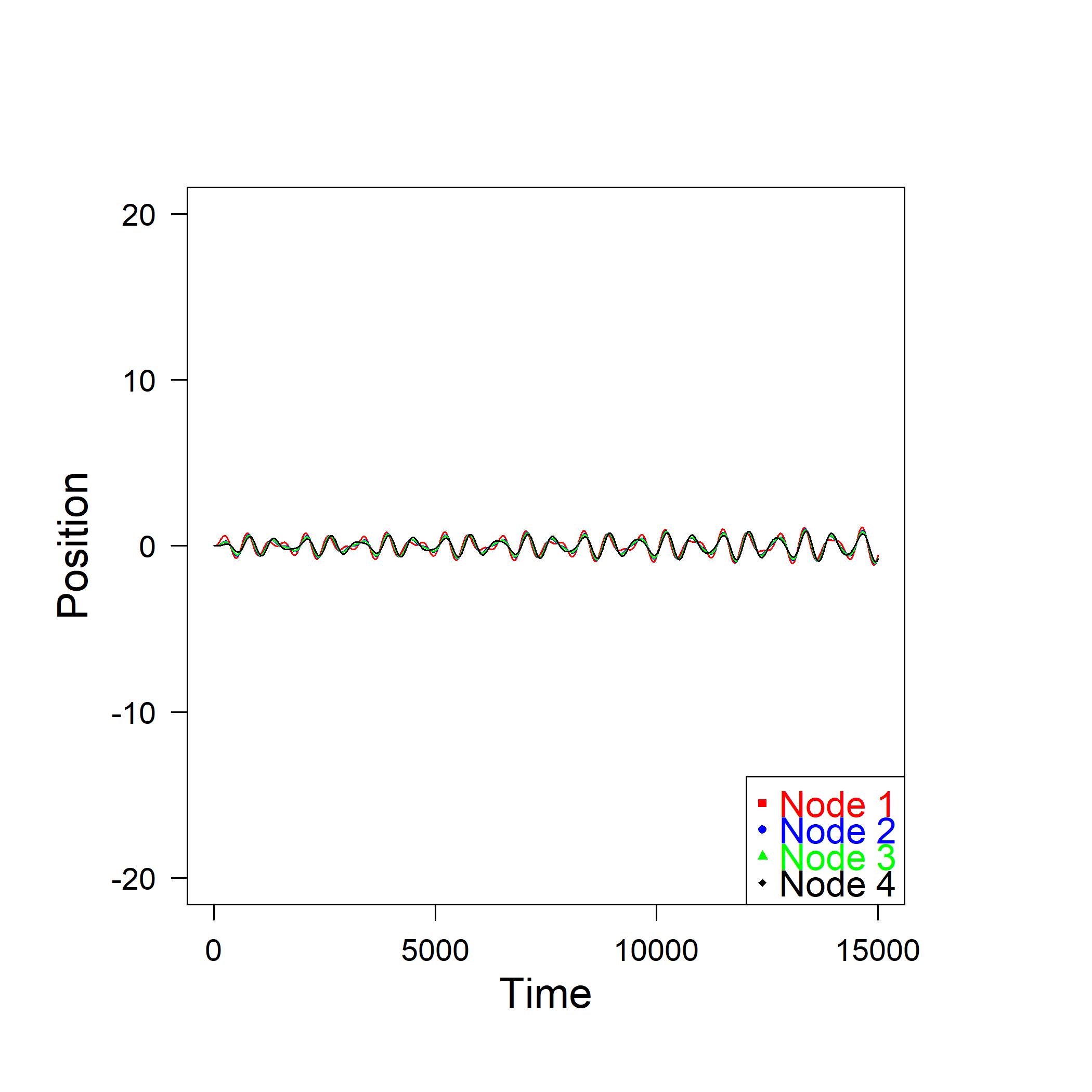}}
\caption{Approaching the non-resonant state with a driving force in node number $1$ of the toy network with dissipative couplings and ${\bf x}_{0}=[0,0,0,0]^T$ and ${\bf v}_{0}=[0,0,0,0]^T$; the frequencies are: panel (a) $\omega=1.20$; panel (b) $\omega=1.35$; and, panel (c) $\omega=1.40$.}
	\label{fig14}  
\end{figure}

In figure \ref{fig13}, we illustrate the position plots of the four nodes in the toy network when an external force $f(t)=\sin \omega t$ is applied to node number $1$ starting from initial conditions ${\bf x}_{0}=[0,0,0,0]^T$ and ${\bf v}_{0}=[0,0,0,0]^T$. As $\omega$ approaches $1$ from below, all the nodes are activated and enter a global resonant state, as proven above.
In figure \ref{fig14}, we illustrate the position plots of the four nodes under the same conditions as before but as $\omega$ approaches, for instance, $\sqrt{2}$ from below; let us notice that, as expected, no resonant state arises within the network, in contrast with what shown in figure \ref{fig12} under the same conditions.

%
%
%
%
%

\section{The Linear Swing Equation in Networks}
\label{section6}

Let us now return to system (\ref{system4}) and consider the case in which we set $c_{1}=0$, $c_{2}=1$, $c^{\prime}_{1}=\gamma \in {\mathbb R}^{+}$, $c^{\prime}_{2}=0$ in the matrix $\bf G$. Instead of considering a driving force $f(t)$ acting on a single node as before, we now give the factor ${\bf b}(t)$ the form of a constant vector. Therefore, let us consider the system 
\begin{equation}
	\left\{ 
	\begin{array}{l}
		{\bf \dot{y}}={\bf G}{\bf y}+{\bf b} \\
		{\bf y}(0)={\bf y}_{0} \\
	\end{array}
	\right.
	\label{system6}
\end{equation}
${\bf G}$ being explicitly the matrix
\begin{equation}
	{\bf G}\coloneqq
	\left[ \begin{array}{cc} 
		{\bf 0} & {\bf I} \\
		-{\bf L}  & -\gamma{\bf I} \\
	\end{array} \right]
	\label{matrixG4}
\end{equation}
and ${\bf b}(t)=[{\bf 0}, {\bf p} ]^{T}$, where ${\bf p}$ is a constant vector in ${\mathbb R}^{n}$. Now, system (\ref{system6}) is equivalent to the equation
\begin{equation}
	{\bf \ddot{x}}=-\gamma 	{\bf \dot{x}}-{\bf L}{\bf{x}}+{\bf{p}}
	\label{swing_equation}
\end{equation}
which is called the linearized swing equation \cite{Sorrentino2021, Sorrentino2022}. The swing equation plays a central role in the analysis of power systems dynamics and in a variety of different contexts, such as the description of superconducting Josephson junctions or nonlinear second-order phase-locking loop in electronics. It governs the rotational dynamics of synchronous machines in stability studies, it is nonlinear in nature and a closed formal solution cannot be explicitly given. Equation (\ref{swing_equation}), in particular, models the dynamics of power grid networks, where nodes can be either generators (generating power) or loads (consuming power) and edges are transmission lines connecting them. The requirement for a synchronous generator is that it does not lose phase synchronization when small temporary changes occur in the output of the machine. Typical examples are small perturbations in the scheduled generation of the machine, which result in a change in its rotor angle, or small loads added to the network. The response of a power system to impacts is oscillatory, and if the oscillations are damped, the system can reach in time a steady state. Transient stability studies determine whether or not synchronization is maintained after the machine has been subjected to a transient disturbance. Each variable $x_i$ in Eq. (\ref{swing_equation}) represents the corresponding node displacement and $v_i=\dot{x}_{i}$ its nodal velocity relative to a rotational reference value; $\gamma$, which we assume common to all the nodes, is the damping coefficient and the vector $\bf p$ has components equal (or, more in general, proportional) to the small power deviation generated (if positive) or consumed (if negative) by nodes. The general solution of Eq. (\ref{system6}) can be written as
\begin{equation}
{\bf y}(t)={\bf \tilde{y}}+e^{{\bf G}t}({\bf y}_{0}-{\bf \tilde{y}})
\label{solution_swing_equation}
\end{equation}
where ${\bf y}(0)={\bf y}_{0}$ (and typically, in this context, ${\bf y}_{0}={\bf 0}$) and ${\bf \tilde{y}}$ is such that ${\bf b}=-{\bf G}{\bf \tilde{y}}$ or, equivalently, such that ${\bf \tilde{y}}=[{\bf \tilde{x}}, {\bf \tilde{v}}]^{T}$ with ${\bf p}={\bf L}{\bf \tilde{x}}$ and ${\bf \tilde{v}}={\bf 0}$.

We now limit the discussion to balanced power grids, that is we assume that $\sum_{i=1}^{n}p_{i}=0$. Furthermore, we suppose that the power grid is connected so that the multiplicity of the null eigenvalue of $\bf L$ is equal to $1$. Under these conditions, the system ${\bf p}={\bf L}{\bf \tilde{x}}$ is consistent since ${\rm rank}({\bf L}|{\bf p})={\rm rank}({\bf L})$. The family of solutions is characterized by the fact that, if ${\bf \tilde{x}}_{1}$ and ${\bf \tilde{x}}_{2}$ are two of them, then ${\bf \tilde{x}}_{1}- {\bf \tilde{x}}_{2}=w{\bf u}$, where $w\in {\mathbb R}$ and ${\bf u}\in {\mathbb R}^{n}$ is the vector of all $1$'s. Since ${\bf G}[{\bf u},{\bf 0}]^{T}=[{\bf 0},{\bf 0}]^{T}$ and $e^{{\bf G}t}[{\bf u},{\bf 0}]^{T}=[{\bf u},{\bf 0}]^{T}$, then the general solution ${\bf y}(t)=e^{{\bf G}t}{\bf y}_{0}+{\bf \tilde{y}}-e^{{\bf G}t}{\bf \tilde{y}}$ is the same for any $w\in {\mathbb R}$. This means that we can use any solution of the linear system ${\bf p}={\bf L}{\bf \tilde{x}}$ to build the general solution in Eq. (\ref{solution_swing_equation}).

\begin{figure}[H]
	\centering
	\subfloat[]{\includegraphics[width=0.25\textwidth]{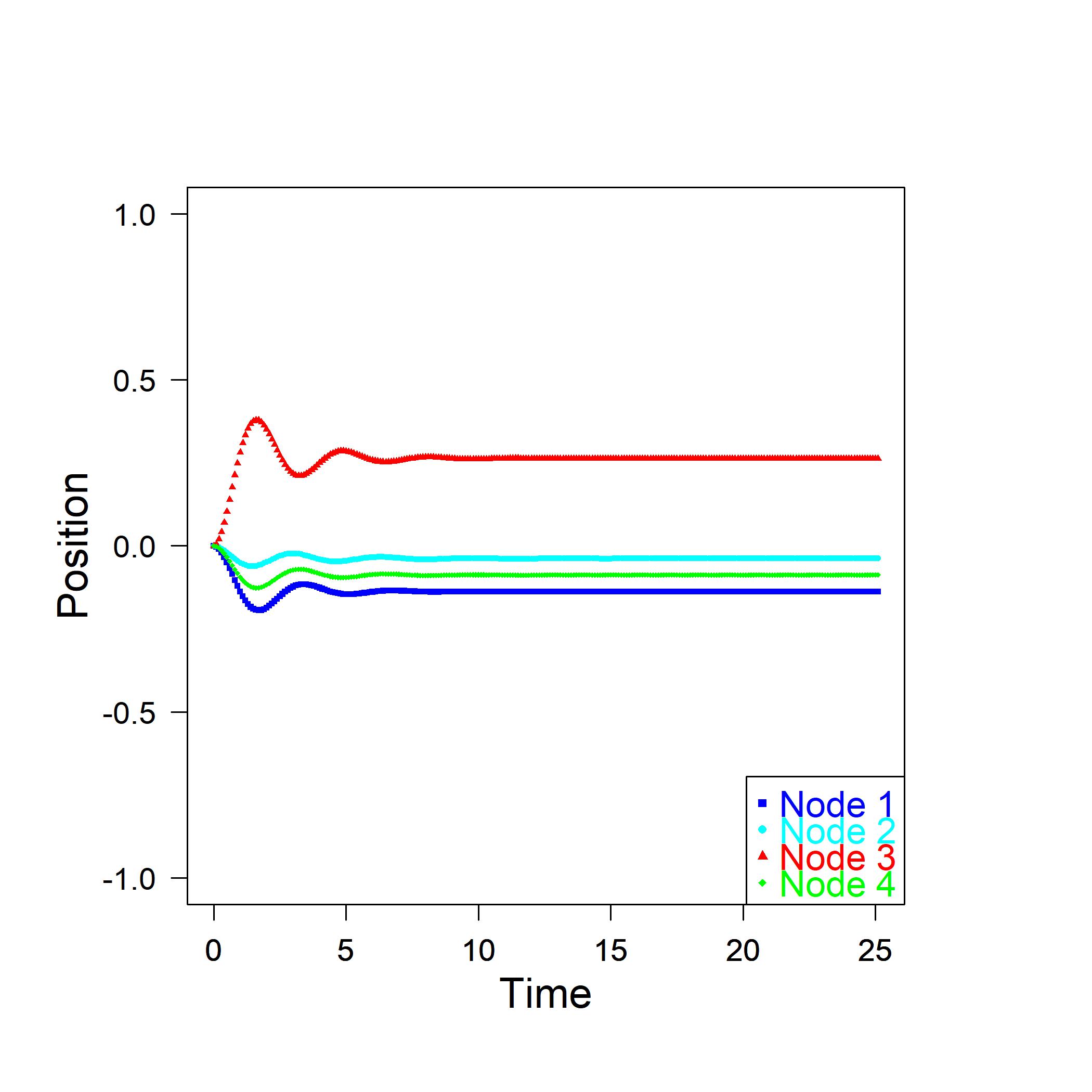}}
	\subfloat[]{\includegraphics[width=0.25\textwidth]{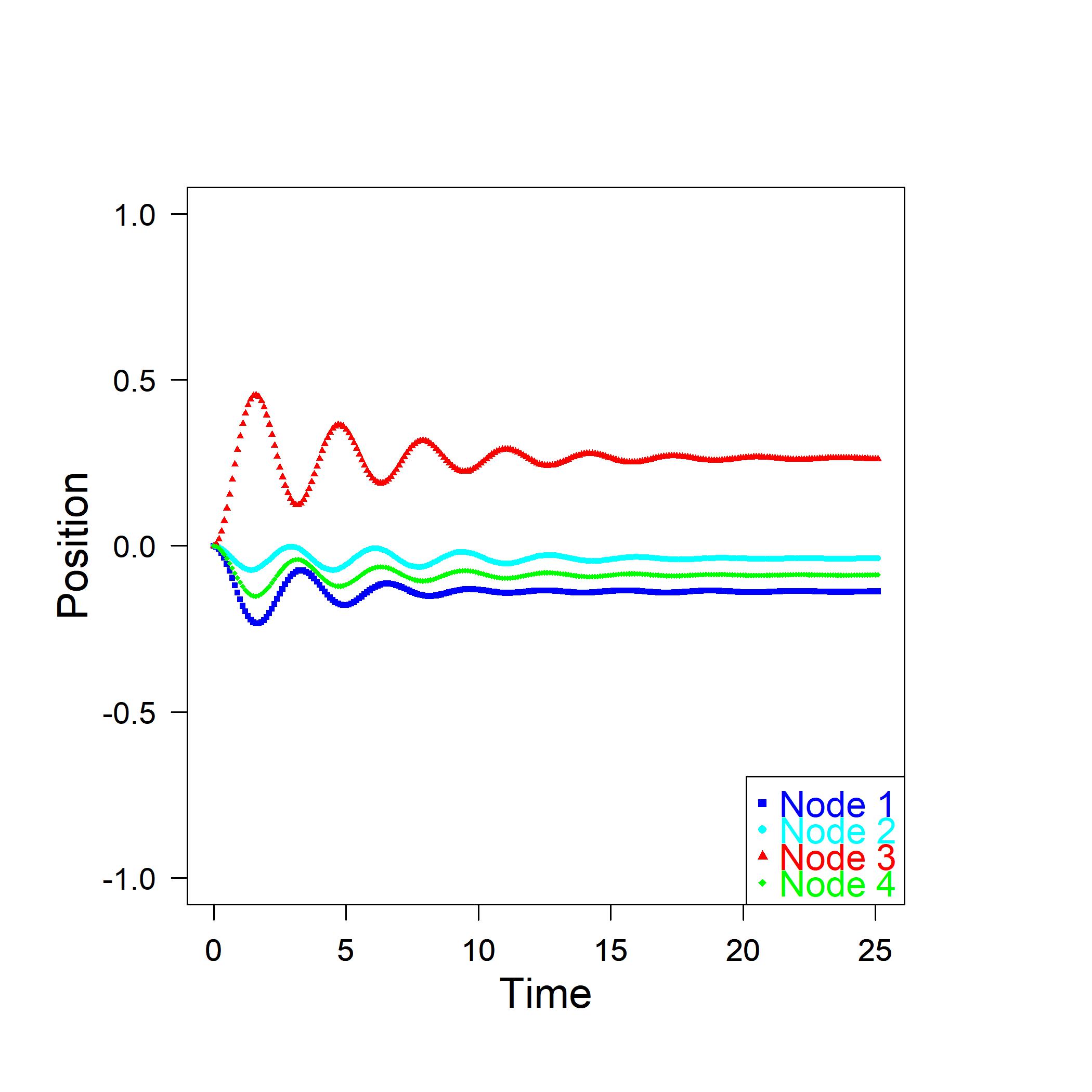}}
	\caption{Evolution of the displacement deviation of the toy network nodes
		for two different damping coefficients: panel (a) $\gamma=1$; panel (b) $\gamma=0.4$}
	\label{fig15} 
\end{figure}

Let us illustrate these ideas on the toy network in figure \ref{fig1}, postponing the study of a real power grid to the next section. Let us suppose that node $3$ is a generator and nodes $1$, $2$, and $4$ are loads and that the corresponding powers are ${\bf p}=[-0.50,-0.20,1.05,-0.35]^{T}$. Then a solution of ${\bf p}={\bf L}{\bf \tilde{x}}$ is, for instance, ${\bf \tilde{x}}=[-0.25,-0.15,0.15,-0.20]^{T}$. By substituting ${\bf \tilde{y}}=[{\bf \tilde{x}}, {\bf 0}]^{T}$ and ${\bf y}_{0}=[{\bf 0}, {\bf 0}]^{T}$  in Eq. (\ref{swing_equation}), we get the results depicted in figure \ref{fig15} for two different values of the damping parameter $\gamma$.

\section{Numerical Experiments}
\label{section7}

\subsection{Application to a social network}

Let us imagine that in a social network there is a debate around a certain topic and different members have different opinions on that topic.
Take, for example, the case in which there are different opinions about the qualities of a particular sportsman or of a given member within the same network.

Each member can shape a positive or negative assessment with all possible gradations from full agreement to full disagreement and he or she can share this opinion with his or her neighbors. We can model this static situation at a fixed time, for instance, by assigning $+1$ to perfect agreement and $-1$ to total disagreement and providing an intermediate graduation of the level of disagreement or agreement continuously from $-1$ to $+1$. The value $0$ would represent the totally neutral position.

However, each member is allowed to change his or her opinion over time, adjusting the assessment made at a given instant. Our goal here is to figure out what entrainment action some nodes may play with respect to others in this opinion formation process.

Following, for instance, \citet{Harary1956} and \citet{Altafini2012}, it is
reasonable to assume that, if a member has a close friend, he is much more likely to be influenced by the latter's opinion than by others' and therefore he will move in accordance with his close neighbors' opinion dynamics. In other words, he will be induced to synchronize the direction and the intensity of the velocity at which he changes his opinion according to that of his friends.

The mapping between the previous theoretical model and the social problem under examination is based on the interpretation of the continuous variables $x_i$ as a measure of the level of agreement or disagreement of each member with the subject. The opinion entrainment effect induced by a node on its neighbors - and, thus, on the network - is then modeled by the couplings described above. What we intend to focus on here is not the numerical value of each individual variable at a given instant. Rather we are interested in the collective behavior of the set of variables $x_i$ over time. Such behavior makes clear, as will be seen, the effect of the topological structure of the network in shaping processes of synchronization of members' ideas and opinions.

\begin{figure}[H]
	\includegraphics[width=0.90\linewidth]{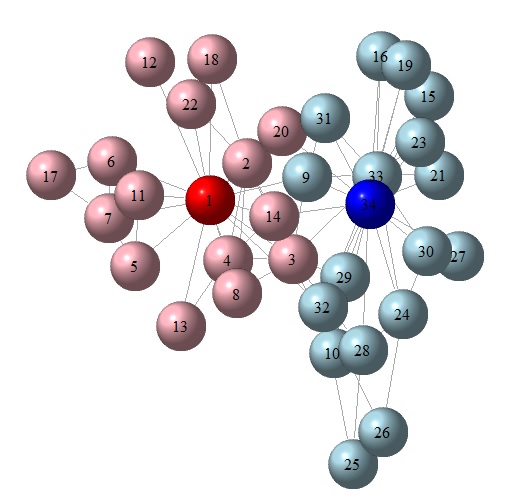}
	\centering
	\caption{Zachary Club network partition}
	\label{fig16}      
\end{figure}

A celebrated example of social network in network theory is provided by the Zachary Club network, see figure \ref{fig16}, in which $34$ members found themselves having to decide whether to agree with the position of the instructor ($+1$) or with the position of the president ($-1$). We can start by assigning an initial 'position' vector in such a way that node number $1$ in red, the instructor, has position $+1$ and node number $34$ in blue, the president, has position $-1$. At the beginning, at time $t=0$, all the other nodes are supposed to be totally neutral, since they still have to build an opinion, and so they are assigned a 'position' equal to $0$.

\begin{figure}[H]
\centering
	\subfloat[]{\includegraphics[width=0.25\textwidth]{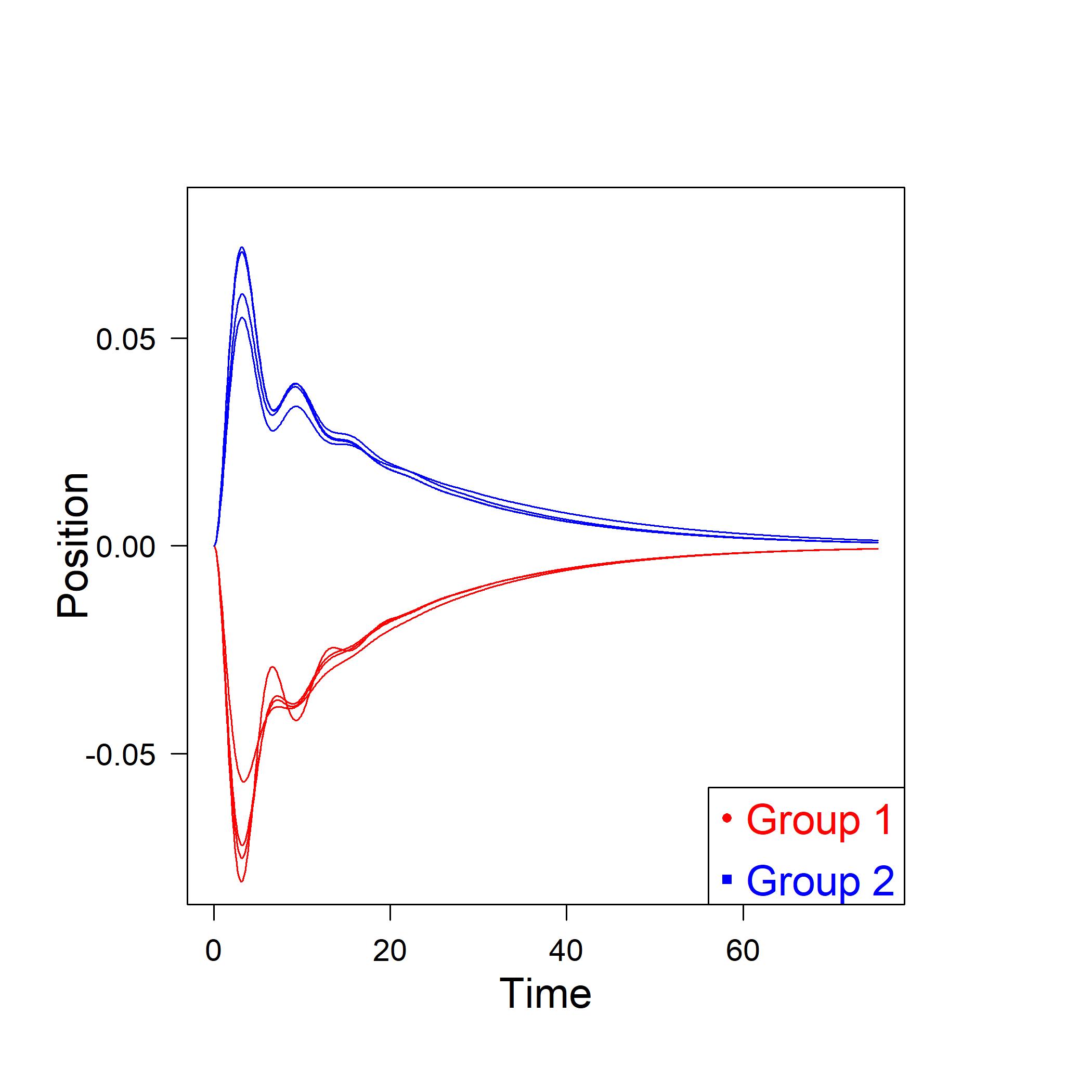}}
	\subfloat[]{\includegraphics[width=0.25\textwidth]{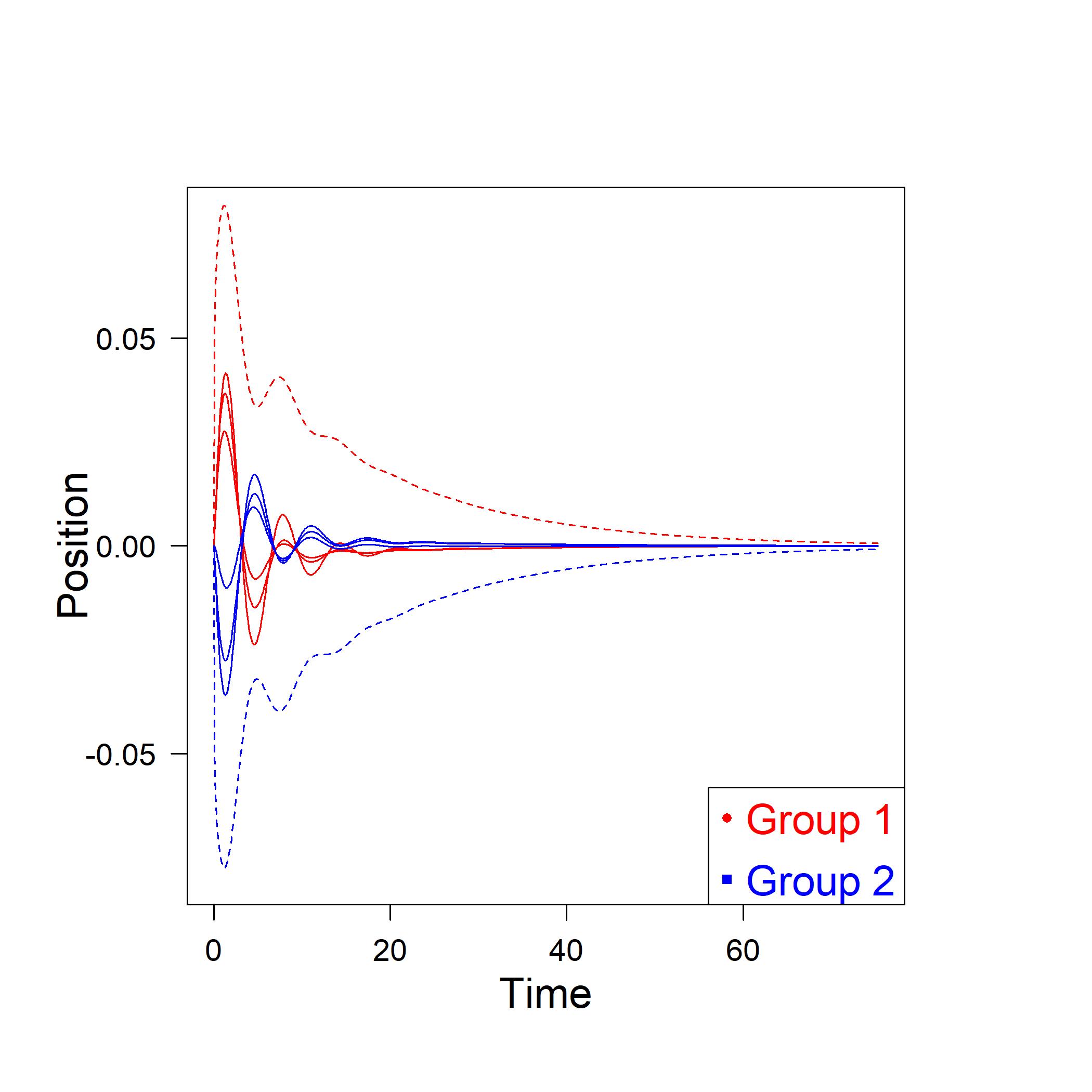}}
\caption{Synchronization in the two different communities: panel (a) initial positions equal to $+1$ for node number $1$, $-1$ for node number $34$, and $0$ for all the other nodes; panel (b) initial velocities equal to $+1$ for node number $1$, $-1$ for node number $34$, and $0$ for all the other nodes.}
	\label{fig17} 
\end{figure}

We start by assuming that the social system evolves according to Eq. (\ref{generalsolution5}) with a matrix $\bf G$ as in Eq. (\ref{matrixG1}).
The plot in figure \ref{fig17}, panel (a), illustrates the position vs time of two different groups of nodes with initial conditions equal to $+1$ for node number $1$ (the instructor), $-1$ for node number $34$ (the president), and $0$ for all the other nodes. Specifically, we choose four nodes for each one of the two communities illustrated in figure \ref{fig16} and colored in red and blue, respectively for the group of the instructor and the group of the president. The four nodes in red are nodes number $4,\ 11,\ 13,\ 18$ and the four nodes in blue are nodes number $23,\ 24,\ 30,\ 33$. As can be seen nodes in the two different groups tend to synchronize each other from the very beginning. As time goes on, both converge to zero due to the balanced initial conditions that see neither leader prevail over the other. What is interesting and unexpected is the fact that each group seems to assume a position which is opposite to that of their respective leader. However, it could be easily explained by saying that if node $1$ is initially in position $+1$ then it can move only in the negative verse toward $0$, dragging in the negative verse also nodes that are under its main influence. Similarly for node $34$.

The plot in figure \ref{fig17}, panel (b), illustrates the position vs time of two groups, including the two respective leaders. In this case, initial conditions see all the nodes having the same initial position equal to $0$ (including nodes $1$ and $34$); but we imagine now that at time $t=0$ node $1$ starts to move with velocity $+1$ (in the positive verse) and node $34$ starts to move with velocity $-1$ (in the negative verse).
In the plot, nodes in red are numbers $1,\ 4,\ 11,\ 18$ and nodes in blue numbers $34,\ 30,\ 32,\ 33$; dashed lines represent the two leader. As can be seen, nodes in the two different groups tend to synchronize each other, to synchronize with their leader and to stay out-of-phase with members of the other group for all the time, until both converge to $0$ due to the balanced initial conditions.

\begin{figure}[H]
\centering
	\subfloat[]{\includegraphics[width=0.25\textwidth]{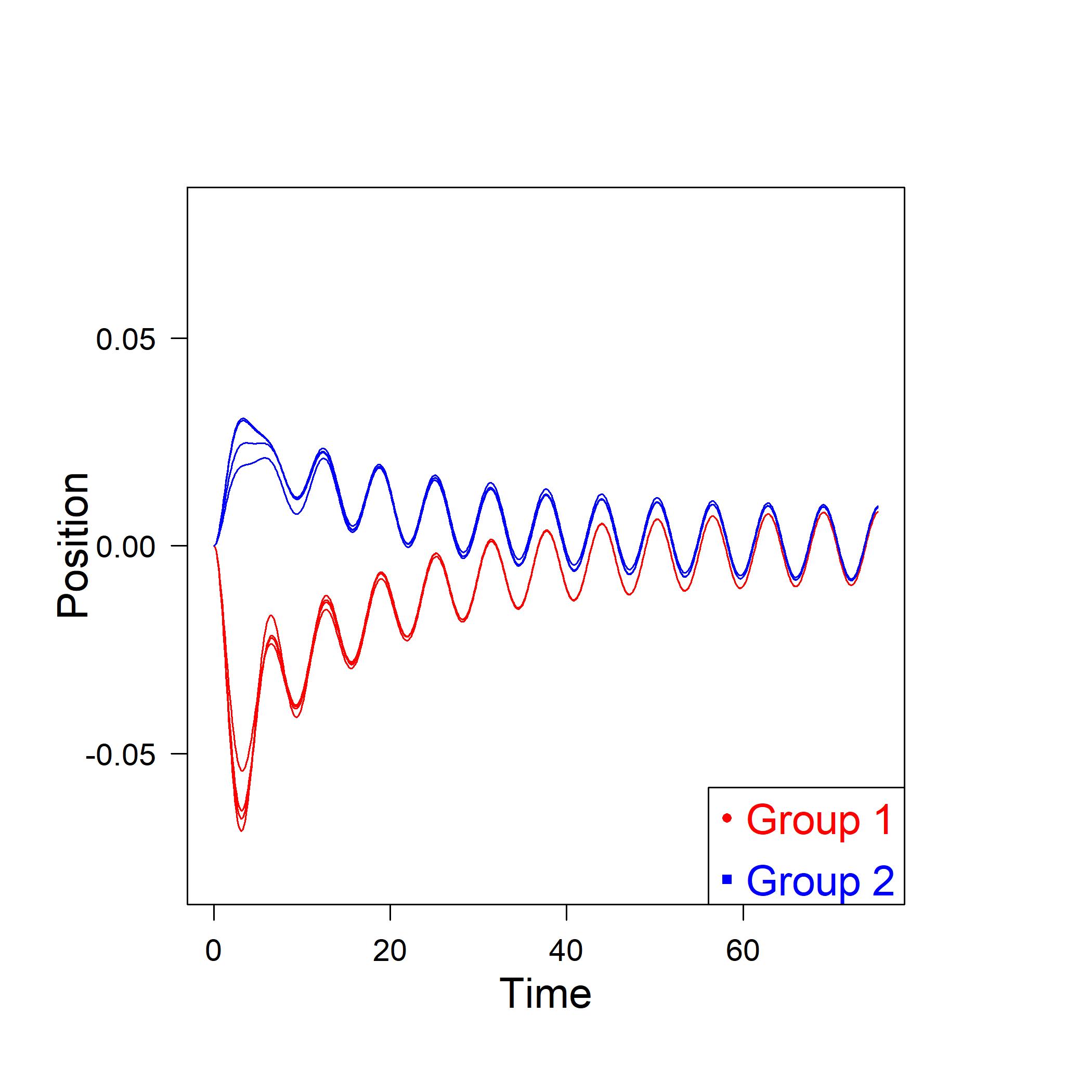}}
	\subfloat[]{\includegraphics[width=0.25\textwidth]{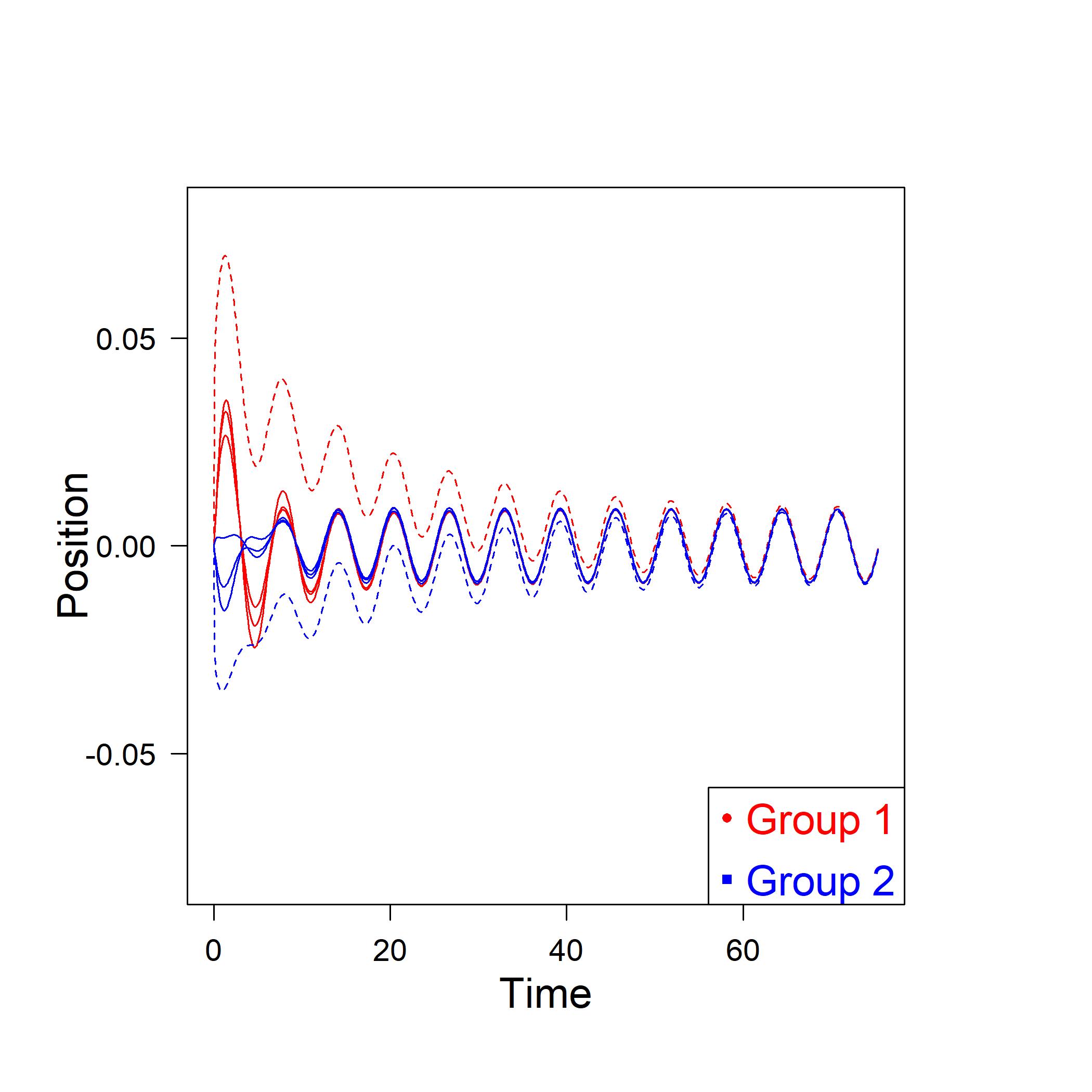}}
\caption{Synchronization in the two different communities: panel (a) node number $1$ and number $34$ initial positions equal to $0.8$ and $-0.5$, respectively; panel (b) node number $1$ and number $34$ initial velocities equal to $0.8$ and $-0.5$, respectively.}
	\label{fig18} 
\end{figure}

Two synchronization processes with unbalanced initial conditions are illustrated in figure \ref{fig18} for the same set of nodes, as in figure \ref{fig17}: in panel (a), node number $1$ and number $34$ have initial positions equal to $0.8$ and $-0.5$ respectively; in panel (b), node number $1$ and number $34$ have initial velocities equal to $0.8$ and $-0.5$ respectively.

We compare now the synchronization times given by Eq. (\ref{mean_time}) under different conditions. In detail, we are interested in understanding which of the two leaders is more effective in terms of speed of synchronization. Keeping in mind the discussion in the last remark in subsection \ref{Synchronization_time}, we start with a strong coupling between nodes $c_{2}^{\prime}=1$. If only the instructor (node number $1$) is endowed with an initial velocity (in our simulation equal to $4$) and all other nodes (including number $34$) are initially held at zero, then the time at which all nodes synchronize with the asymptotic solution in Eq. (\ref{asymptoticsolution}) below a threshold $\varepsilon=0.001$ is equal to $t=91.70$ with a mean time equal to $t=5.29$. If we now assign the same initial velocity to the president (node $34$) while keeping all others (including node number $1$) at zero, then the same time becomes $t=96.20$ with a mean time equal to $t=5.47$. The two times are basically comparable with slightly higher effectiveness for the instructor than for the president, confirming the fact that the network is divided into two balanced communities. The same conclusion is validated using other initial conditions on both the velocity and the position of the two leaders.

Let us now move to the resonant phenomena. We assume now that the social system evolves according to Eq. (\ref{synchronized_solution1}) with a matrix $\bf G$ as in Eq. (\ref{matrixG2}). The resonant frequencies given in proposition \ref{proposition3} range from $1$ to $4.37$ and they are listed in the Table \ref{table4}.
\begin{table}[H]
	\begin{center}
		\begin{tabular}{||c|c|c||c|c|c||c|c|c||}
			\hline
	 $i$		& ${\mu_i}$ &  $\omega_{i}$ & $i$ &  ${\mu_i}$ &  $\omega_{i}$  & $i$	& ${\mu_i}$ & $\omega_{i}$ \tabularnewline \hline
$1$ & $18.137$ & $4.37455$ & $11$ & $4.581$ & $2.36237$ & $21$ & $2.000$	 & $1.73205$   \tabularnewline \hline
$2$ & $17.055$ & $4.24914$  & $12$ & $4.480$ & $2.34094$ & $22$ & $1.955$	 & $1.71903$ \tabularnewline \hline
$3$ & $13.306$ & $3.78234$  & $13$ & $4.276$ & $2.29693$ & $23$ & $1.826$	 & $1.68109$  \tabularnewline \hline
$4$ & $10.921$ & $3.45269$  & $14$ & $3.472$ & $2.11476$
&$24$ & $1.762$	 & $1.66190$ \tabularnewline \hline
$5$ & $9.777$	 & $3.28287$  & $15$ & $3.382$	& $2.09332$ & $25$ & $1.599$	 & $1.61223$ \tabularnewline \hline
$6$ & $6.996$	 & $2.82776$  & $16$ & $3.376$	 & $2.09193$ & $26$ & $1.259$	 & $1.50313$  \tabularnewline \hline
$7$ & $6.516$	 & $2.74145$  & $17$ & $3.242$	 & $2.05963$ & $27$ & $1.125$	 & $1.45774$ \tabularnewline \hline
$8$ & $6.332$	 & $2.70769$  & $18$ & $3.014$	 & $2.00349$ & $28$ & $0.909$	 & $1.38176$ \tabularnewline \hline 
$9$ & $5.618$	 & $2.57255$    &$19$ & $2.749$	 & $1.93627$  &$29$ & $0.468$	 & $1.21183$ \tabularnewline \hline
$10$ & $5.379$	 & $2.52559$ & $20$ & $2.487$	 & $1.86738$   & $30$ & $0$  	 & $1.00000$   \tabularnewline \hline
\end{tabular}
\end{center}
\caption{Eigenvalues of the matrix $\bf L$ and corresponding resonance frequencies for the Zachary Club network. Note that eigenvalue $2.000$ and the corresponding frequency have multiplicity equal to $5$.}
\label{table4}
\end{table}

Let us imagine now that an external force with frequency $\omega$ is applied to node number $1$, or, in other words, that the first leader acts as an influencer by applying the same driving force in that point of the network. Let us consider two nodes belonging to the two different clusters previously considered: node number $11$ in the first cluster (red colored) and node number $30$ in the second cluster (blue colored). Let us consider their behavior as $\omega$ approaches two resonant frequencies, specifically $\omega_{30}=1$ and $\omega_{22}=1.719026$. Results are illustrated in the next figures \ref{fig19} and \ref{fig20}.

As discussed in the remark after proposition \ref{proposition3}, all nodes resonate at the ground frequency $\omega=1$. Node number $11$, which belongs to the sphere of influence of the first leader, certainly resonates at the ground frequency and so does at the second frequency under consideration; on the contrary, node number $30$, which falls within the sphere of influence of the second leader, while behaving similarly to node $11$ at the ground frequency, does not respond to the second frequency from node $1$. This is easily explained if we look at the component of the corresponding eigenvector $\phi_{22}$. Indeed, for node number $11$, it is $\phi_{22}(11)=-0.540899038$, whereas, for node number 30, it is $\phi_{22}(30)=0.002692861$. The first component is about 200 times larger than the second one, which justifies the very low influence of node 1 on node 30 at that frequency.

\begin{figure}[H]
	\centering
	\subfloat[]{\includegraphics[width=0.30\textwidth]{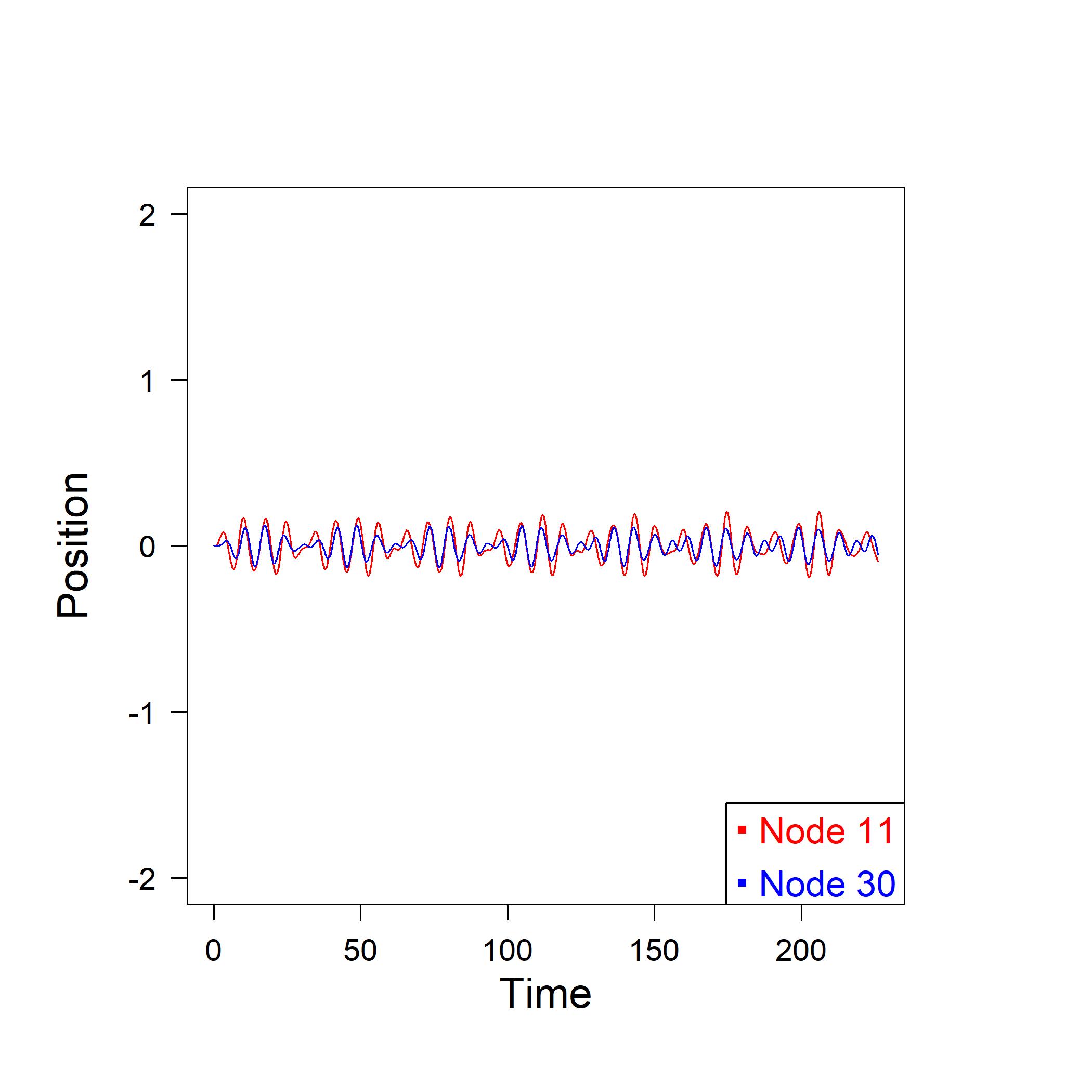}}
	\\
	\subfloat[]{\includegraphics[width=0.30\textwidth]{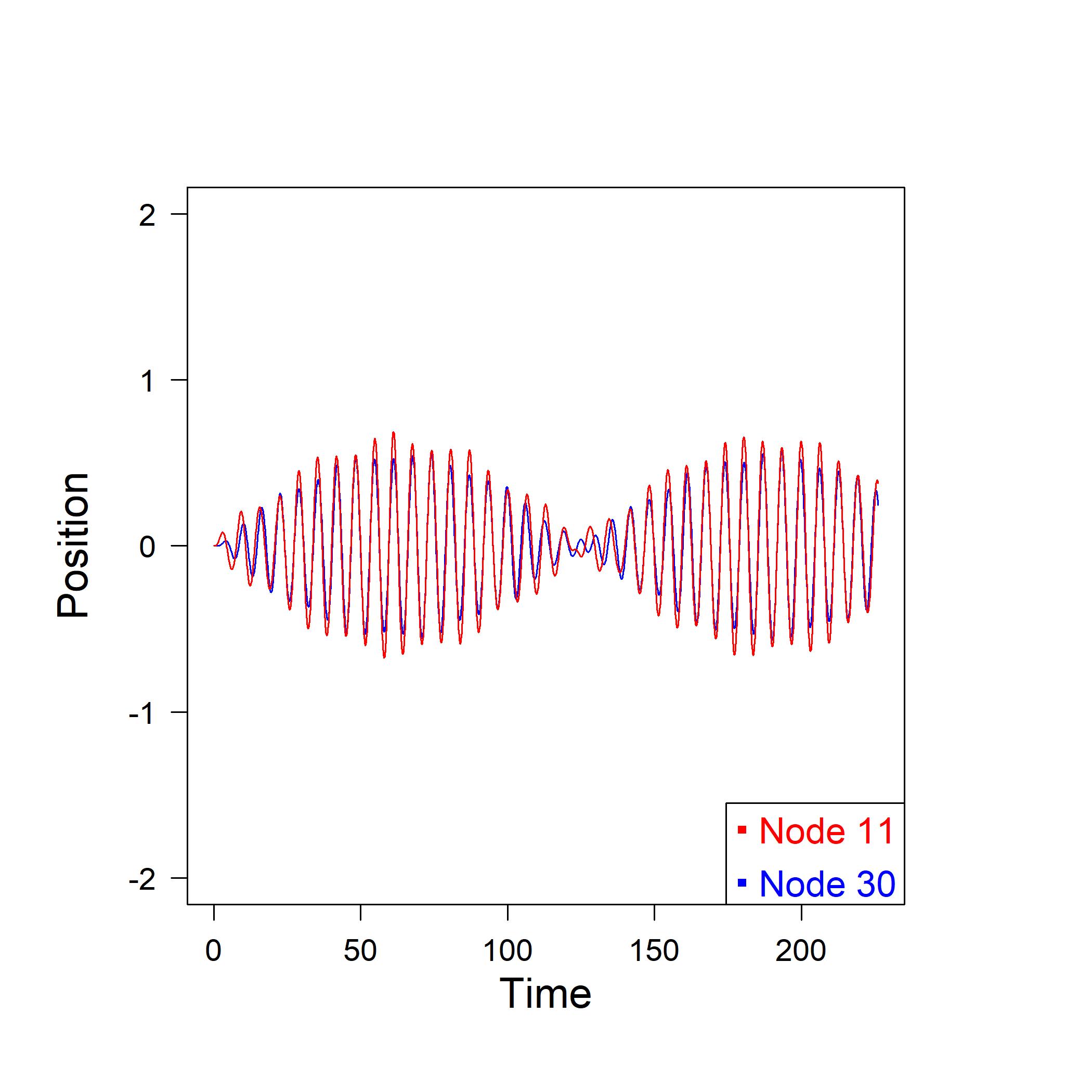}}
	\\
	\subfloat[]{\includegraphics[width=0.30\textwidth]{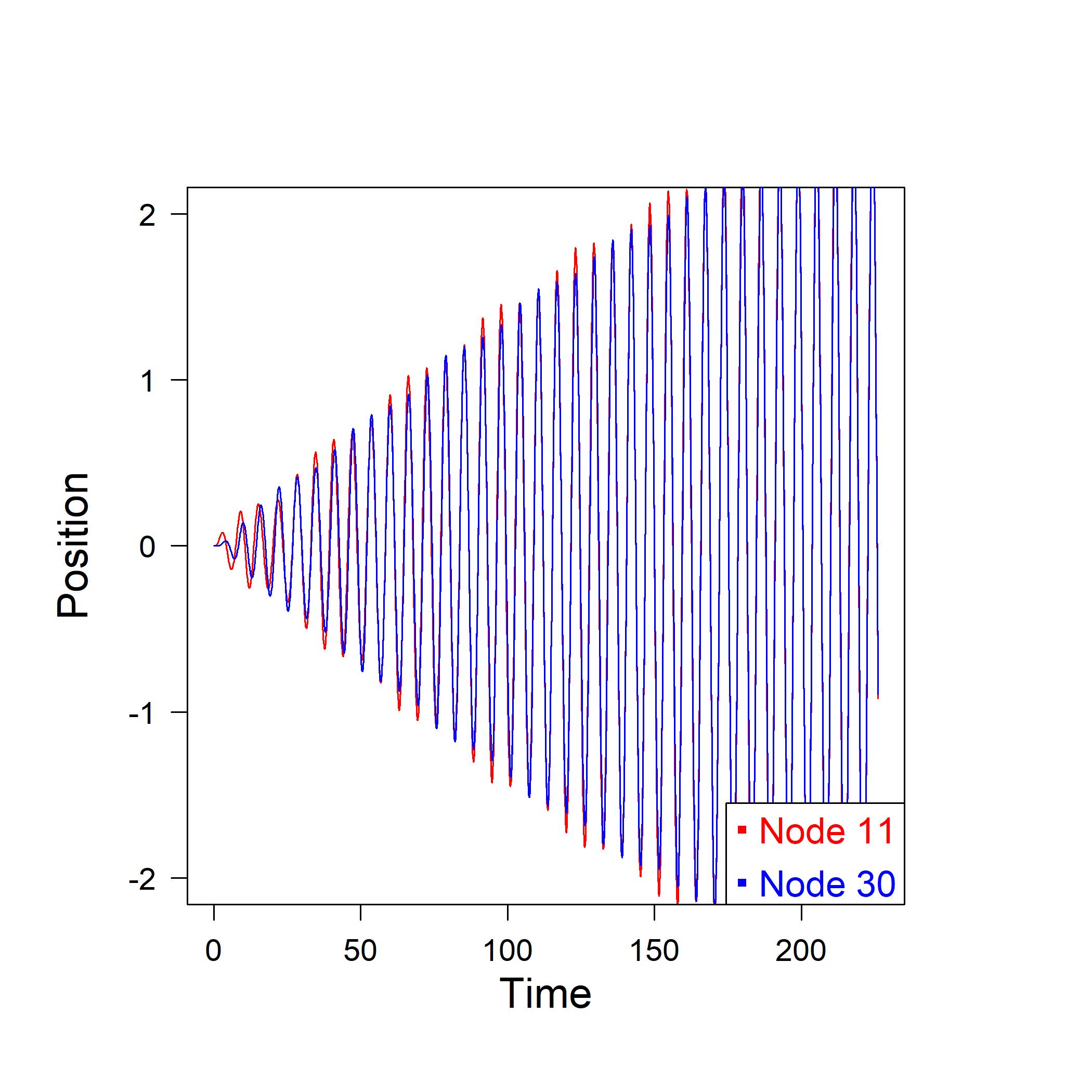}}
	\caption{Approaching the resonant state with a driving force in node number $1$ of the Zachary network, and with ${\bf x}_{0}={\bf 0}$ and ${\bf v}_{0}={\bf 0}$; the frequencies are: panel (a) $\omega=0.80$; panel (b) $\omega=0.95$; and panel (c) $\omega=0.99$.}
	\label{fig19} 
\end{figure}

\begin{figure}[H]
	\centering
	\subfloat[]{\includegraphics[width=0.30\textwidth]{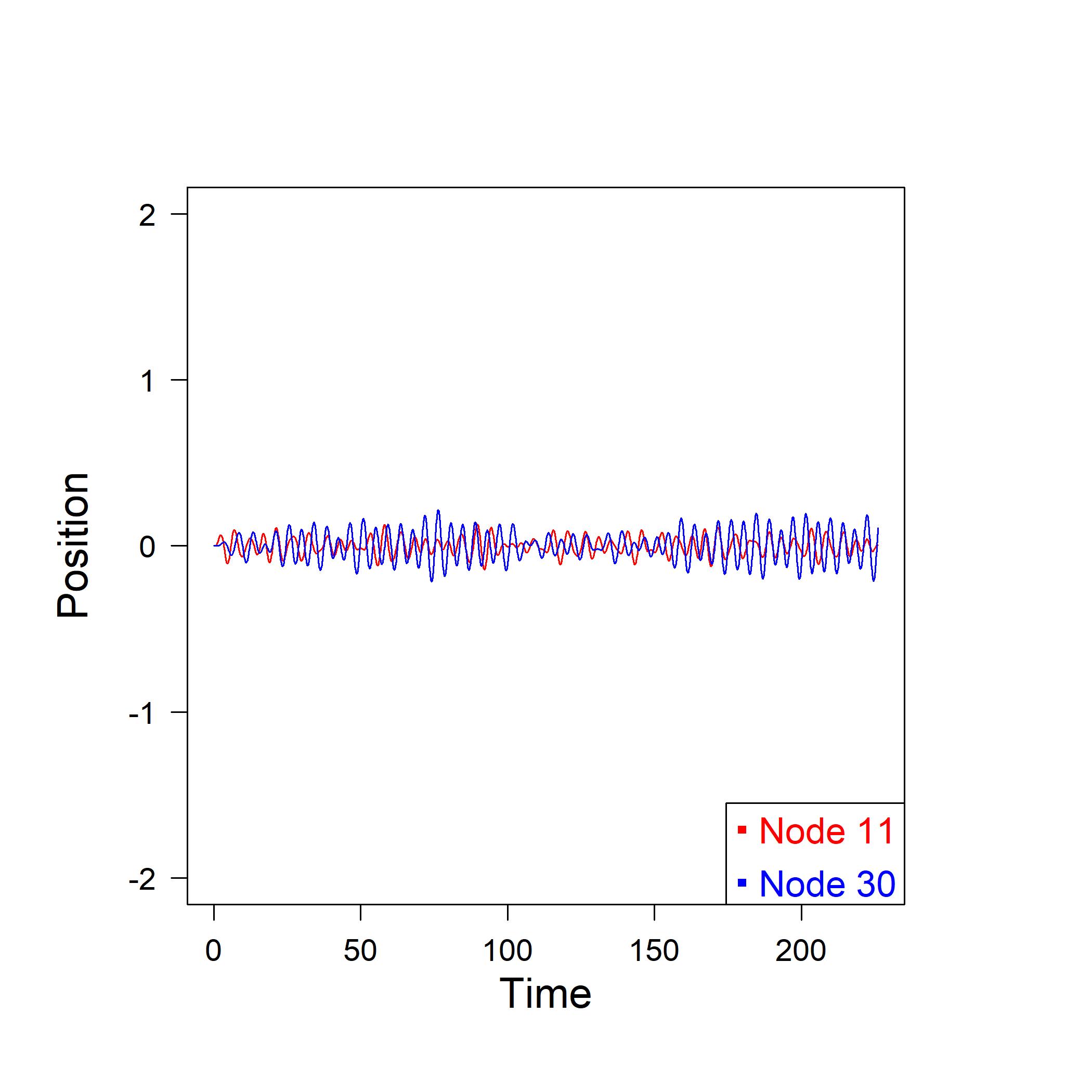}}
	\\
	\subfloat[]{\includegraphics[width=0.30\textwidth]{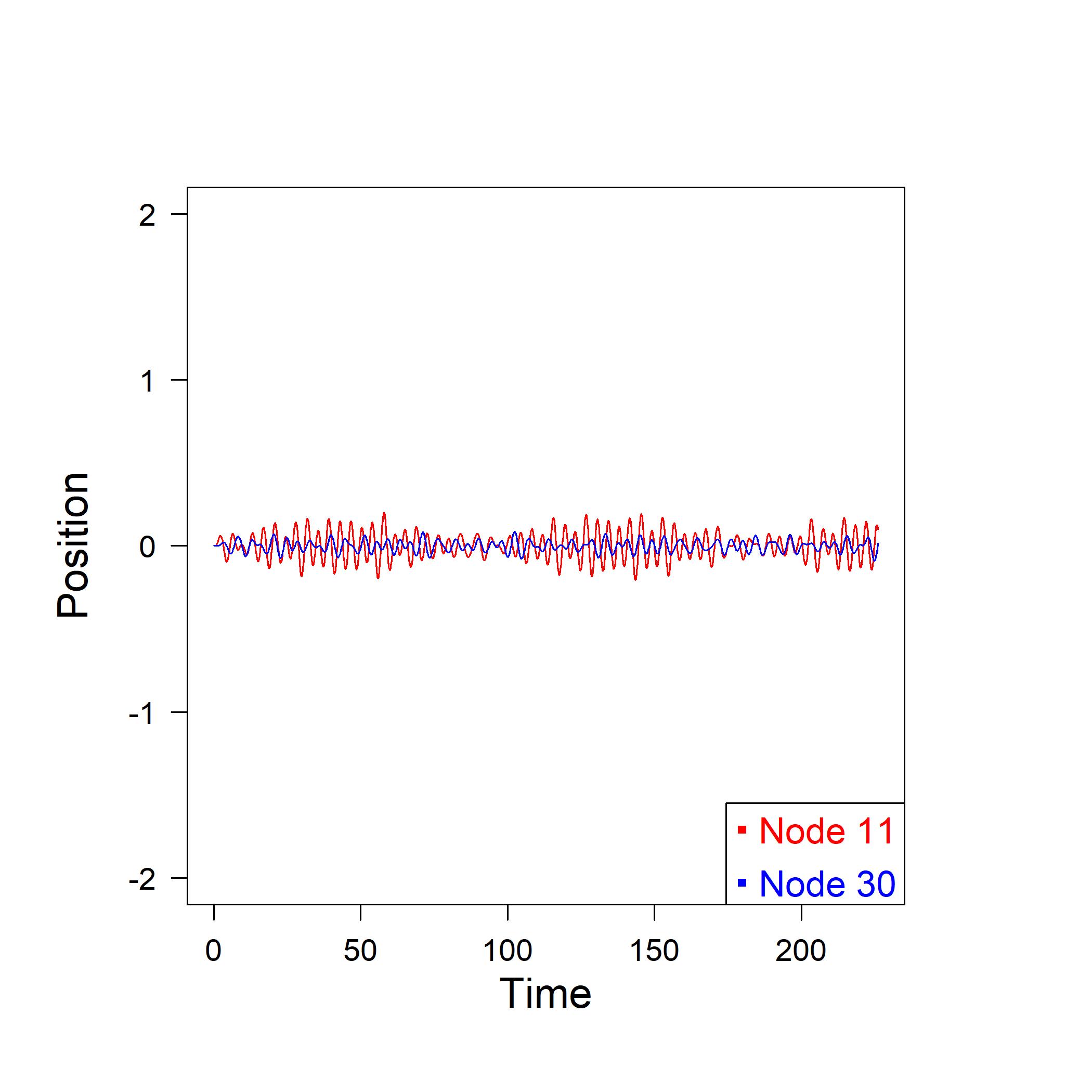}}
	\\
	\subfloat[]{\includegraphics[width=0.30\textwidth]{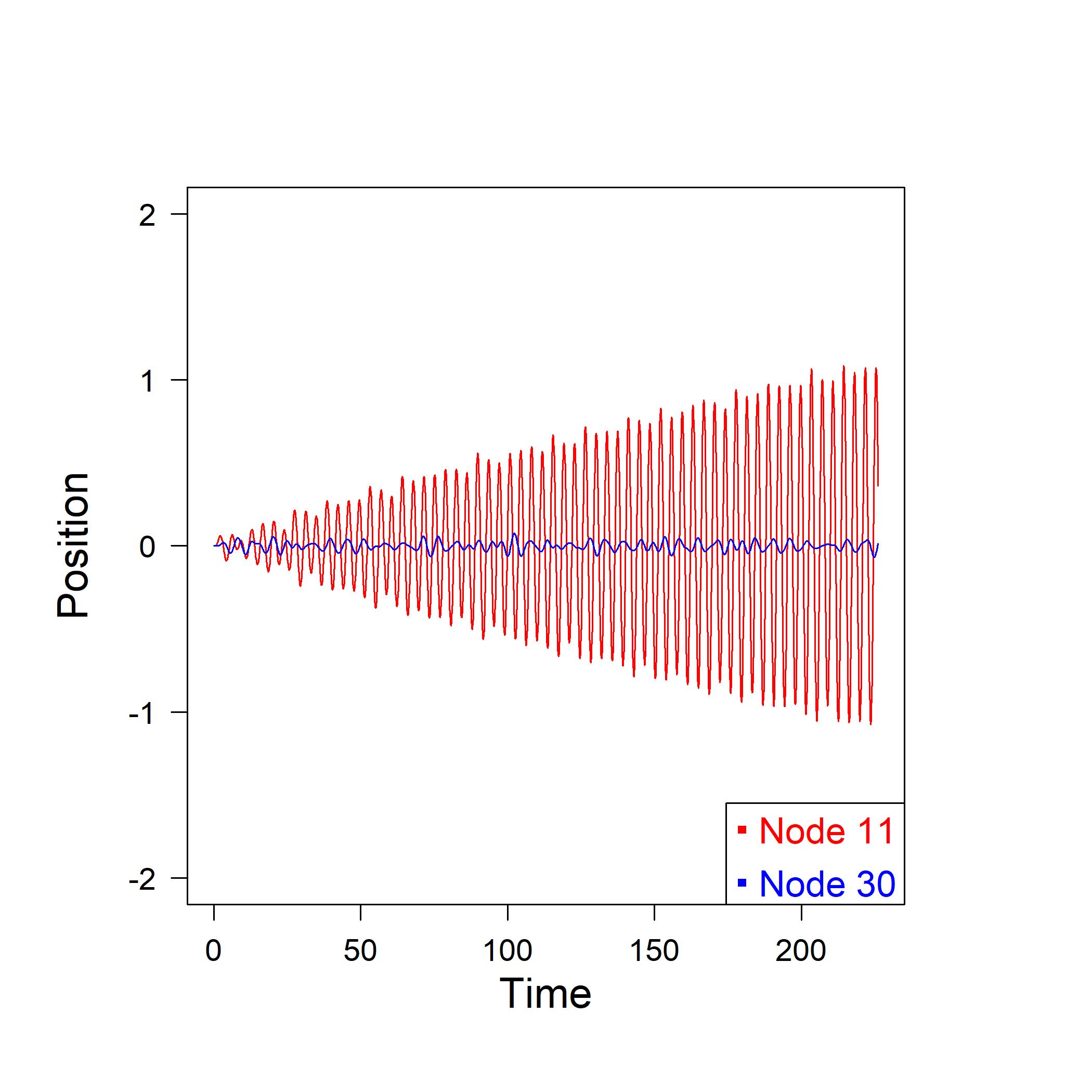}}
	\caption{Approaching the resonant state with a driving force in node number $1$ of the Zachary network, and with ${\bf x}_{0}={\bf 0}$ and ${\bf v}_{0}={\bf 0}$; the frequencies are: panel (a) $\omega=1.50$; panel (b) $\omega=1.65$; and panel (c) $\omega=1.71$.}
	\label{fig20}
\end{figure}

\begin{figure}[H]
\centering
	\subfloat[]{\includegraphics[width=0.30\textwidth]{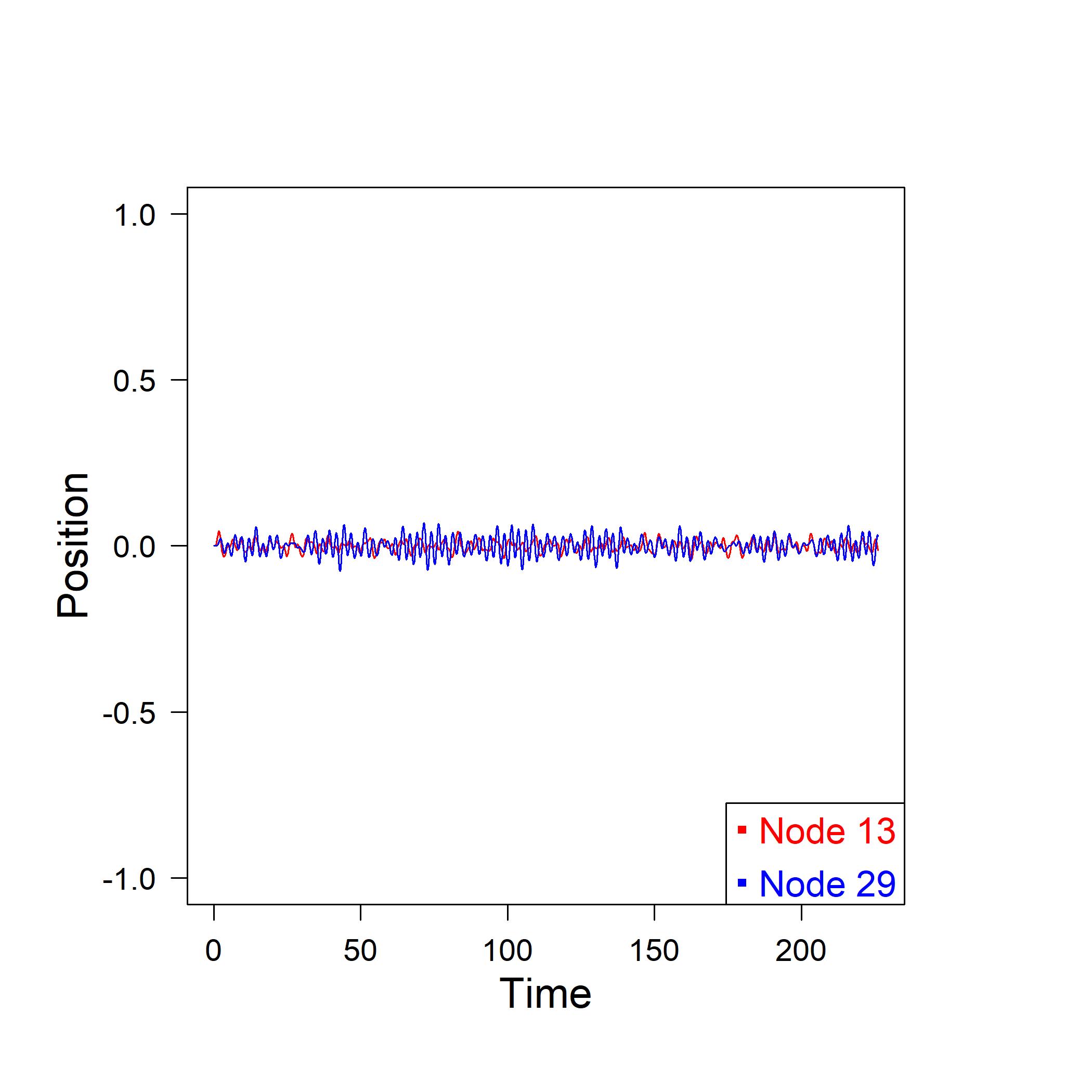}}
	\\
	\subfloat[]{\includegraphics[width=0.30\textwidth]{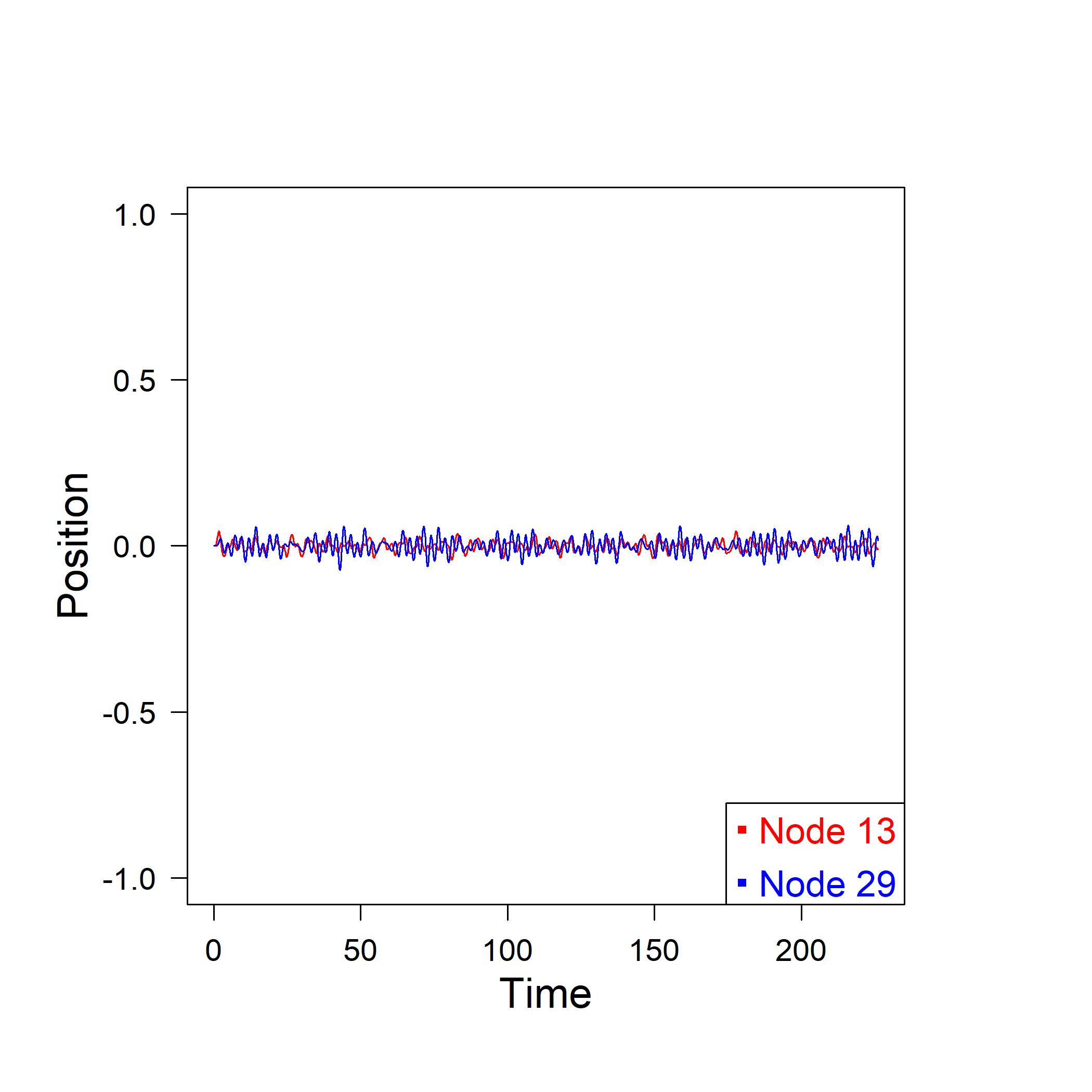}}
	\\
	\subfloat[]{\includegraphics[width=0.30\textwidth]{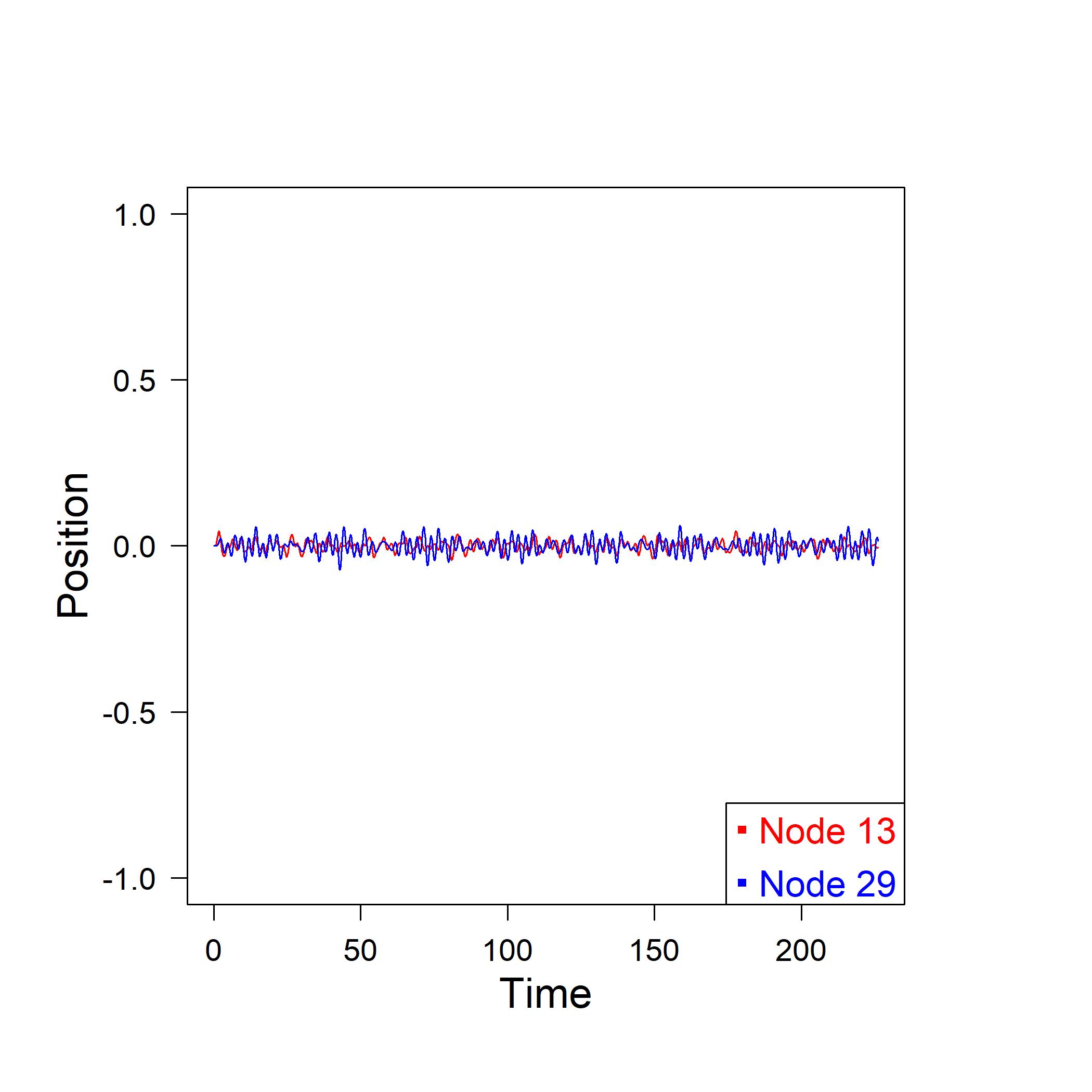}}
\caption{Approaching the resonant state with a driving force in node number $1$ of the Zachary network, and with ${\bf x}_{0}={\bf 0}$ and ${\bf v}_{0}={\bf 0}$; the frequencies are: panel (a) $\omega=2.560$; panel (b) $\omega=2.570$; and panel (c) $\omega=2.572$.}
	\label{fig21} 
\end{figure}

Let us consider now the frequency $\omega_{9}=2.572554$. Since $\phi_{9}(1)=0$, we expect that the oscillation of node $1$ with that frequency does not affect the network. This is confirmed by figure \ref{fig21}, where we selected two arbitrary nodes (numbers 13 and 29) in the two different clusters, which show no amplification of their own oscillation when the frequency of node 1 approaches the indicated value.

\begin{figure}[H]
\centering
	\subfloat[]{\includegraphics[width=0.30\textwidth]{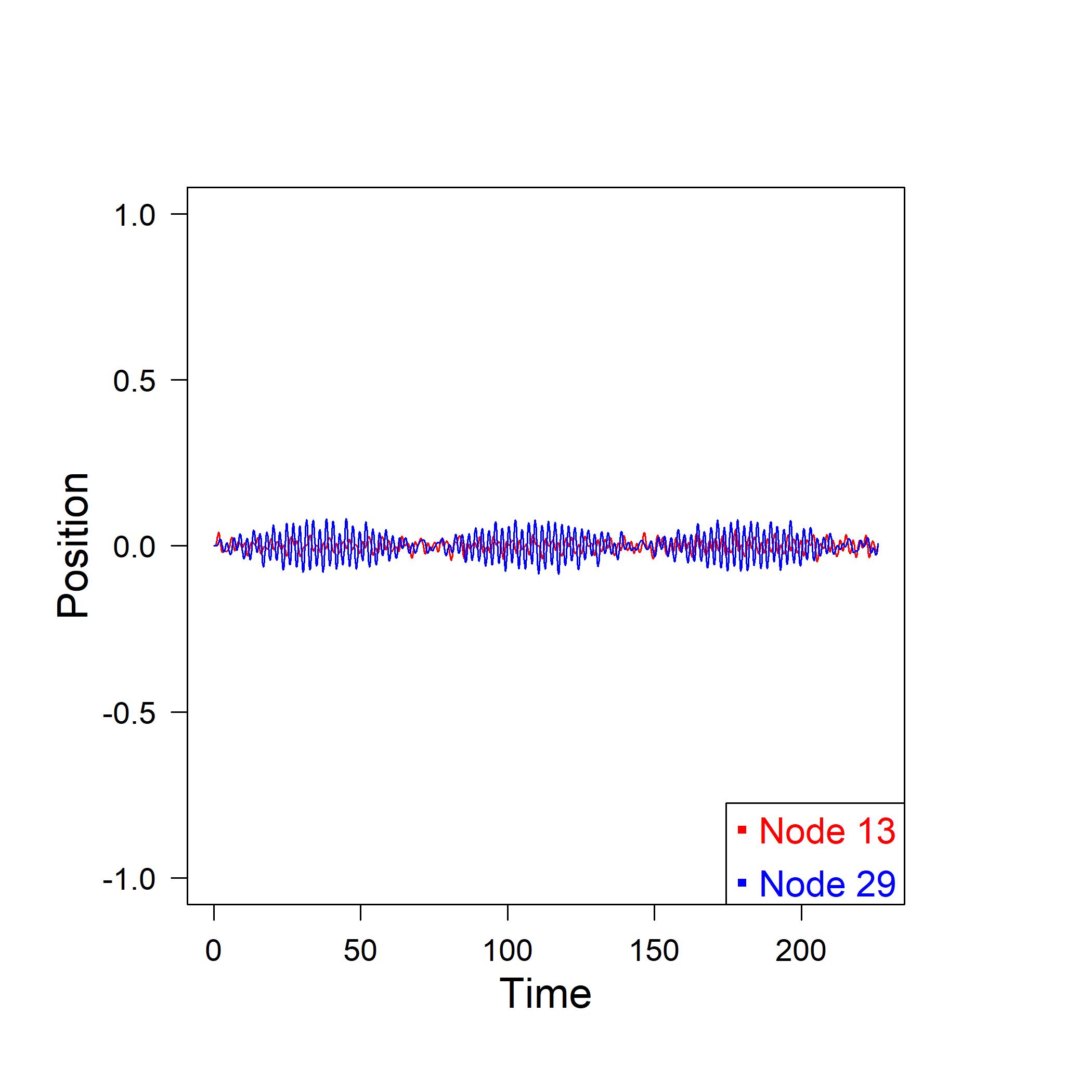}}
	\\
	\subfloat[]{\includegraphics[width=0.30\textwidth]{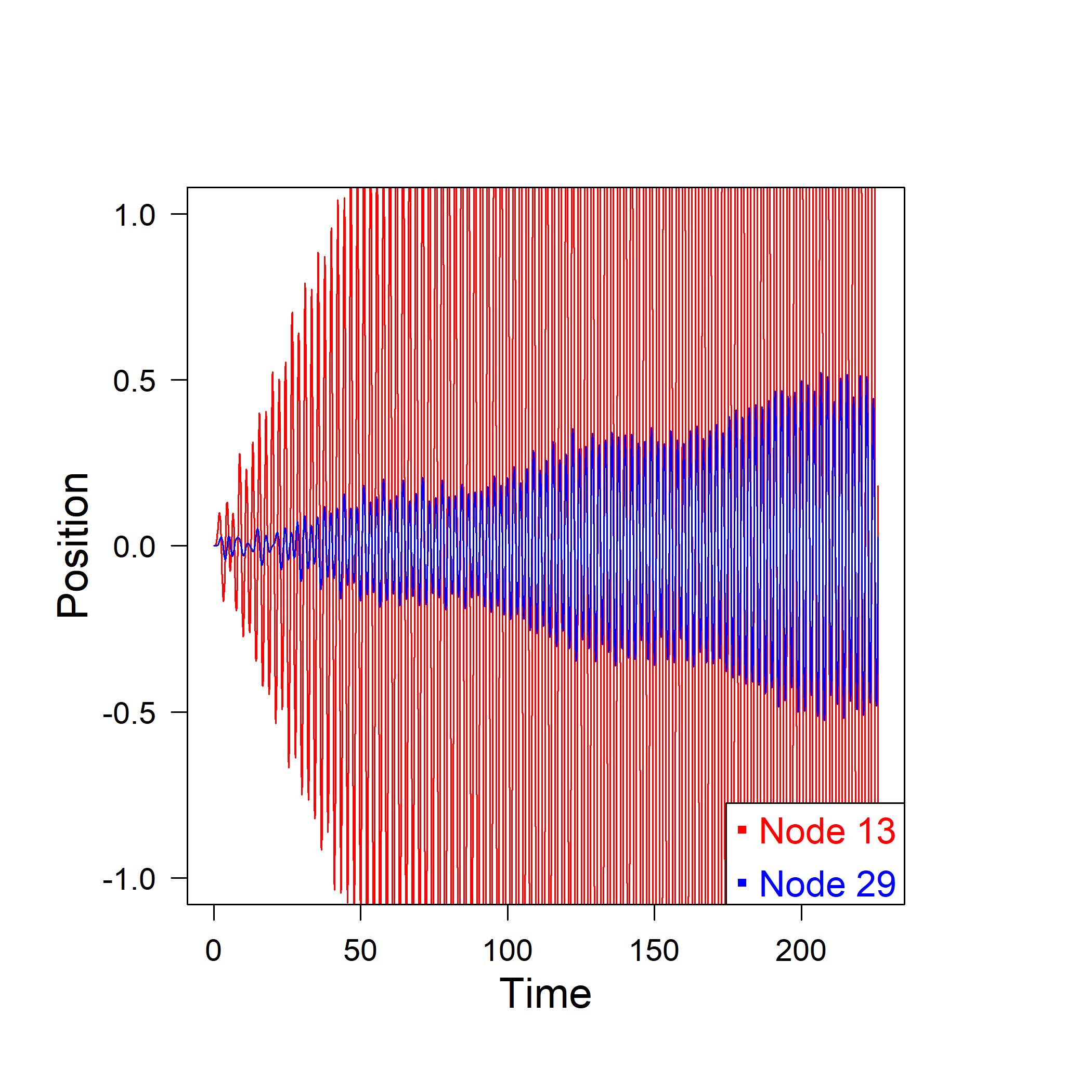}}
	\\
	\subfloat[]{\includegraphics[width=0.30\textwidth]{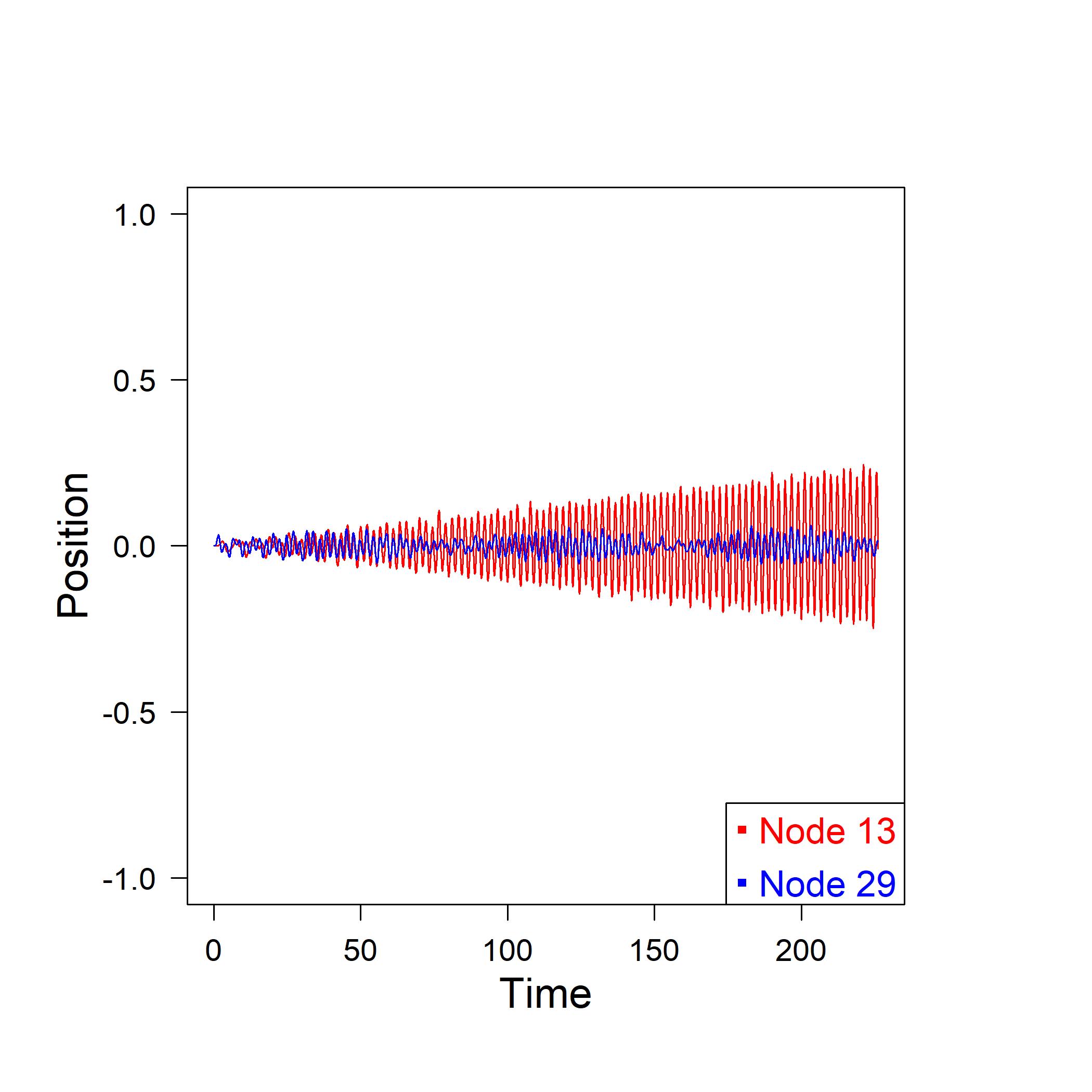}}
\caption{Approaching the resonant state with a driving force of fixed frequency $\omega=2.82775$, and with ${\bf x}_{0}={\bf 0}$ and ${\bf v}_{0}={\bf 0}$; the source nodes are: node number $1$ in panel (a), node number $4$ in panel (b) and node number $34$ in panel (c).}
	\label{fig22} 
\end{figure}

Further confirmation comes from the frequency $\omega_{6}=2.827755$. Since $\phi_{6}(1)=-0.0027808652$ and $\phi_{6}(4)=-0.8231724260$, we expect that the network globally responds much more if that oscillation starts from node $4$ than from node $1$. Figure \ref{fig22} shows the different behavior of two randomly selected nodes (again numbers 13 and 29) when such a frequency is applied on the network starting from different source nodes, specifically nodes 1, 4 and 34. The amplitude of the oscillation of nodes $13$ and $29$ grows much more quickly when the source is placed in node $4$ than when it is placed in node $34$ or, worse, in node $1$. The two leaders prove to be not so effective at this specific frequency.

\begin{figure}[H]
\centering
	\subfloat[]{\includegraphics[width=0.30\textwidth]{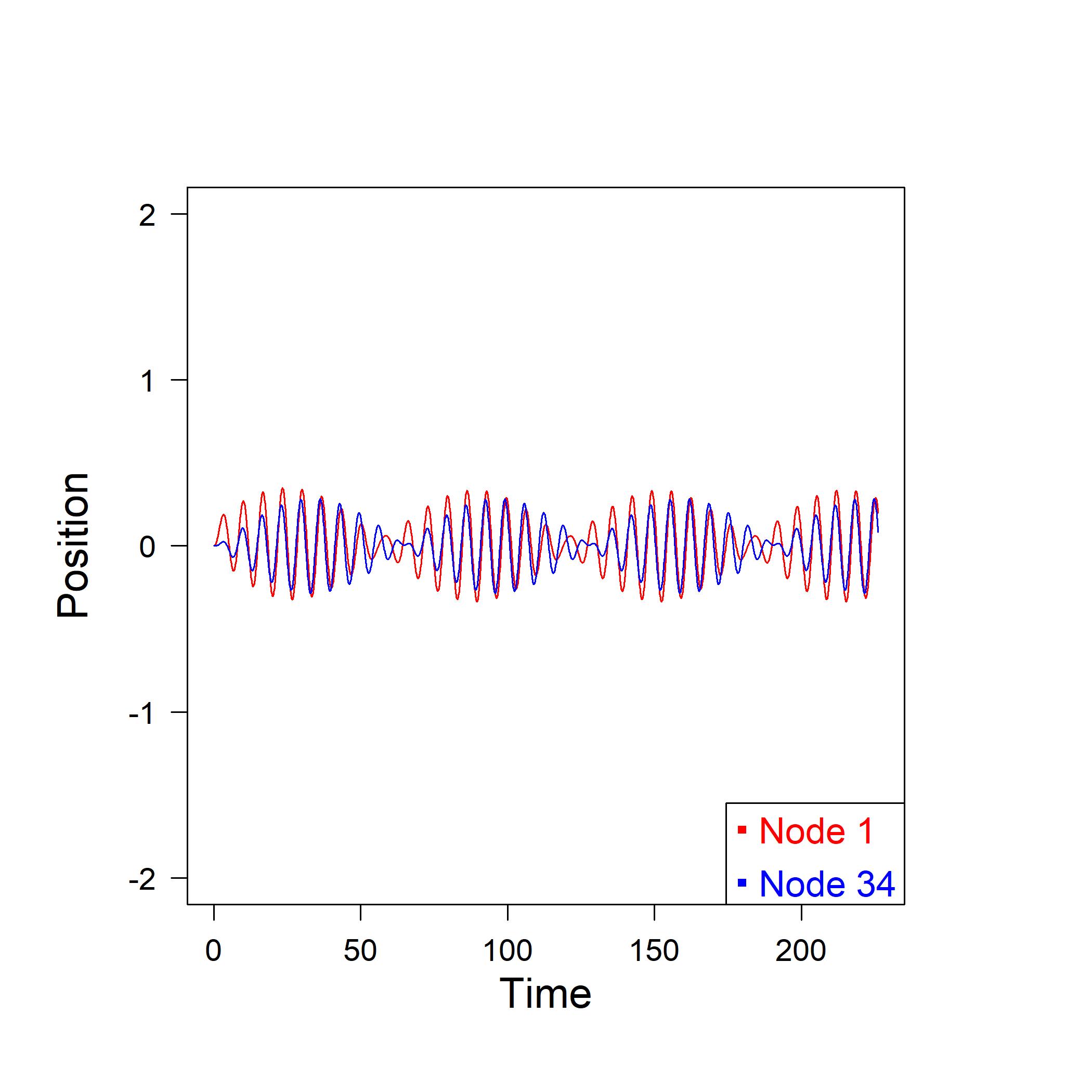}}
	\\
	\subfloat[]{\includegraphics[width=0.30\textwidth]{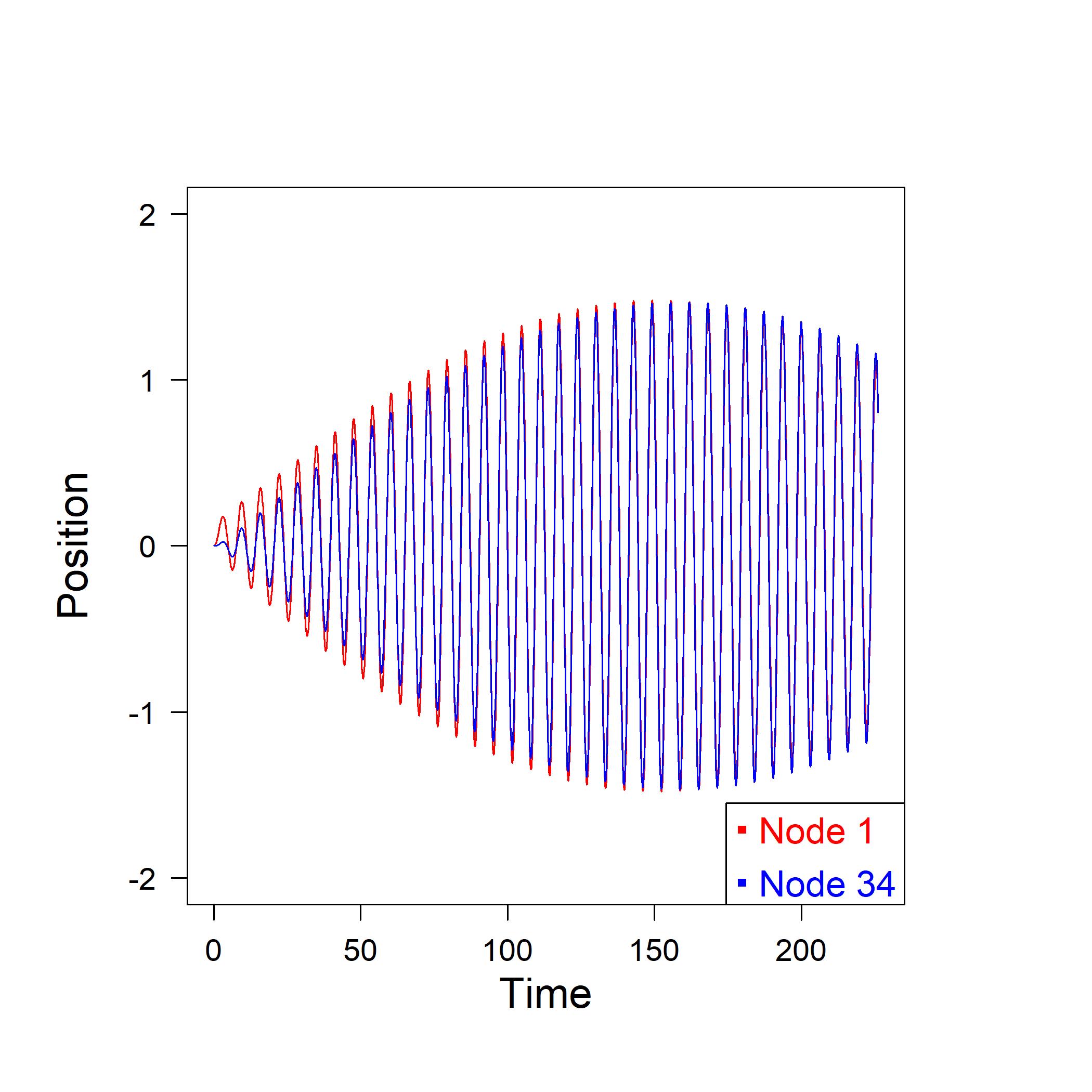}}
	\\
	\subfloat[]{\includegraphics[width=0.30\textwidth]{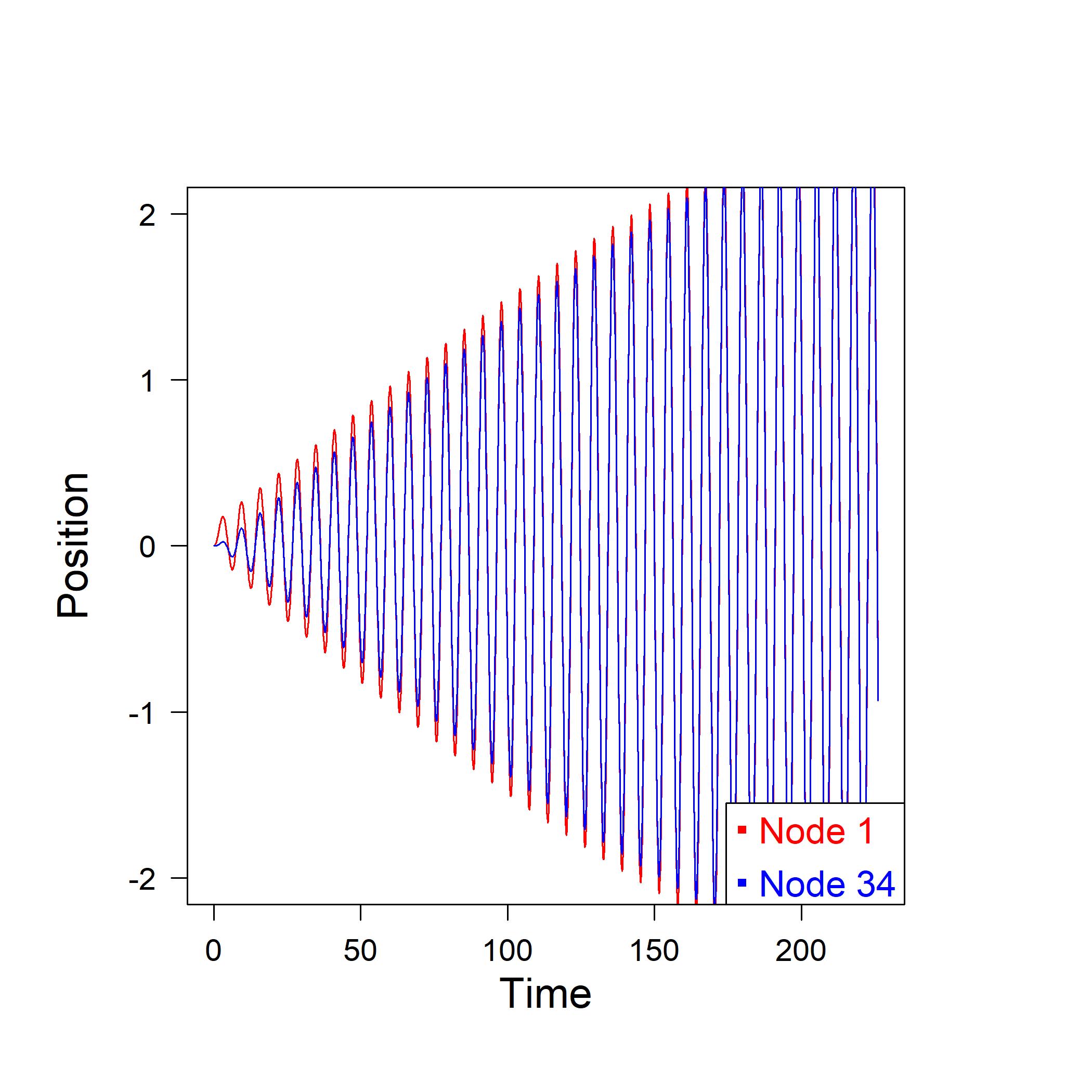}}
\caption{Approaching the synchronized resonant state with a driving force in node number $1$ of the Zachary network and ${\bf x}_{0}={\bf 0}$ and ${\bf v}_{0}={\bf 0}$; the frequencies are: panel (a) $\omega=0.90$; panel (b) $\omega=0.98$; and panel (c) $\omega=0.99$.}
	\label{fig23} 
\end{figure}

Finally, let us consider the case of the resonant synchronization described by proposition \ref{proposition4}, by assuming now that the social system evolves according to Eq. (\ref{separatesolution2}) with a matrix $\bf G$ as in Eq. (\ref{matrixG3}). Figure \ref{fig23} shows the behavior of the two leaders when the first one approaches frequency $\omega =1$ from below. All the other nodes display a similar behavior and enter a synchronized resonant state with the two leaders. No other frequency is able to produce a state of synchronized resonance on all nodes in the network.

\subsection{Application to a power grid}

\begin{figure}[H]
	\centering
	\subfloat[]{\includegraphics[width=0.45\textwidth]{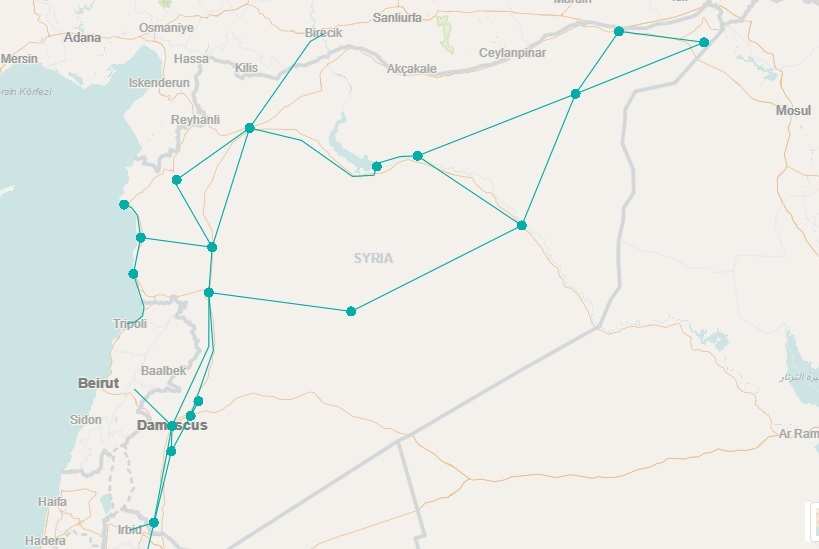}}
	\\
	\subfloat[]{\includegraphics[width=0.45\textwidth]{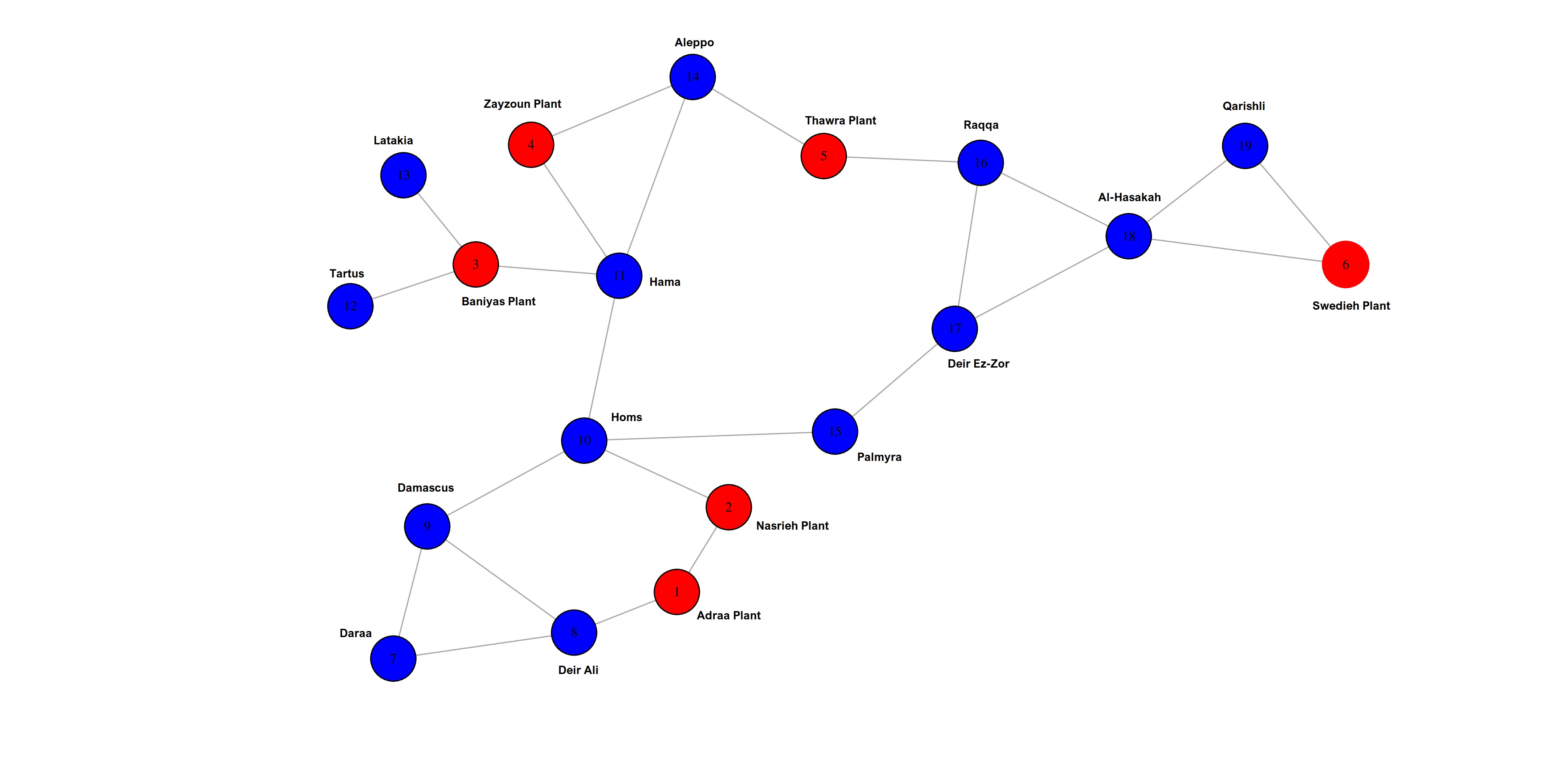}}
	\caption{(a) The main structure of the Syrian power grid; (b) the topology of the same network: nodes in red are generators and nodes in blue are loads.}
	\label{fig24} 
\end{figure}

We apply now results discussed in Section \ref{section6} to a small power grid, namely to the main electricity transmission network of Syria, represented in figure \ref{fig24} panel (a).\footnote{Credits. Source: Arab Union of Electricity;  Author: Haroun Beltaifa; Last Updated: April 13, 2020;  Country: Syrian Arab Republic; License: Creative Commons Attribution 4.0; url: https://energydata.info/dataset/syria-electricity-transmission-network-2017; Disclaimer: the author did not modify the original data.} A power grid consists of a network whose nodes are rotating machines and edges are transmission lines. Nodes can produce power (generating plants) or consume power (loads). The network topology of the Syrian grid is depicted in figure \ref{fig24} panel (b). It consists of $19$ nodes with $6$ generators (nodes 1-6, in red), $13$ loads (nodes 7-19, in blue), and $24$ edges.

\begin{figure}[H]
	\centering
	\subfloat[]{\includegraphics[width=0.36\textwidth]{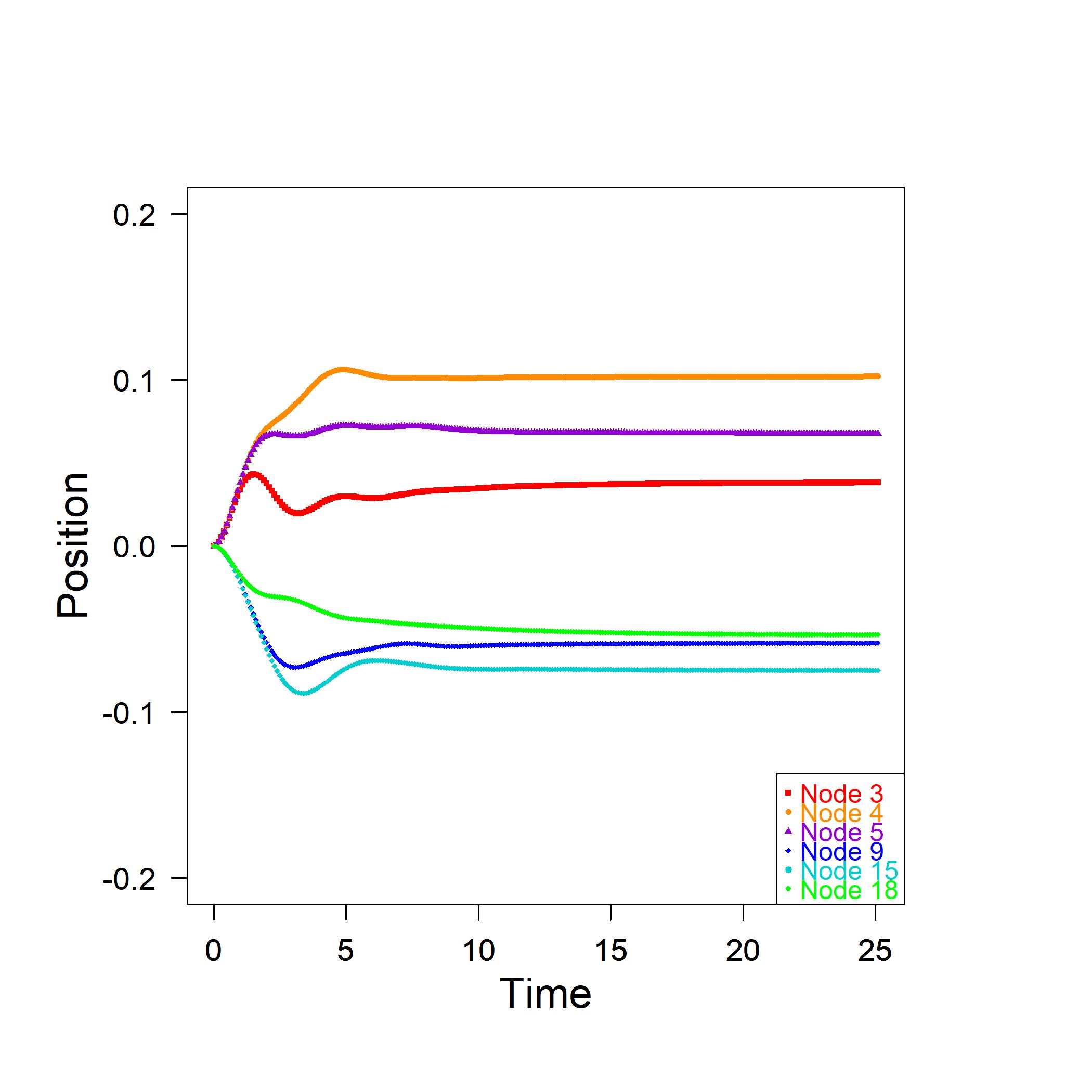}}
	\\
	\subfloat[]{\includegraphics[width=0.36\textwidth]{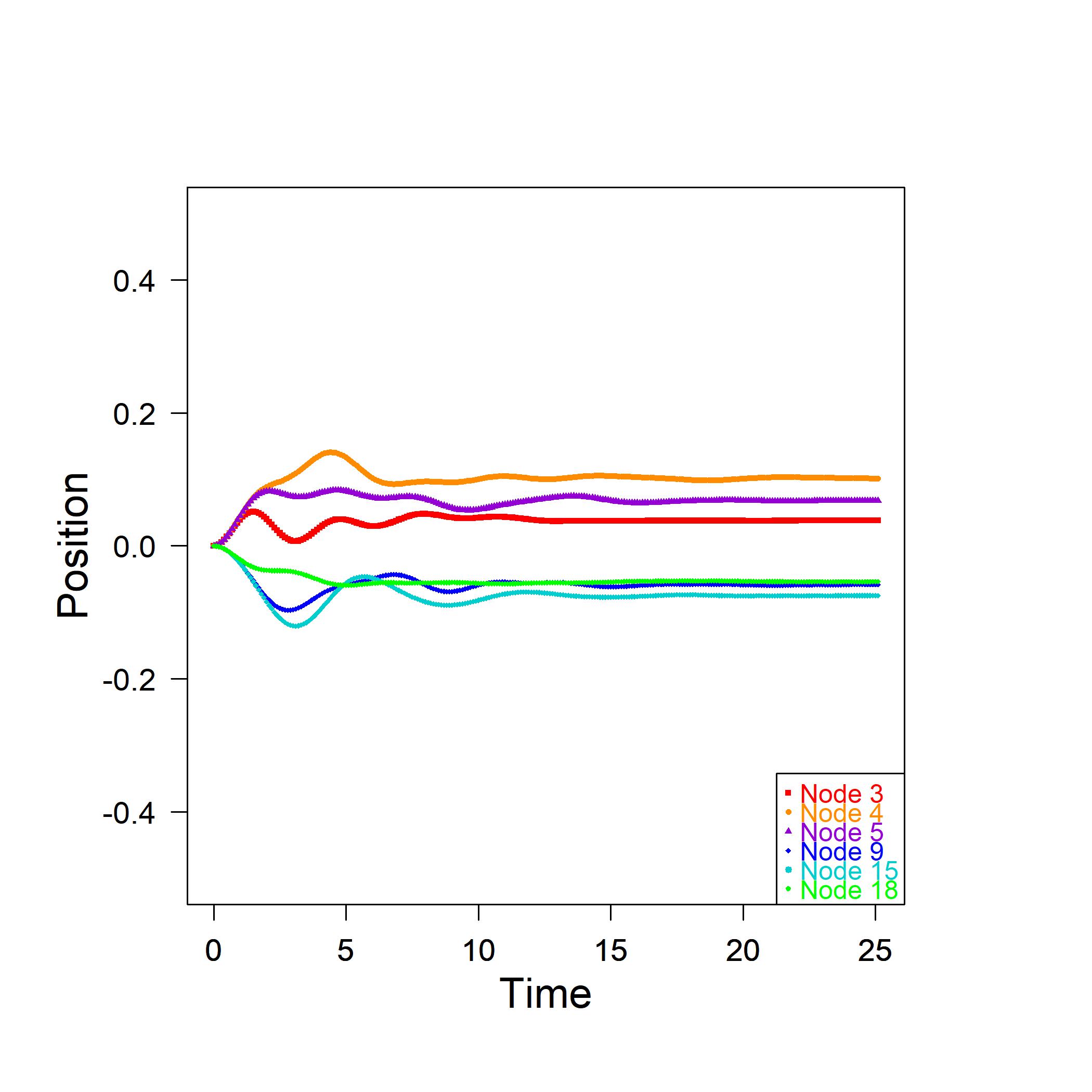}}
	\\
	\subfloat[]{\includegraphics[width=0.36\textwidth]{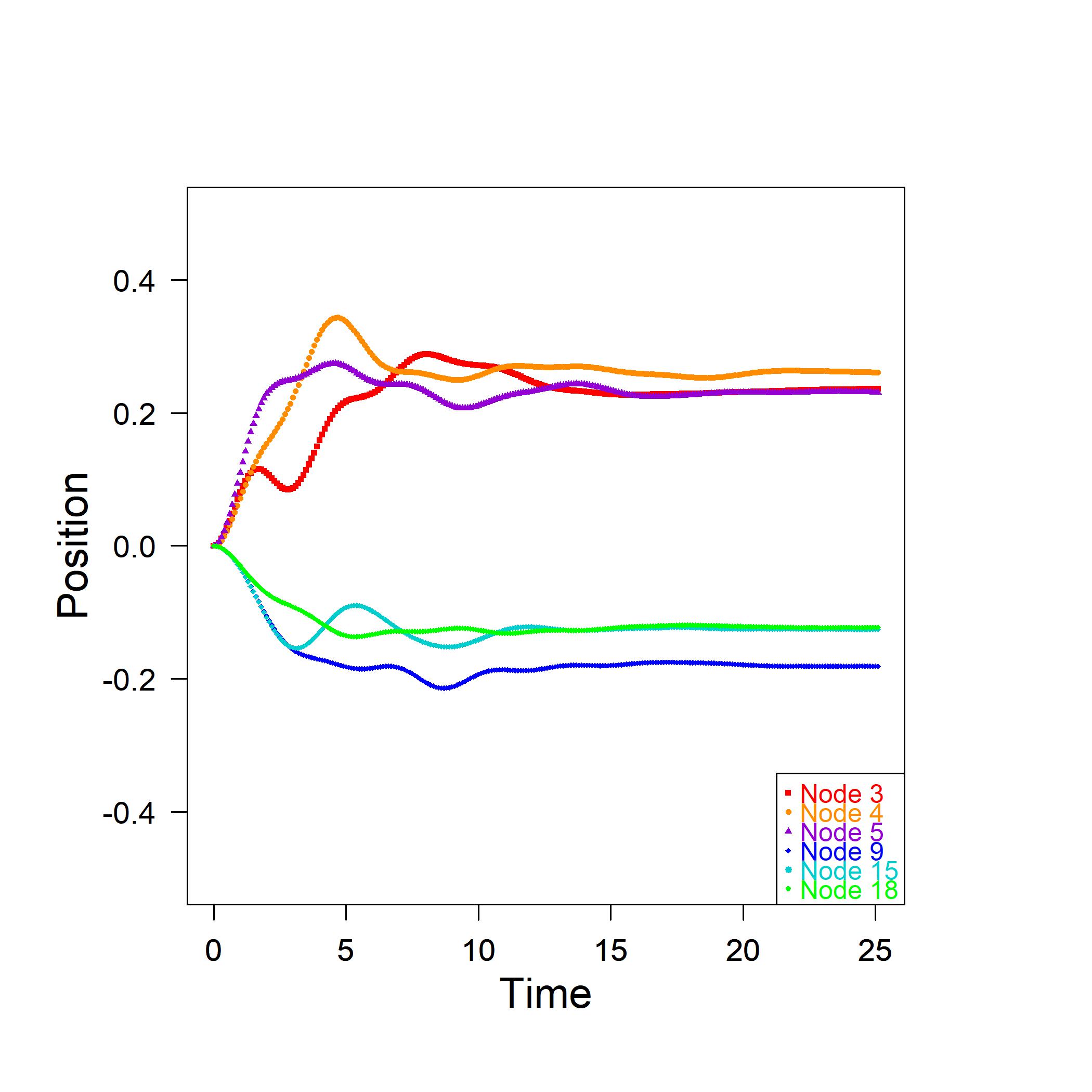}}
	\caption{Dynamic responses of selected nodes in the Syrian power grid under different conditions: panel (a) $\gamma=1$; panel (b) $\gamma=0.4$; panel (c) $\gamma=0.4$. See the text for the different values of the $\bf p$ vector.}
	\label{fig25} 
\end{figure}

We aim primarily to show the qualitative behavior of the solution described in Eq. (\ref{solution_swing_equation}) for specific nodes and for different values of the involved parameters. In figure \ref{fig25} we represent the dynamic response of selected nodes in the Syrian power grid: nodes $3$, $4$ and $5$ correspond to the main power plants in terms of expected capacity, namely Baniyas, Zayzoun and Al-Thawrah plants, whereas nodes $9$, $15$ and $18$ are located in three major Syrian cities, namely Damascus, Palmyra and Al-Hasakah.
In panel (a), we consider a distribution of nominal balanced powers such that ${p_i}=0.13$ for generators and ${p_i}=-0.06$ for loads with damping coefficient $\gamma=1$; in panel (b), we keep the same powers but with damping coefficient $\gamma=0.4$; and, in panel (c), we assume ${\bf p}=[0.004,0.181,0.252,0.202,0.306,0.056]$ for the six generators and ${p_i}=-0.077$ for all the motors. Figure \ref{fig25} qualitatively illustrates some important features of the transient dynamics of $x_{i}(t)$, such as the first peak value, the first peak time and the value at steady state under different conditions (see \citet{Sorrentino2021} for a comprehensive discussion about these values). Table \ref{table5} collects the values of these parameters for two representative nodes in the examined simulation and under conditions corresponding to panels in figure \ref{fig25}.
\begin{table}[H]
		\begin{center}
			\begin{tabular}{||c|c|c|c||}
				\hline
				\multicolumn{4}{||c||}{\textbf{Node} $\bf 4$}\tabularnewline \hline
				{\rm \bf Panel} & {\rm \bf First peak value} & {\rm \bf First peak time} & {\rm \bf Steady state value} \tabularnewline \hline
				(a) & $0.1061903$ & $4.9$ & $0.1020552$ \tabularnewline \hline
				(b) & $0.1409010$ & $4.4$ & $0.1016554$  \tabularnewline \hline
				(c) & $ 0.3436105$ & $4.7$ & $0.2613574$  \tabularnewline \hline
				\multicolumn{4}{||c||}{\textbf{Node} $\bf 15$}\tabularnewline \hline
			{\rm \bf Panel} & {\rm \bf First peak value} & {\rm \bf First peak time} & {\rm \bf Steady state value} \tabularnewline \hline
				(a) & $-0.08873022$ & $3.4$ & $-0.07487655$ \tabularnewline \hline
				(b) & $-0.12053190$ & $3.1$ & $-0.07480745$  \tabularnewline \hline
				(c) & $-0.15393090$ & $3.1$ & $-0.12581260$   \tabularnewline \hline
			\end{tabular}
		\end{center}
		\caption{First peak value, first peak time and steady state value for two representative nodes under conditions corresponding to panels in figure \ref{fig25}.}
		\label{table5}
\end{table}
For instance, as can be observed, the plant corresponding to node $4$ typically has the highest peak value, while node $3$ has the shortest peak time among the examined nodes. For the sake of brevity, we do not consider a power distribution that respects the symmetries of the network, which is such that all nodes within the same symmetry cluster have the same power. In this case, it has been shown that if the vector $\bf p$ respects the symmetries, the response of the complete network is completely encoded in that of the quotient network \cite{SanchezGarcia2020}. Similarly, we do not consider the effects of heterogeneity in the values of the damping coefficients.

\section{Conclusions and future perspectives}

This paper proposes a comprehensive analysis of the potential responses of a complex network when it is host to coupled oscillatory phenomena. It adds several contributions to the existing literature, by providing, for instance, a closed expression for the average synchronization times and the resonance frequencies of the whole network.
Interesting extensions for future works may be the following: (a) a study of the effect of symmetries in the topological structure of the network and a comparison between the response of the full network and of the so-called quotient network under network automorphisms; (b) the analysis of the effect of topologies described by non-normal matrices on the transient dynamics and the relationship between the degree of non-normality and the resilience properties of the networked dynamical system, or, finally, (c) the introduction of non-local interactions through generalizations of the Laplacian and of an external stress factor acting on the network, in order to assess their impact on resonance and synchronization outcomes.

\hfill

{\bf ACKNOWLEDGMENTS}

The author would like to thank the two anonymous reviewers for their helpful advice and insightful suggestions.

\appendix
\section{Mathematical Properties and Polar Decomposition of the $\bf G$ matrices}
\label{Appendix A}
We propose here some theoretical results about the polar decomposition of the matrix $\bf G$ in Eq. (\ref{matrixG1}), with implications in the synchronization process. In particular, we prove that the asymptotic behavior is entirely encoded in the unitary component $\bf U$ of the matrix ${\bf G}$, as stated in Proposition \ref{proposition2} in the main text. Since $\bf G$ is a real non-singular matrix, its polar decomposition can be directly obtained as ${\bf G}={\bf G}({\bf G}^{T}{\bf G})^{-1/2}({\bf G}^{T}{\bf G})^{1/2}$. Therefore, we set
\begin{equation}
{\bf G}={\bf U}{\bf P}
\label{polardecomposition}
\end{equation}
where
\begin{equation}
\label{definitions_UP}
{\bf U}={\bf G}({\bf G}^{T}{\bf G})^{-1/2} \qquad {\rm and}\qquad  {\bf P}=({\bf G}^{T}{\bf G})^{1/2}
\end{equation}
Matrices $\bf U$ and $\bf P$ have the following immediate properties:
\begin{enumerate}
\item[$\bf U$:] \begin{enumerate}
		\item $\bf U$ is a symplectic matrix, ${\bf U}^{T}{\bf J}{\bf U}={\bf J}$; ${\bf U}\in {\rm Sp}(2n)$;
		\item $\bf U$ is an orthogonal matrix, ${\bf U}{\bf U}^{T}={\bf U}^{T}{\bf U}={\bf I}$;
		\item ${\bf U}\in {\rm SO}(2n)$.
	\end{enumerate}
\item[$\bf P$:] \begin{enumerate}
	\item $\bf P$ is a symplectic matrix, ${\bf P}^{T}{\bf J}{\bf P}={\bf J}$;
	\item $\bf P$ is a symmetric matrix, ${\bf P}^{T}={\bf P}$;
	\item ${\bf P}$ is positive definite.
\end{enumerate}
\end{enumerate}

In the following propositions, we provide the expressions for both the eigenvalues and eigenvectors of the two matrices $\bf U$ and $\bf P$. They play an important role in proposition \ref{proposition2}.

\begin{proposition}
	\label{proposition5}
	The eigenvalues of $\bf P$ are the $n$ couples of reciprocal values given by
	\begin{equation}
	\lambda^{\rm P\pm}_i=\frac{1}{2}\left[\sqrt{\mu_i^2+4}\pm \mu_i\right]
	\label{eigenvalues_P}
	\end{equation}
and the corresponding normalized eigenvectors are given by
	\begin{equation}
	\psi^{\rm P \pm}_i=
	\frac{1}{\sqrt{1+{\lambda^{\rm P\pm}_i}^{2}}}
	\left[ \begin{array}{c} 
	\phi_i \\
	\pm \lambda^{\rm P\pm}_i\phi_i \\
	\end{array} \right]
	\label{eigenvectors_P}
	\end{equation}
\end{proposition}

\begin{remark}
All the eigenvalues ${\lambda^{P \pm}_i}$ are real and positive: ${\lambda^{P \pm}_i}\in {\mathbb R}\ {and}\ {\lambda^{P \pm}_i}>0, \forall i=1,\dots, n$. They are in couples of reciprocal values: $1/\lambda^{\rm P-}_i=\lambda^{\rm P+}_i$. They satisfy $\lambda^{\rm P+}_i>1$ and $0<\lambda^{\rm P-}_i<1$. Finally, for $\mu_n=0$, ${\lambda^{P \pm}_n}=1$ and the corresponding eigenspace has dimension $2$. 
\end{remark}

Let us analyze now the operator $\bf U$ defined in Eq. (\ref{definitions_UP}). The first step is to prove that it can be given the form stated in the following proposition.

\begin{proposition}
\label{proposition6}
The unitary operator $\bf U$ can be expressed as
\begin{equation}
\label{U}
{\bf U}=
\left[ \begin{array}{cc} 
{\bf \cal A} & {\bf \cal B} \\
-{\bf \cal B} & {\bf \cal A} \\
\end{array} \right]
\end{equation}
where ${\bf \cal A}$ and ${\bf \cal B}$ are the $n-$square matrices
\begin{equation}
\label{AB}
{\bf \cal A}=\sum_{i=1}^{n}\frac{1-{\lambda^{\rm P+}_i}^{2}}{1+{\lambda^{\rm P+}_i}^{2}}{\phi_i}{\phi_i}^{T}\ {\it and}\ {\bf \cal B}=\sum_{i=1}^{n}\frac{2{\lambda^{\rm P+}_i}}{1+{\lambda^{\rm P+}_i}^{2}}{\phi_i}{\phi_i}^{T}
\end{equation}
\end{proposition}

Coefficients in formulae (\ref{AB}) suggest introducing suitable angles which will play a significant role in the interpretation of matrices ${\bf \cal A}$ and ${\bf \cal B}$.

Let us define angles $\theta_{i}$ as
\begin{equation}
\label{angles}
{\lambda^{\rm P +}_i}=\tan\left(\frac{\theta_i}{2}\right), \ i=1,\dots, n
\end{equation}

By introducing angles $\theta_i$, matrices ${\bf \cal A}$ and ${\bf \cal B}$ in (\ref{AB}) can be expressed as
\begin{equation}
\label{ABangles}
{\bf \cal A}=\sum_{i=1}^{n}\cos(\theta_i){\phi_i}{\phi_i}^{T}\quad {\rm and}\quad {\bf \cal B}=\sum_{i=1}^{n}\sin(\theta_i){\phi_i}{\phi_i}^{T}
\end{equation}

so that $\bf U$ becomes
\begin{equation}
\label{Uangles}
{\bf U}=
\left[ \begin{array}{cc} 
\sum_{i=1}^{n}\cos(\theta_i){\phi_i}{\phi_i}^{T} & \sum_{i=1}^{n}\sin(\theta_i){\phi_i}{\phi_i}^{T} \\
-\sum_{i=1}^{n}\sin(\theta_i){\phi_i}{\phi_i}^{T} & \sum_{i=1}^{n}\cos(\theta_i){\phi_i}{\phi_i}^{T} \\
\end{array} \right].
\end{equation}

Expression (\ref{Uangles}) makes it clear the ${\rm SO}(2n)$-nature of the operator $\bf U$.
Now we can focus on the eigenvalues and eigenvectors of $\bf U$, in the next proposition.

\begin{proposition}
	\label{proposition7}
The eigenvalues of $\bf U$ are the $n$ couples
	\begin{equation}
	\lambda^{\rm U\pm}_j=\cos\theta_j\pm i \sin\theta_j=e^{\pm i \theta_j}
	\label{eigenvalues_U}
	\end{equation}
where $\theta_j=2{\arctan}\, {\lambda^{\rm P +}_j}$ and the corresponding normalized eigenvectors are given by
	\begin{equation}
	\psi^{\rm U \pm}_j=
	\frac{1}{\sqrt{2}}
	\left[ \begin{array}{c} 
	\phi_j \\
	\pm i\phi_j \\
	\end{array} \right]
	\label{eigenvectors_U}
	\end{equation}
\end{proposition}

\begin{remark}
The angles ${\theta_j}=2\arctan {\lambda^{\rm P +}_j}$, introduced in formula (\ref{angles}), are the arguments of the unitary complex eigenvalues of $\bf U$. More precisely, ${\lambda^{\rm P +}_j}>1$ determine the angles $\frac{\pi}{2}<{\theta_j}<\pi$ of $n$ eigenvalues of $\bf U$ and $0<{\lambda^{\rm P -}_j}<1$, which are reciprocals of ${\lambda^{\rm P +}_j}$, determine the angles $\pi<{\theta_j}<\frac{3}{2}\pi$ of the other $n$ eigenvalues of $\bf U$. In particular, $\cos\theta_j<0$, $\forall j=1,\dots, n-1$. Finally, $\cos\theta_n=0$ and the two $n-$th eigenvalues are  $\lambda^{\rm U+}_n=i$ and $\lambda^{\rm U-}_n=-i$.
\end{remark}

\begin{remark}
The unitary operator $\bf U$ takes the standard form
\begin{equation}
{\bf U}=\sum_{j=1}^{n}\left( e^{i\theta_{j}}\psi^{\rm U +}_j\psi^{\rm U + \dagger}_j+e^{-i\theta_{j}}\psi^{\rm U -}_j\psi^{\rm U - \dagger}_j \right)
\end{equation}
where $\psi^{\rm U \pm \dagger}_j$ denotes the adjoint vector, which is the transpose of the complex conjugate vector.
\end{remark}

\section{Proofs of the propositions}
\label{Appendix B}

\subsection{Proof of Proposition \ref{proposition1}}

\begin{proof}
The equation ${\bf G}{\psi^{\rm G}_i}=\lambda^{\rm G}_i{\psi^{\rm G}_i}$ is equivalent to
\begin{equation*}
\left[ \begin{array}{cc} 
{\bf 0} & {\bf I} \\
- {\bf I}   & -{\bf L} \\
\end{array} \right]
\left(\begin{array}{c}
{\bf x}_i \\
{\bf v}_i \\
\end{array} \right)=
\lambda^{\rm G}_i
\left(\begin{array}{c}
{\bf x}_i \\
{\bf v}_i \\
\end{array} \right)
\  \Longrightarrow \
\left\{ 
\begin{array}{l}
{\bf v}_i=\lambda^{\rm G}_i {\bf x}_i \\
-{\bf x}_i-{\bf L}{\bf v}_i=\lambda^{\rm G}_i{\bf v}_i\\
\end{array}
\right.
\end{equation*}
which brings to
\begin{equation*}
{\bf L}{\bf x}_i=-\frac{1+{\lambda^{\rm G}_i}^{2}}{\lambda^{\rm G}_i}{\bf x}_i
\end{equation*}
Then ${\bf x}_i$ are eigenvectors of the Laplacian $\bf L$ and the eigenvalues $\lambda^{\rm G}_i$ in Eq. (\ref{eigenvalues_G}) are the roots of the equation ${\lambda^{\rm G}_i}^{2}+\mu_i\lambda^{\rm G}_i+1=0$. It is straightforward to show that the eigenvalues $\lambda^{\rm G \pm}_i$ are reciprocal: $1/\lambda^{\rm G-} _i=\lambda^{\rm G+}_i$. Moreover
\begin{equation*}
\begin{split}
&\left[ \begin{array}{cc} 
{\bf 0} & {\bf I} \\
- {\bf I}   & -{\bf L} \\
\end{array} \right]
\psi_i=
\left[ \begin{array}{cc} 
{\bf 0} & {\bf I} \\
- {\bf I}   & -{\bf L} \\
\end{array} \right]
\left(\begin{array}{c}
\phi_i \\
\lambda^{\rm G}_i \phi_i\\
\end{array} \right)=\\
&\left(\begin{array}{c}
\lambda^{\rm G}_i \phi_i\\
-\phi_i-\lambda^{\rm G}_i {\bf L} \phi_i\\
\end{array} \right)=
\left(\begin{array}{c}
\lambda^{\rm G}_i \phi_i\\
-\phi_i-\lambda^{\rm G}_i \mu_i \phi_i\\
\end{array} \right)\\
&=
\left(\begin{array}{c}
\lambda^{\rm G}_i \phi_i\\
-\left(1+\mu_i\lambda^{\rm G}_i\right)\phi_i\\
\end{array} \right)=
\left(\begin{array}{c}
\lambda^{\rm G}_i \phi_i\\
{\lambda^{\rm G}_i}^2\phi_i\\
\end{array} \right)\\
&=
{\lambda^{\rm G}_i}
\left(\begin{array}{c}
\phi_i\\
{\lambda^{\rm G}_i}\phi_i\\
\end{array} \right)=
{\lambda^{\rm G}_i}\psi_i
\end{split}
\end{equation*}
Normalization follows trivially.
\end{proof}

\subsection{Proof of Proposition \ref{proposition2}}

\begin{proof}
Since ${\bf U}={\bm \Psi}_{\rm U}{\bm \Lambda}_{\rm U}{\bm \Psi}_{\rm U}^{\dagger}$ and $e^{t{\bf U}}={\bm \Psi}_{\rm U}e^{t{\bm  \Lambda}_{\rm U}}{\bm \Psi}_{\rm U}^{\dagger}$ and, according to the remark before the proposition statement, ${\rm Re}(\lambda^{\rm U+}_j)=\cos\theta_j<0$ for $\forall j=1,\dots, n-1$; and ${\rm Re}(\lambda^{\rm U+}_n)=\cos\theta_n=0$, we have, for $t\to +\infty$:
\begin{equation*}
\begin{split}
e^{t{\bf U}}
&\sim
\left[ \begin{array}{ccc} 
\frac{\phi_n}{\sqrt{2}} &
\frac{\phi_n}{\sqrt{2}} &
\dots \\
\frac{i\phi_n}{\sqrt{2}} &
\frac{-i\phi_n}{\sqrt{2}} &
\dots \\
\end{array} \right]
\left( \begin{array}{ccccc} 
e^{it} & & & & \\
& e^{-it} & & & \\
& & 0 & & \\
& & & \ddots &  \\
& & & & 0 \\
\end{array} \right)
\left[ \begin{array}{cc} 
\frac{\phi_{n}^{T}}{\sqrt{2}} &
\frac{-i\phi_{n}^{T}}{\sqrt{2}}\\
\frac{\phi_{n}^{T}}{\sqrt{2}} &
\frac{+i\phi_{n}^{T}}{\sqrt{2}}\\
\vdots &
\vdots\\
\end{array} \right]
\\
&=\frac{1}{2} 
\left[ \begin{array}{ccc} 
{\phi_n} &
{\phi_n} &
\dots \\
{i\phi_n} &
{-i\phi_n} &
\dots \\
\end{array} \right]
\left[ \begin{array}{cc} 
e^{it}\phi_{n}^{T} &
-ie^{it}\phi_{n}^{T}\\
e^{-it}\phi_{n}^{T} &
+ie^{-it}\phi_{n}^{T}\\
0 &
0\\
\vdots &
\vdots\\
0 &
0\\
\end{array} \right]
\\
&=\frac{1}{2}
\left[ \begin{array}{cc} 
e^{it}\phi_{n}\phi_{n}^{T}+e^{-it}\phi_{n}\phi_{n}^{T} &
-ie^{it}\phi_{n}\phi_{n}^{T}+ie^{-it}\phi_{n}\phi_{n}^{T}\\
ie^{it}\phi_{n}\phi_{n}^{T}-ie^{-it}\phi_{n}\phi_{n}^{T} &
e^{it}\phi_{n}\phi_{n}^{T}+e^{-it}\phi_{n}\phi_{n}^{T}\\
\end{array} \right]
\\
&=
\left[ \begin{array}{cc} 
\frac{e^{it}+e^{-it}}{2}\phi_{n}\phi_{n}^{T} &
\frac{e^{it}-e^{-it}}{2i}\phi_{n}\phi_{n}^{T} \\
-\frac{e^{it}-e^{-it}}{2i}\phi_{n}\phi_{n}^{T}  &
\frac{e^{it}+e^{-it}}{2}\phi_{n}\phi_{n}^{T} \\
\end{array} \right]\\
&=
\left[ \begin{array}{cc} 
\cos t\, \phi_{n}\phi_{n}^{T} &
\sin t\, \phi_{n}\phi_{n}^{T} \\
-\sin t\, \phi_{n}\phi_{n}^{T}  &
\cos t\, \phi_{n}\phi_{n}^{T} \\
\end{array} \right]\\
&=\frac{1}{n}
\left[ \begin{array}{cc} 
\cos t\, {\bf 1} &
\sin t\, {\bf 1} \\
-\sin t\, {\bf 1} &
\cos t\, {\bf 1} \\
\end{array} \right]\\
\end{split}
\end{equation*}

Similarly, for the operator $\bf G$. In fact, since $\psi^{\rm G \pm}_i$ is given by Eq. (\ref{eigenvectors_G}) then, for $\mu_n =0$,
\begin{equation*}
\psi^{\rm G \pm}_n=\frac{1}{\sqrt 2}\left( \begin{array}{c} 
\phi_n \\
\pm i\phi_n \\
\end{array} \right)
=\psi^{\rm U \pm}_n
\end{equation*}

Now, let ${\bm \Psi}_{\rm G}=\left[\psi^{\rm G +}_n, \psi^{\rm G -}_n, \dots \right]$. 
Since ${\bm \Psi}_{\rm G}^{-1}{\bm \Psi}_{\rm G}={\bf I}$ and $\left({\psi^{\rm G \pm}_n}\right)^{\dagger}{\psi^{\rm G \pm}_n}=1$ and $\left({\psi^{\rm G \pm}_n}\right)^{\dagger}{\psi^{\rm G \mp}_n}=0$, we have
\begin{equation*}
{\bm \Psi}_{\rm G}^{-1}=\left[ \begin{array}{c} 
{\psi^{\rm G +}_n}^{\dagger} \\
{\psi^{\rm G -}_n}^{\dagger}\\
\vdots\\
\end{array} \right]
\end{equation*}

and, through the same steps as before, for $t\to +\infty$:
\begin{equation*}
\begin{split}
e^{t{\bf G}}
&={\bm \Psi}_{\rm G}\, e^{t{\bm \Lambda}_{G}}\, {\bm \Psi}_{\rm G}^{-1}\\
&\sim
\left[ \begin{array}{ccc} 
\psi^{\rm G +}_n &
\psi^{\rm G -}_n &
\dots \\
\end{array} \right]
\left( \begin{array}{ccccc} 
e^{it} & & & & \\
& e^{-it} & & & \\
& & 0 & & \\
& & & \ddots &  \\
& & & & 0 \\
\end{array} \right)
\left[ \begin{array}{c} 
{\psi^{\rm G +}_n}^{\dagger} \\
{\psi^{\rm G -}_n}^{\dagger}\\
\vdots\\
\end{array} \right]
\\
&=\frac{1}{n}
\left[ \begin{array}{cc} 
\cos t\, {\bf 1} &
\sin t\, {\bf 1} \\
-\sin t\, {\bf 1} &
\cos t\, {\bf 1} \\
\end{array} \right]\\
\end{split}
\end{equation*}
\end{proof}

\subsection{Proof of Proposition \ref{proposition3}}

\begin{proof}
The integrand in ${\bf y}(t)=e^{{\bf G}t}\left[ \int_{0}^{t}e^{-{\bf G}u}{\bf b}(u)du +{\bf c} \right]$ takes the form
\begin{equation}
e^{-{\bf G}u}{\bf b}(u)=f(u)\left[ \begin{array}{c} 
-(\sqrt{\bf H})^{-1}\sin\sqrt{\bf H}u\, {\bf e}_{h} \\
\cos\sqrt{\bf H}u \, {\bf e}_{h}\\
\end{array} \right]\\
\label{integrand}
\end{equation}

Let us observe that ${\bf H}={\bf I}+{\bf L}$ is nonsingular with eigenvalues $1+\mu_{i},\ i=1,\dots, n$, so that the two components in Eq. (\ref{integrand}) can be re-written as
\begin{equation*}
\begin{split}
&-(\sqrt{\bf H})^{-1}\sin\sqrt{\bf H}u\, {\bf e}_{h}=\sum_{i=1}^{n}\lambda_{i}^{(1)}\phi_{i}(h)\phi_{i}\\
&
{\rm with}
\quad
\lambda_{i}^{(1)}=-\frac{1}{\sqrt{1+\mu_{i}}}\sin \left( \sqrt{1+\mu_{i}}\, u \right)
\end{split}
\end{equation*}

and 
\begin{equation*}
\begin{split}
&\cos\sqrt{\bf H}u\, {\bf e}_{h}=\sum_{i=1}^{n}\lambda_{i}^{(2)}\phi_{i}(h)\phi_{i}\\
&{\rm with}
\quad
\lambda_{i}^{(2)}(u)=\cos \left( \sqrt{1+\mu_{i}}\, u \right)
\end{split}
\end{equation*}

For both components, the integral in Eq. (\ref{generalsolution4}) becomes
\begin{equation}
\int_{0}^{t}e^{-{\bf G}u}{\bf b}(u)du= \sum_{i=1}^{n}\left( \int_{0}^{t}f(u)\lambda_{i}^{(1,2)}(u)du\right) \phi_{i}(h)\phi_{i}
\label{integrals}
\end{equation}

We compute separately the two integrals in Eq. (\ref{integrals}):
\begin{equation*}
\begin{split}
&\int_{0}^{t}f(u)\lambda_{i}^{(1)}(u)du=\\
&\frac{F_{0}}{2\sqrt{1+\mu_{i}}}\left[ \frac{\sin\left( \sqrt{1+\mu_{i}}+\omega \right)t}{\sqrt{1+\mu_{i}}+\omega} - \frac{\sin\left( \sqrt{1+\mu_{i}}-\omega \right)t}{\sqrt{1+\mu_{i}}-\omega} \right]
\end{split}
\end{equation*}

and 
\begin{equation*}
\begin{split}
&\int_{0}^{t}f(u)\lambda_{i}^{(2)}(u)du=\\
&\frac{F_{0}}{2}\left[
\frac{2\omega}{\omega^2-(1+\mu_{i})}-
\frac{\cos\left(\omega+ \sqrt{1+\mu_{i}} \right)t}{\omega +\sqrt{1+\mu_{i}}}-
\frac{\cos\left(\omega- \sqrt{1+\mu_{i}} \right)t}{\omega-\sqrt{1+\mu_{i}}} \right]
\end{split}
\end{equation*}

Then Eq. (\ref{integrals}) gives
\begin{equation*}
\int_{0}^{t}e^{-{\bf G}u}{\bf b}(u)du=
\frac{F_{0}}{2} \sum_{i=1}^{n}
\left[ \begin{array}{c} 
a_{i}(t) \phi_{i}(h)\phi_{i} \\
b_{i}(t) \phi_{i}(h)\phi_{i}\\
\end{array} \right]\\
\end{equation*}

with
\begin{equation*}
\begin{split}
&a_{i}(t)=\\
&\frac{1}{\sqrt{1+\mu_{i}}}
\left[ \frac{\sin\left( \sqrt{1+\mu_{i}}+\omega \right)t}{\sqrt{1+\mu_{i}}+\omega} - \frac{\sin\left( \sqrt{1+\mu_{i}}-\omega \right)t}{\sqrt{1+\mu_{i}}-\omega} \right]
\end{split}
\end{equation*}

and
\begin{equation*}
\begin{split}
&b_{i}(t)=\\
&\left[
\frac{2\omega}{\omega^2-(1+\mu_{i})}-
\frac{\cos\left(\omega+ \sqrt{1+\mu_{i}} \right)t}{\omega +\sqrt{1+\mu_{i}}}-
\frac{\cos\left(\omega- \sqrt{1+\mu_{i}} \right)t}{\omega-\sqrt{1+\mu_{i}}} \right]
\end{split}
\end{equation*}

Finally, Eq. (\ref{generalsolution4}) becomes
\begin{equation*}
\begin{split}
&{\bf y}(t)=\\
&\frac{F_{0}}{2} \sum_{i=1}^{n}
\left[ \begin{array}{c} 
a_{i}(t) \phi_{i}(h)\cos \sqrt{\bf{H}}\, t \phi_{i}+b_{i}(t)\phi_{i}(h)(\sqrt{\bf{H}})^{-1}\sin \sqrt{\bf{H}}\, t \phi_{i} \\
-a_{i}(t) \phi_{i}(h) \sqrt{\bf{H}}\sin \sqrt{\bf{H}}\, t \phi_{i}+b_{i}(t)\phi_{i}(h)\cos \sqrt{\bf{H}}\, t \phi_{i}\\
\end{array} \right]
\end{split}
\end{equation*}

and, after some algebraic manipulations, it can be expressed as
\begin{equation*}
{\bf y}(t)=F_{0} \sum_{i=1}^{n}
\frac{\omega}{\omega_{i}^{2}-\omega^2}\phi_{i}(h)
\left[\begin{array}{c} 
\left(\frac{\sin \omega t}{\omega}- \frac{\sin \omega_{i} t}{\omega_{i}} \right) \phi_{i} \\
\left(\cos \omega t- \cos \omega_{i} t \right) \phi_{i}\\
\end{array} \right]\\
\end{equation*}

with $\omega_{i}\coloneqq\sqrt{1+\mu_{i}}$.
\end{proof}

\subsection{Proof of Proposition \ref{proposition4}}

\begin{proof}
In this case we have
\begin{equation*}
\begin{split}
&e^{-{\bf G}u}{\bf b}(u)=\\
&f(u)\left[ \begin{array}{c} 
{\bf B}^{-1}\left(e^{-\frac{{\bf B}-{\bf L}}{2}u}-
e^{\frac{{\bf B}+{\bf L}}{2}u}
\right)\, {\bf e}_{h} \\
\frac{1}{2}{\bf B}^{-1}\left(
({\bf B}-{\bf L})e^{-\frac{{\bf B}-{\bf L}}{2}u}+
({\bf B}+{\bf L})e^{\frac{{\bf B}+{\bf L}}{2}u}
\right) \, {\bf e}_{h}\\
\end{array} \right]
\end{split}
\end{equation*}

Let us consider the two components separately. The first one is
\begin{equation*}
{\bf B}^{-1}\left(e^{-\frac{{\bf B}-{\bf L}}{2}u}-
e^{\frac{{\bf B}+{\bf L}}{2}u}
\right)\, {\bf e}_{h}=
\sum_{i=1}^{n}\lambda_{i}^{(1)}(u)\phi_{i}^{\star}(h)\phi_{i}
\end{equation*}

where 
\begin{equation}
\lambda_{i}^{(1)}(u)=\frac{1}{\sqrt{\mu_{i}^2-4}}
\left[ e^{-\frac{1}{2}(\sqrt{\mu_{i}^{2}-4}-\mu_{i})u}-e^{\frac{1}{2}(\sqrt{\mu_{i}^{2}-4}+\mu_{i})u}\right]
\label{eigenvalues_Ui}
\end{equation}

Let us observe that, by using Eq. (\ref{eigenvalues_G}) for the eigenvalues $\lambda^{\rm G\pm}_i$, Eq. (\ref{eigenvalues_Ui}) can be simplified as
\begin{equation*}
\lambda_{i}^{(1)}(u)=\frac{e^{-\lambda_{i}^{G+}u}-e^{-\lambda_{i}^{G-}u}}{\sqrt{\mu_{i}^2-4}}
\end{equation*}

Similarly, the second one is
\begin{equation*}
\begin{split}
&\frac{1}{2}{\bf B}^{-1}\left(
({\bf B}-{\bf L})e^{-\frac{{\bf B}-{\bf L}}{2}u}+
({\bf B}+{\bf L})e^{\frac{{\bf B}+{\bf L}}{2}u}
\right) \, {\bf e}_{h}\\
&=\sum_{i=1}^{n}\lambda_{i}^{(2)}(u)\phi_{i}^{\star}(h)\phi_{i}
\end{split}
\end{equation*}

where 
\begin{equation*}
\lambda_{i}^{(2)}(u)=\frac{\lambda_{i}^{G+}e^{-\lambda_{i}^{G+}u}-\lambda_{i}^{G-}e^{-\lambda_{i}^{G-}u}}{\sqrt{\mu_{i}^2-4}}
\end{equation*}

For both components, the integral in the general solution has the form
\begin{equation*}
\int_{0}^{t}e^{-{\bf G}u}{\bf b}(u)du= \sum_{i=1}^{n}\left( \int_{0}^{t}f(u)\lambda_{i}^{(1,2)}(u)du\right) \phi_{i}^{\star}(h)\phi_{i}
\end{equation*}

The integral of the first component is
\begin{align}
\footnotesize
\begin{split}
&\int_{0}^{t}f(u)\lambda_{i}^{(1)}(u)du\\
&=\frac{F_{0}}{\sqrt{\mu_{i}^2-4}}
\int_{0}^{t}\sin \omega u \left( e^{-\lambda_{i}^{G+}u}-e^{-\lambda_{i}^{G-}u} \right) du\\
&= \frac{F_{0}}{\sqrt{\mu_{i}^2-4}}\Bigg[
\bigg(-\frac{\lambda_{i}^{G+}}{\omega^2+(\lambda_{i}^{G+})^{2}}e^{-\lambda_{i}^{G+}t}+\frac{\lambda_{i}^{G-}}{\omega^2+(\lambda_{i}^{G-})^{2}}e^{-\lambda_{i}^{G-}t}\Big)\sin \omega t+\\
&-\Big(\frac{1}{\omega^2+(\lambda_{i}^{G+})^{2}}e^{-\lambda_{i}^{G+}t}-\frac{1}{\omega^2+(\lambda_{i}^{G-})^{2}}e^{-\lambda_{i}^{G-}t}\Big)\omega \cos \omega t+\\
&+\Big(\frac{1}{\omega^2+(\lambda_{i}^{G+})^{2}}-\frac{1}{\omega^2+(\lambda_{i}^{G-})^{2}}\bigg) \omega \Bigg]=\\
&=\frac{F_{0}}{\sqrt{\mu_{i}^2-4}}\Bigg[
\bigg(-\frac{\lambda_{i}^{G+}}{\omega^2+(\lambda_{i}^{G+})^{2}}e^{-\lambda_{i}^{G+}t}+\frac{\lambda_{i}^{G-}}{\omega^2+(\lambda_{i}^{G-})^{2}}e^{-\lambda_{i}^{G-}t}\Big)\sin \omega t+\\
&-\Big(\frac{1}{\omega^2+(\lambda_{i}^{G+})^{2}}e^{-\lambda_{i}^{G+}t}-\frac{1}{\omega^2+(\lambda_{i}^{G-})^{2}}e^{-\lambda_{i}^{G-}t}\Big)\omega \cos \omega t\Bigg] +\\
&+\frac{F_{0}\omega \mu_{i}}{\left(\omega^2+(\lambda_{i}^{G+})^{2}\right)\left(\omega^2+(\lambda_{i}^{G-})^{2}\right)}\coloneqq a_{i}\\
\end{split}
\label{a_i}
\end{align}

The integral of the second component is
\begin{equation}
\footnotesize
\begin{split}
&\int_{0}^{t}f(u)\lambda_{i}^{(2)}(u)du \\
&=\frac{F_{0}}{\sqrt{\mu_{i}^2-4}}
\int_{0}^{t}\sin \omega u \left( \lambda_{i}^{G+} e^{-\lambda_{i}^{G+}u}-\lambda_{i}^{G-} e^{-\lambda_{i}^{G-}u} \right) du\\
&=\frac{F_{0}}{\sqrt{\mu_{i}^2-4}}\Bigg[
\bigg(-\frac{(\lambda_{i}^{G+})^2}{\omega^2+(\lambda_{i}^{G+})^{2}}e^{-\lambda_{i}^{G+}t}+\frac{(\lambda_{i}^{G-})^2}{\omega^2+(\lambda_{i}^{G-})^{2}}e^{-\lambda_{i}^{G-}t}\Big)\sin \omega t+\\
&-\Big(\frac{\lambda_{i}^{G+}}{\omega^2+(\lambda_{i}^{G+})^{2}}e^{-\lambda_{i}^{G+}t}-\frac{\lambda_{i}^{G-}}{\omega^2+(\lambda_{i}^{G-})^{2}}e^{-\lambda_{i}^{G-}t}\Big)\omega \cos \omega t\Bigg] +\\
&+\frac{F_{0}\omega (\omega^2 -1)}{\left(\omega^2+(\lambda_{i}^{G+})^{2}\right)\left(\omega^2+(\lambda_{i}^{G-})^{2}\right)}\coloneqq b_{i}\\
\end{split}
\label{b_i}
\end{equation}

Finally, by setting $\xi_{i}=\sqrt{\mu_{i}^2-4}$:

\begin{widetext}
\begin{equation*}
{\bf y}(t)=e^{{\bf G}t}\left[ \int_{0}^{t}e^{-{\bf G}u}{\bf b}(u)du  \right]
=\sum_{i=1}^{n} \frac{\phi_{i}^{\star}(h)}{2\xi_{i}} 
\left[ \begin{array}{c} 
\left\{ a_{i} \left[ (\xi_{i}+\mu_{i})e^{\frac{\xi_{i}-\mu_{i}}{2}t}+(\xi_{i}-\mu_{i})e^{-\frac{\xi_{i}+\mu_{i}}{2}t} \right]+
2b_{i}\left[ e^{\frac{\xi_{i}-\mu_{i}}{2}t}-e^{-\frac{\xi_{i}+\mu_{i}}{2}t} \right] \right\} \phi_{i} \\
\left\{ -2 a_{i} \left[ e^{\frac{\xi_{i}-\mu_{i}}{2}t}-
e^{-\frac{\xi_{i}+\mu_{i}}{2}t} \right]+
b_{i}\left[ (\xi_{i}-\mu_{i}) e^{\frac{\xi_{i}-\mu_{i}}{2}t}+(\xi_{i}+\mu_{i})e^{-\frac{\xi_{i}+\mu_{i}}{2}t} \right] \right\} \phi_{i}\\
\end{array} \right]
\end{equation*}
	
By observing that $\xi_{i}-\mu_{i}=2\lambda_{i}^{G+}$ and $\xi_{i}+\mu_{i}=-2\lambda_{i}^{G-}$, this solution can be re-written as
\begin{equation*}
{\bf y}(t)=
\left[ \begin{array}{c} 
{\bf x}(t) \\
{\bf v}(t) \\
\end{array} \right]=
\sum_{i=1}^{n} \frac{\phi_{i}^{\star}(h)}{\xi_{i}} 
\left[ \begin{array}{c} 
\big[ \left(\lambda_{i}^{G+}a_{i}-b_{i}\right) e^{\lambda_{i}^{G-}t}+ \left(b_{i}- \lambda_{i}^{G-}a_{i} \right) e^{\lambda_{i}^{G+}t} \big] \phi_{i} \\
\big[ \left(a_{i}-\lambda_{i}^{G-}b_{i}\right) e^{\lambda_{i}^{G-}t}+ \left(\lambda_{i}^{G+}b_{i}-a_{i}\right) e^{\lambda_{i}^{G+}t} \big] \phi_{i}\\
\end{array} \right]
\end{equation*}
\end{widetext}

\hfill\eject

\hfill\eject

Finally, by substituting the explicit expressions for $a_{i}=a_{i}(t)$ and $b_{i}=b_{i}(t)$ in Eqs. (\ref{a_i}) and (\ref{b_i}), we get the final expressions for ${\bf x}(t)$ and ${\bf v}(t)$:
\begin{equation*}
\footnotesize
\begin{split}
{\bf x}(t)&=
\sum_{i=1}^{n} \frac{F_{0}\phi_{i}^{\star}(h)}{\xi_{i}} 
\Bigg[ \left( 
\frac{1-(\lambda_{i}^{-})^2}{\xi_{i}(\omega^2+(\lambda_{i}^{-})^2)}+
\frac{1-(\lambda_{i}^{+})^2}{\xi_{i}(\omega^2+(\lambda_{i}^{+})^2}
\right) \sin \omega t \\
&+ \left( 
\frac{1}{\omega^2+(\lambda_{i}^{-})^2}-
\frac{1}{\omega^2+(\lambda_{i}^{+})^2}
\right) \omega \cos \omega t \\
&+
\left( 
\frac{\omega}{(\omega^2+(\lambda_{i}^{+})^2}e^{\lambda_{i}^{+} t}-\frac{\omega}{(\omega^2+(\lambda_{i}^{-})^2)}e^{\lambda_{i}^{-} t}
\right)\Bigg]\phi_{i}\\
\end{split}
\end{equation*}

and 
\begin{equation*}
\footnotesize
\begin{split}
{\bf v}(t)&=
\sum_{i=1}^{n} \frac{F_{0}\phi_{i}^{\star}(h)}{\xi_{i}} 
\Bigg[ \left( 
\frac{\lambda_{i}^{-}(1-(\lambda_{i}^{-})^2)}{\xi_{i}(\omega^2+(\lambda_{i}^{-})^2)}+
\frac{\lambda_{i}^{+}(1-(\lambda_{i}^{+})^2)}{\xi_{i}(\omega^2+(\lambda_{i}^{+})^2}
\right) \sin \omega t \\
&+ \left( 
\frac{\lambda_{i}^{-}}{\omega^2+(\lambda_{i}^{-})^2}-
\frac{\lambda_{i}^{+}}{\omega^2+(\lambda_{i}^{+})^2}
\right) \omega \cos \omega t \\
&+
\left( 
\frac{\omega\lambda_{i}^{+}}{(\omega^2+(\lambda_{i}^{+})^2}e^{\lambda_{i}^{+} t}-\frac{\omega\lambda_{i}^{-}}{(\omega^2+(\lambda_{i}^{-})^2)}e^{\lambda_{i}^{-} t}
\right)\Bigg]\phi_{i}\\
\end{split}
\end{equation*}
\end{proof}

\subsection{Proof of Proposition \ref{proposition5}}                       
\begin{proof}
Since ${\bf P}=({\bf G}^{T}{\bf G})^{1/2}$, we look for the eigenvalues of ${\bf G}^{T}{\bf G}$ which can be expressed as\footnote{Note that  ${\bf G}^{T}{\bf G}\neq {\bf G}{\bf G}^{T}$ and ${\bf G}{\bf G}^{T}=\left[ \begin{array}{cc} 
	{\bf I} & -{\bf L} \\
	-{\bf L}   & {\bf I}+{\bf L}^2 \\
	\end{array} \right]$}:
\begin{equation*}
{\bf G}^{T}{\bf G}=
\left[ \begin{array}{cc} 
{\bf I} & {\bf L} \\
{\bf L}   & {\bf I}+{\bf L}^2 \\
\end{array} \right]
\end{equation*}
The equation
\begin{equation*}
\begin{split}
&\left[ \begin{array}{cc} 
{\bf I} & {\bf L} \\
{\bf L}   & {\bf I}+{\bf L}^2 \\
\end{array} \right]
\left(\begin{array}{c}
{\bf x}_i \\
{\bf v}_i \\
\end{array} \right)=
{\lambda_i}
\left(\begin{array}{c}
{\bf x}_i \\
{\bf v}_i \\
\end{array} \right)
\Longrightarrow \\
&\left\{ 
\begin{array}{l}
{\bf x}_i+{\bf L}{\bf v}_i=\lambda_i {\bf x}_i \\
{\bf L}{\bf x}_i+({\bf I}+{\bf L}^{2}){\bf v}_i=\lambda_i{\bf v}_i\\
\end{array}
\right.
\end{split}
\end{equation*}
is equivalent to
\begin{equation*}
{\bf L}^{2}{\bf v}_i=\frac{\left({\lambda_i}-1\right)^{2}}{\lambda_i}{\bf v}_i
\end{equation*}
Then ${\bf v}_i$ is an eigenvector of the squared Laplacian ${\bf L}^2$ and the values $\lambda_i$ are roots of the equation ${\lambda_i}^{2}-(2+{\mu_i}^2){\lambda_i}+1=0$:
\begin{equation*}
\lambda^{\pm }_i=\frac{1}{2}\left[2+{\mu_i}^{2}\pm\mu_i \sqrt{{\mu_i}^{2}+4} \right]
\end{equation*}

The eigenvalues of $\bf P$ are then given by $\lambda^{{\rm P}\pm }_i=\sqrt{\lambda^{\pm }_i}$, equal to expression (\ref{eigenvalues_P}). As for the eigenvectors, let us consider $\psi^{\rm P +}_i$ in (\ref{eigenvectors_P}), the proof being equivalent for $\psi^{\rm P -}_i$. We have
\begin{equation*}
\begin{split}
&
{\bf P}\psi^{\rm P +}_i=\lambda^{\rm P + }_i\psi^{\rm P +}_i
\Rightarrow
{\bf P}^{2}\psi^{\rm P +}_i={\lambda^{\rm P + }_i}^{2}\psi^{\rm P +}_i\\
&
\left[ \begin{array}{cc} 
{\bf I} & {\bf L} \\
{\bf L}   & {\bf I}+{\bf L}^{2} \\
\end{array} \right]
\left(\begin{array}{c}
{\phi_i} \\
\lambda^{\rm P +}_i \phi_i\\
\end{array} \right)=
{\lambda^{\rm P + }_i}^{2}
\left(\begin{array}{c}
\phi_i \\
\lambda^{\rm P +}_i \phi_i\\
\end{array} \right)\Rightarrow\\
&
\left\{ 
\begin{array}{l}
{\phi_i}+\lambda^{\rm P + }_i{\mu_i}{\phi_i}={\lambda^{\rm P + }_i}^{2}{\phi_i}\\
{\mu_i}{\phi_i}+\lambda^{\rm P + }_i{\phi_i}+\lambda^{\rm P + }_i{\mu_i}^{2}{\phi_i}={\lambda^{\rm P + }_i}^{3}{\phi_i}\\
\end{array}
\right.\Rightarrow\\
&
\left\{ 
\begin{array}{l}
1+\lambda^{\rm P + }_i{\mu_i}={\lambda^{\rm P + }_i}^{2}\\
{\mu_i}+\lambda^{\rm P + }_i+{\mu_i}^{2}\lambda^{\rm P + }_i={\lambda^{\rm P + }_i}^{3}\\
\end{array}
\right.\\
\end{split}
\end{equation*}

Now, the right-hand and left-hand sides of the last two equations are equal, since ${\lambda^{P +}_i}^{2}=\frac{1}{2}\left[2+{\mu_i}^{2}+\mu_i \sqrt{{\mu_i}^{2}+4} \right]$ and ${\lambda^{P +}_i}^{3}=\frac{1}{2}\left[3{\mu_i}+{\mu_i}^{3}+(1+{\mu_i}^{2}) \sqrt{{\mu_i}^{2}+4} \right]$ and this ends the proof.
\end{proof}

\subsection{Proof of Proposition \ref{proposition6}}

\begin{proof}
As a matter of simplicity and without loss of generality, we conduct the proof in the case $n=2$. According to its definition, $\bf U$ is given by
\begin{equation*}
{\bf U}={\bf G}({\bf G}^{T}{\bf G})^{-1/2}={\bf G}{\bf P}^{-1}={\bf G}{\bm \Psi}_{\rm P}{\bm \Lambda}_{\rm P}^{-1}{\bm \Psi}_{\rm P}^{T}
\end{equation*}

where ${\bm \Psi}_{\rm P}$ is the matrix whose columns are eigenvectors of ${\bf P}$ and ${\bm \Lambda}_{\rm P}$ is the diagonal matrix of the eigenvalues of ${\bf P}$:
\begin{equation*}
{\bm \Psi}_{\rm P}=
\left[ \begin{array}{cccc} 
\frac{\phi_1}{\sqrt{1+{\lambda^{\rm P + }_1}^{2}}} &
\frac{\phi_2}{\sqrt{1+{\lambda^{\rm P + }_2}^{2}}} &
\frac{\phi_2}{\sqrt{1+{\lambda^{\rm P - }_2}^{2}}} &
\frac{\phi_1}{\sqrt{1+{\lambda^{\rm P - }_1}^{2}}} \\
\frac{{\lambda^{\rm P + }_1}\phi_1}{\sqrt{1+{\lambda^{\rm P + }_1}^{2}}} &
\frac{{\lambda^{\rm P + }_2}\phi_2}{\sqrt{1+{\lambda^{\rm P + }_2}^{2}}} &
\frac{-{\lambda^{\rm P - }_2}\phi_2}{\sqrt{1+{\lambda^{\rm P - }_2}^{2}}} &
\frac{-{\lambda^{\rm P - }_1}\phi_1}{\sqrt{1+{\lambda^{\rm P - }_1}^{2}}} \\
\end{array} \right]
\end{equation*}

and ${\bm \Lambda}_{\rm P}^{-1}={\rm diag}\, \left\{\frac{1}{\lambda^{\rm P + }_1},\frac{1}{\lambda^{\rm P + }_2},\frac{1}{\lambda^{\rm P - }_2},\frac{1}{\lambda^{\rm P - }_1} \right\}$.

Then ${\bf P}^{-1}={\bm \Psi}_{\rm P}{\bm \Lambda}_{\rm P}^{-1}{\bm  \Psi}_{\rm P}^{T}$ is given by
\begin{widetext}
\begin{equation*}
{\bf P}^{-1}=
\left[ \begin{array}{cc} 
\frac{{\phi_1}{\phi_1}^{T}}{{\lambda^{\rm P + }_1}(1+{\lambda^{\rm P + }_1}^{2})}+
\frac{{\phi_2}{\phi_2}^{T}}{{\lambda^{\rm P + }_2}(1+{\lambda^{\rm P + }_2}^{2})}+
\frac{{\phi_2}{\phi_2}^{T}}{{\lambda^{\rm P - }_2}(1+{\lambda^{\rm P - }_2}^{2})}+
\frac{{\phi_1}{\phi_1}^{T}}{{\lambda^{\rm P - }_1}(1+{\lambda^{\rm P - }_1}^{2})}
&
\frac{{\phi_1}{\phi_1}^{T}}{1+{\lambda^{\rm P + }_1}^{2}}+
\frac{{\phi_2}{\phi_2}^{T}}{1+{\lambda^{\rm P + }_2}^{2}}-
\frac{{\phi_2}{\phi_2}^{T}}{1+{\lambda^{\rm P - }_2}^{2}}-
\frac{{\phi_1}{\phi_1}^{T}}{1+{\lambda^{\rm P - }_1}^{2}}
\\
\frac{{\phi_1}{\phi_1}^{T}}{1+{\lambda^{\rm P + }_1}^{2}}+
\frac{{\phi_2}{\phi_2}^{T}}{1+{\lambda^{\rm P + }_2}^{2}}-
\frac{{\phi_2}{\phi_2}^{T}}{1+{\lambda^{\rm P - }_2}^{2}}-
\frac{{\phi_1}{\phi_1}^{T}}{1+{\lambda^{\rm P - }_1}^{2}}
&
\frac{{\lambda^{\rm P + }_1}{\phi_1}{\phi_1}^{T}}{1+{\lambda^{\rm P + }_1}^{2}}+
\frac{{\lambda^{\rm P + }_2}{\phi_2}{\phi_2}^{T}}{1+{\lambda^{\rm P + }_2}^{2}}+
\frac{{\lambda^{\rm P - }_2}{\phi_2}{\phi_2}^{T}}{1+{\lambda^{\rm P - }_2}^{2}}+
\frac{{\lambda^{\rm P - }_1}{\phi_1}{\phi_1}^{T}}{1+{\lambda^{\rm P - }_1}^{2}}
\end{array} \right]
\end{equation*}

Finally $\bf U$ can be written in the block form (\ref{U}) with
\begin{equation*}
{\bf \cal A}=
\frac{{\phi_1}{\phi_1}^{T}}{1+{\lambda^{\rm P + }_1}^{2}}+
\frac{{\phi_2}{\phi_2}^{T}}{1+{\lambda^{\rm P + }_2}^{2}}-
\frac{{\phi_2}{\phi_2}^{T}}{1+{\lambda^{\rm P - }_2}^{2}}-
\frac{{\phi_1}{\phi_1}^{T}}{1+{\lambda^{\rm P - }_1}^{2}}=
\frac{1-{\lambda^{\rm P + }_1}^{2}}{1+{\lambda^{\rm P + }_1}^{2}}{\phi_1}{\phi_1}^{T}+
\frac{1-{\lambda^{\rm P + }_2}^{2}}{1+{\lambda^{\rm P + }_2}^{2}}{\phi_2}{\phi_2}^{T}
\end{equation*}

and
\begin{equation*}
{\bf \cal B}=
\frac{{\lambda^{\rm P + }_1}{\phi_1}{\phi_1}^{T}}{1+{\lambda^{\rm P + }_1}^{2}}+
\frac{{\lambda^{\rm P + }_2}{\phi_2}{\phi_2}^{T}}{1+{\lambda^{\rm P + }_2}^{2}}+
\frac{{\lambda^{\rm P - }_2}{\phi_2}{\phi_2}^{T}}{1+{\lambda^{\rm P - }_2}^{2}}+
\frac{{\lambda^{\rm P - }_1}{\phi_1}{\phi_1}^{T}}{1+{\lambda^{\rm P - }_1}^{2}}=
\frac{2{\lambda^{\rm P + }_1}}{1+{\lambda^{\rm P + }_1}^{2}}{\phi_1}{\phi_1}^{T}+
\frac{2{\lambda^{\rm P + }_2}}{1+{\lambda^{\rm P + }_2}^{2}}{\phi_2}{\phi_2}^{T}
\end{equation*}
\end{widetext}

In the computations above, we made use of the identities: $1+{\mu_i}{\lambda^{\rm P +}_i}={\lambda^{\rm P +}_i}^{2}$ and $1-{\mu_i}{\lambda^{\rm P -}_i}={\lambda^{\rm P -}_i}^{2}$.
\end{proof}

\subsection{Proof of Proposition \ref{proposition7}}

\begin{proof}
Let us refer to the sign $+$. We have to check that ${\bf U}{\psi^{\rm U+}_j}=\lambda^{\rm U+}_j{\psi^{\rm U+}_j}$:
\begin{equation*}
{\bf U}{\psi^{\rm U+}_j}=\left[ \begin{array}{cc} 
{\bf \cal A} & {\bf \cal B} \\
-{\bf \cal B} & {\bf \cal A} \\
\end{array} \right]
\left( \begin{array}{c} 
\phi_j \\
i\phi_j \\
\end{array} \right)=
\left( \begin{array}{c} 
{\bf \cal A}\phi_j+i{\bf \cal B}\phi_j \\
-{\bf \cal B}\phi_j + i{\bf \cal A}\phi_j \\
\end{array} \right)
\end{equation*}

Since ${\phi_k}^{T}{\phi_j}=\delta_{kj}$, we have
\begin{equation*}
\begin{split}
&{\bf \cal A}\phi_j+i{\bf \cal B}\phi_j\\ =&\sum_{k=1}^{n}\cos(\theta_k){\phi_k}{\phi_k}^{T}{\phi_j}+
i\sum_{k=1}^{n}\sin(\theta_k){\phi_k}{\phi_k}^{T}{\phi_j}\\
=&\cos\theta_j{\phi_j}+i\sin\theta_j{\phi_j}\\
=&(\cos\theta_j+i\sin\theta_j){\phi_j}=\lambda^{\rm U+}_j \phi_j\\
\end{split}
\end{equation*}
\begin{equation*}
\begin{split}
&-{\bf \cal B}\phi_j+i{\bf \cal A}\phi_j\\ =&-\sum_{k=1}^{n}\sin(\theta_k){\phi_k}{\phi_k}^{T}{\phi_j}+
i\sum_{k=1}^{n}\cos(\theta_k){\phi_k}{\phi_k}^{T}{\phi_j}\\
=&-\sin\theta_j{\phi_j}+i\cos\theta_j{\phi_j}\\
=&(\cos\theta_j +i\sin\theta_j)i\phi_j=\lambda^{\rm U+}_j (i\phi_j)\\
\end{split}
\end{equation*}
\end{proof}

{\bf DATA AVAILABILITY}

The original datasets used in this paper are publicly available. Data generated and analyzed during the current study are available from the corresponding author on request.

{\bf AUTHOR DECLARATIONS}

\textbf{Conflict of Interest}

The author has no conflicts to disclose.

\textbf{Author Contributions}

Paolo Bartesaghi: Conceptualization; Investigation; Analysis; Visualization; Writing - original draft; Writing - review \& editing.

\nocite{*}
\bibliography{References}

\begin{thebibliography}{40}%
\makeatletter
\providecommand \@ifxundefined [1]{%
 \@ifx{#1\undefined}
}%
\providecommand \@ifnum [1]{%
 \ifnum #1\expandafter \@firstoftwo
 \else \expandafter \@secondoftwo
 \fi
}%
\providecommand \@ifx [1]{%
 \ifx #1\expandafter \@firstoftwo
 \else \expandafter \@secondoftwo
 \fi
}%
\providecommand \natexlab [1]{#1}%
\providecommand \enquote  [1]{``#1''}%
\providecommand \bibnamefont  [1]{#1}%
\providecommand \bibfnamefont [1]{#1}%
\providecommand \citenamefont [1]{#1}%
\providecommand \href@noop [0]{\@secondoftwo}%
\providecommand \href [0]{\begingroup \@sanitize@url \@href}%
\providecommand \@href[1]{\@@startlink{#1}\@@href}%
\providecommand \@@href[1]{\endgroup#1\@@endlink}%
\providecommand \@sanitize@url [0]{\catcode `\\12\catcode `\$12\catcode
  `\&12\catcode `\#12\catcode `\^12\catcode `\_12\catcode `\%12\relax}%
\providecommand \@@startlink[1]{}%
\providecommand \@@endlink[0]{}%
\providecommand \url  [0]{\begingroup\@sanitize@url \@url }%
\providecommand \@url [1]{\endgroup\@href {#1}{\urlprefix }}%
\providecommand \urlprefix  [0]{URL }%
\providecommand \Eprint [0]{\href }%
\providecommand \doibase [0]{http://dx.doi.org/}%
\providecommand \selectlanguage [0]{\@gobble}%
\providecommand \bibinfo  [0]{\@secondoftwo}%
\providecommand \bibfield  [0]{\@secondoftwo}%
\providecommand \translation [1]{[#1]}%
\providecommand \BibitemOpen [0]{}%
\providecommand \bibitemStop [0]{}%
\providecommand \bibitemNoStop [0]{.\EOS\space}%
\providecommand \EOS [0]{\spacefactor3000\relax}%
\providecommand \BibitemShut  [1]{\csname bibitem#1\endcsname}%
\let\auto@bib@innerbib\@empty
\bibitem [{\citenamefont {Perrot}\ \emph {et~al.}(2016)\citenamefont {Perrot},
  \citenamefont {Béchet}, \citenamefont {Hanzen}, \citenamefont {Arnaud},
  \citenamefont {Pradel},\ and\ \citenamefont {Cézilly}}]{Perrot2016}%
  \BibitemOpen
  \bibfield  {author} {\bibinfo {author} {\bibfnamefont {C.}~\bibnamefont
  {Perrot}}, \bibinfo {author} {\bibfnamefont {A.}~\bibnamefont {Béchet}},
  \bibinfo {author} {\bibfnamefont {C.}~\bibnamefont {Hanzen}}, \bibinfo
  {author} {\bibfnamefont {A.}~\bibnamefont {Arnaud}}, \bibinfo {author}
  {\bibfnamefont {R.}~\bibnamefont {Pradel}}, \ and\ \bibinfo {author}
  {\bibfnamefont {F.}~\bibnamefont {Cézilly}},\ }\href@noop {} {\bibfield
  {journal} {\bibinfo  {journal} {Scientific Reports}\ }\textbf {\bibinfo
  {volume} {6}} (\bibinfo {year} {2016})}\BibitemShut {NoStop}%
\bibitem [{\citenamefont {Stone}\ \emph {et~al.}(2002)\citenamefont {Stone},
  \citenamefont {Olinky}, \citenamefont {Blasius}, \citenamefont {Huppert},\
  and\ \citenamefont {Cazelles}}]{Stone2002}%
  \BibitemOpen
  \bibfield  {author} {\bibinfo {author} {\bibfnamefont {L.}~\bibnamefont
  {Stone}}, \bibinfo {author} {\bibfnamefont {R.}~\bibnamefont {Olinky}},
  \bibinfo {author} {\bibfnamefont {B.}~\bibnamefont {Blasius}}, \bibinfo
  {author} {\bibfnamefont {A.}~\bibnamefont {Huppert}}, \ and\ \bibinfo
  {author} {\bibfnamefont {B.}~\bibnamefont {Cazelles}},\ }\href@noop {}
  {\bibfield  {journal} {\bibinfo  {journal} {AIP Conference Proceedings}\
  }\textbf {\bibinfo {volume} {622}},\ \bibinfo {pages} {476} (\bibinfo {year}
  {2002})}\BibitemShut {NoStop}%
\bibitem [{\citenamefont {Earn}\ \emph {et~al.}(1998)\citenamefont {Earn},
  \citenamefont {Rohani},\ and\ \citenamefont {Grenfell}}]{Earn1998}%
  \BibitemOpen
  \bibfield  {author} {\bibinfo {author} {\bibfnamefont {D.~J.~D.}\
  \bibnamefont {Earn}}, \bibinfo {author} {\bibfnamefont {P.}~\bibnamefont
  {Rohani}}, \ and\ \bibinfo {author} {\bibfnamefont {B.~T.}\ \bibnamefont
  {Grenfell}},\ }\href@noop {} {\bibfield  {journal} {\bibinfo  {journal}
  {Proceedings of the Royal Society of London. Series B: Biological Sciences}\
  }\textbf {\bibinfo {volume} {265}},\ \bibinfo {pages} {7} (\bibinfo {year}
  {1998})}\BibitemShut {NoStop}%
\bibitem [{\citenamefont {Brumley}\ \emph {et~al.}(2012)\citenamefont
  {Brumley}, \citenamefont {Polin}, \citenamefont {Pedley},\ and\ \citenamefont
  {Goldstein}}]{Brumley2012}%
  \BibitemOpen
  \bibfield  {author} {\bibinfo {author} {\bibfnamefont {D.~R.}\ \bibnamefont
  {Brumley}}, \bibinfo {author} {\bibfnamefont {M.}~\bibnamefont {Polin}},
  \bibinfo {author} {\bibfnamefont {T.~J.}\ \bibnamefont {Pedley}}, \ and\
  \bibinfo {author} {\bibfnamefont {R.~E.}\ \bibnamefont {Goldstein}},\
  }\href@noop {} {\bibfield  {journal} {\bibinfo  {journal} {Phys. Rev. Lett.}\
  }\textbf {\bibinfo {volume} {109}},\ \bibinfo {pages} {268102} (\bibinfo
  {year} {2012})}\BibitemShut {NoStop}%
\bibitem [{\citenamefont {Hammond}\ \emph {et~al.}(2007)\citenamefont
  {Hammond}, \citenamefont {Bergman},\ and\ \citenamefont
  {Brown}}]{Hammond2007}%
  \BibitemOpen
  \bibfield  {author} {\bibinfo {author} {\bibfnamefont {C.}~\bibnamefont
  {Hammond}}, \bibinfo {author} {\bibfnamefont {H.}~\bibnamefont {Bergman}}, \
  and\ \bibinfo {author} {\bibfnamefont {P.}~\bibnamefont {Brown}},\
  }\href@noop {} {\bibfield  {journal} {\bibinfo  {journal} {Trends Neurosci.}\
  }\textbf {\bibinfo {volume} {30}},\ \bibinfo {pages} {357} (\bibinfo {year}
  {2007})}\BibitemShut {NoStop}%
\bibitem [{\citenamefont {Fujisaka}\ and\ \citenamefont
  {Yamada}(1983)}]{Yamada1983}%
  \BibitemOpen
  \bibfield  {author} {\bibinfo {author} {\bibfnamefont {H.}~\bibnamefont
  {Fujisaka}}\ and\ \bibinfo {author} {\bibfnamefont {T.}~\bibnamefont
  {Yamada}},\ }\href@noop {} {\bibfield  {journal} {\bibinfo  {journal}
  {Progress of Theoretical Physics}\ }\textbf {\bibinfo {volume} {69}},\
  \bibinfo {pages} {32} (\bibinfo {year} {1983})}\BibitemShut {NoStop}%
\bibitem [{\citenamefont {Pecora}\ and\ \citenamefont
  {Carroll}(1990)}]{Pecora1990}%
  \BibitemOpen
  \bibfield  {author} {\bibinfo {author} {\bibfnamefont {L.~M.}\ \bibnamefont
  {Pecora}}\ and\ \bibinfo {author} {\bibfnamefont {T.~L.}\ \bibnamefont
  {Carroll}},\ }\href@noop {} {\bibfield  {journal} {\bibinfo  {journal} {Phys.
  Rev. Lett.}\ }\textbf {\bibinfo {volume} {64}},\ \bibinfo {pages} {821}
  (\bibinfo {year} {1990})}\BibitemShut {NoStop}%
\bibitem [{\citenamefont {Pecora}\ \emph {et~al.}(2000)\citenamefont {Pecora},
  \citenamefont {Carroll}, \citenamefont {Johnson}, \citenamefont {Mar},\ and\
  \citenamefont {Fink}}]{Pecora2000}%
  \BibitemOpen
  \bibfield  {author} {\bibinfo {author} {\bibfnamefont {L.}~\bibnamefont
  {Pecora}}, \bibinfo {author} {\bibfnamefont {T.}~\bibnamefont {Carroll}},
  \bibinfo {author} {\bibfnamefont {G.}~\bibnamefont {Johnson}}, \bibinfo
  {author} {\bibfnamefont {D.}~\bibnamefont {Mar}}, \ and\ \bibinfo {author}
  {\bibfnamefont {K.~S.}\ \bibnamefont {Fink}},\ }\href@noop {} {\bibfield
  {journal} {\bibinfo  {journal} {International Journal of Bifurcation and
  Chaos}\ }\textbf {\bibinfo {volume} {10}},\ \bibinfo {pages} {273} (\bibinfo
  {year} {2000})}\BibitemShut {NoStop}%
\bibitem [{\citenamefont {Barahona}\ and\ \citenamefont
  {Pecora}(2002)}]{Barahona2002}%
  \BibitemOpen
  \bibfield  {author} {\bibinfo {author} {\bibfnamefont {M.}~\bibnamefont
  {Barahona}}\ and\ \bibinfo {author} {\bibfnamefont {L.}~\bibnamefont
  {Pecora}},\ }\href@noop {} {\bibfield  {journal} {\bibinfo  {journal}
  {Physical Review Letters}\ }\textbf {\bibinfo {volume} {89}},\ \bibinfo
  {pages} {054101} (\bibinfo {year} {2002})}\BibitemShut {NoStop}%
\bibitem [{\citenamefont {Chen}\ \emph {et~al.}(2012)\citenamefont {Chen},
  \citenamefont {Lu}, \citenamefont {Zhan},\ and\ \citenamefont
  {Chen}}]{Chen2012}%
  \BibitemOpen
  \bibfield  {author} {\bibinfo {author} {\bibfnamefont {J.}~\bibnamefont
  {Chen}}, \bibinfo {author} {\bibfnamefont {J.-a.}\ \bibnamefont {Lu}},
  \bibinfo {author} {\bibfnamefont {C.}~\bibnamefont {Zhan}}, \ and\ \bibinfo
  {author} {\bibfnamefont {G.}~\bibnamefont {Chen}},\ }\enquote {\bibinfo
  {title} {Laplacian spectra and synchronization processes on complex
  networks},}\ in\ \href@noop {} {\emph {\bibinfo {booktitle} {Handbook of
  Optimization in Complex Networks: Theory and Applications}}},\ \bibinfo
  {editor} {edited by\ \bibinfo {editor} {\bibfnamefont {M.~T.}\ \bibnamefont
  {Thai}}\ and\ \bibinfo {editor} {\bibfnamefont {P.~M.}\ \bibnamefont
  {Pardalos}}}\ (\bibinfo  {publisher} {Springer US},\ \bibinfo {address}
  {Boston, MA},\ \bibinfo {year} {2012})\ pp.\ \bibinfo {pages}
  {81--113}\BibitemShut {NoStop}%
\bibitem [{\citenamefont {Ren}(2008)}]{Ren2008}%
  \BibitemOpen
  \bibfield  {author} {\bibinfo {author} {\bibfnamefont {W.}~\bibnamefont
  {Ren}},\ }\href@noop {} {\bibfield  {journal} {\bibinfo  {journal}
  {Automatica}\ }\textbf {\bibinfo {volume} {44}},\ \bibinfo {pages} {3195 }
  (\bibinfo {year} {2008})}\BibitemShut {NoStop}%
\bibitem [{\citenamefont {Zhan}\ \emph {et~al.}(2013)\citenamefont {Zhan},
  \citenamefont {Liu},\ and\ \citenamefont {He}}]{Zhan2013}%
  \BibitemOpen
  \bibfield  {author} {\bibinfo {author} {\bibfnamefont {M.}~\bibnamefont
  {Zhan}}, \bibinfo {author} {\bibfnamefont {S.}~\bibnamefont {Liu}}, \ and\
  \bibinfo {author} {\bibfnamefont {Z.}~\bibnamefont {He}},\ }\href@noop {}
  {\bibfield  {journal} {\bibinfo  {journal} {PLOS ONE}\ }\textbf {\bibinfo
  {volume} {8}},\ \bibinfo {pages} {1} (\bibinfo {year} {2013})}\BibitemShut
  {NoStop}%
\bibitem [{\citenamefont {Shang}(2012)}]{Shang2012}%
  \BibitemOpen
  \bibfield  {author} {\bibinfo {author} {\bibfnamefont {Y.}~\bibnamefont
  {Shang}},\ }\href@noop {} {\bibfield  {journal} {\bibinfo  {journal} {Annals
  of the Academy of Romanian Scientists: Series on Mathematics and its
  Applications}\ }\textbf {\bibinfo {volume} {4}} (\bibinfo {year}
  {2012})}\BibitemShut {NoStop}%
\bibitem [{\citenamefont {Pietras}\ and\ \citenamefont
  {Daffertshofer}(2019)}]{Pietras2019}%
  \BibitemOpen
  \bibfield  {author} {\bibinfo {author} {\bibfnamefont {B.}~\bibnamefont
  {Pietras}}\ and\ \bibinfo {author} {\bibfnamefont {A.}~\bibnamefont
  {Daffertshofer}},\ }\href@noop {} {\bibfield  {journal} {\bibinfo  {journal}
  {Physics Reports}\ }\textbf {\bibinfo {volume} {819}},\ \bibinfo {pages} {1}
  (\bibinfo {year} {2019})}\BibitemShut {NoStop}%
\bibitem [{\citenamefont {Arenas}\ \emph {et~al.}(2008)\citenamefont {Arenas},
  \citenamefont {Díaz-Guilera}, \citenamefont {Kurths}, \citenamefont
  {Moreno},\ and\ \citenamefont {Zhou}}]{Arenas2008}%
  \BibitemOpen
  \bibfield  {author} {\bibinfo {author} {\bibfnamefont {A.}~\bibnamefont
  {Arenas}}, \bibinfo {author} {\bibfnamefont {A.}~\bibnamefont
  {Díaz-Guilera}}, \bibinfo {author} {\bibfnamefont {J.}~\bibnamefont
  {Kurths}}, \bibinfo {author} {\bibfnamefont {Y.}~\bibnamefont {Moreno}}, \
  and\ \bibinfo {author} {\bibfnamefont {C.}~\bibnamefont {Zhou}},\ }\href@noop
  {} {\bibfield  {journal} {\bibinfo  {journal} {Physics Reports}\ }\textbf
  {\bibinfo {volume} {469}},\ \bibinfo {pages} {93} (\bibinfo {year}
  {2008})}\BibitemShut {NoStop}%
\bibitem [{\citenamefont {Arola-Fernández}\ \emph {et~al.}(2021)\citenamefont
  {Arola-Fernández}, \citenamefont {Skardal},\ and\ \citenamefont
  {Arenas}}]{Arenas2021}%
  \BibitemOpen
  \bibfield  {author} {\bibinfo {author} {\bibfnamefont {L.}~\bibnamefont
  {Arola-Fernández}}, \bibinfo {author} {\bibfnamefont {P.~S.}\ \bibnamefont
  {Skardal}}, \ and\ \bibinfo {author} {\bibfnamefont {A.}~\bibnamefont
  {Arenas}},\ }\href@noop {} {\bibfield  {journal} {\bibinfo  {journal} {Chaos:
  An Interdisciplinary Journal of Nonlinear Science}\ }\textbf {\bibinfo
  {volume} {31}},\ \bibinfo {pages} {061105} (\bibinfo {year}
  {2021})}\BibitemShut {NoStop}%
\bibitem [{\citenamefont {Eroglu}\ \emph {et~al.}(2017)\citenamefont {Eroglu},
  \citenamefont {Lamb},\ and\ \citenamefont {Pereira}}]{Eroglu2017}%
  \BibitemOpen
  \bibfield  {author} {\bibinfo {author} {\bibfnamefont {D.}~\bibnamefont
  {Eroglu}}, \bibinfo {author} {\bibfnamefont {J.~S.~W.}\ \bibnamefont {Lamb}},
  \ and\ \bibinfo {author} {\bibfnamefont {T.}~\bibnamefont {Pereira}},\
  }\href@noop {} {\bibfield  {journal} {\bibinfo  {journal} {Contemporary
  Physics}\ }\textbf {\bibinfo {volume} {58}},\ \bibinfo {pages} {207}
  (\bibinfo {year} {2017})}\BibitemShut {NoStop}%
\bibitem [{\citenamefont {Bhatta}\ \emph {et~al.}(2021)\citenamefont {Bhatta},
  \citenamefont {Hayat},\ and\ \citenamefont {Sorrentino}}]{Sorrentino2021}%
  \BibitemOpen
  \bibfield  {author} {\bibinfo {author} {\bibfnamefont {K.}~\bibnamefont
  {Bhatta}}, \bibinfo {author} {\bibfnamefont {M.~M.}\ \bibnamefont {Hayat}}, \
  and\ \bibinfo {author} {\bibfnamefont {F.}~\bibnamefont {Sorrentino}},\
  }\href {\doibase 10.1109/TNSE.2021.3097972} {\bibfield  {journal} {\bibinfo
  {journal} {IEEE Transactions on Network Science and Engineering}\ }\textbf
  {\bibinfo {volume} {8}},\ \bibinfo {pages} {2482} (\bibinfo {year}
  {2021})}\BibitemShut {NoStop}%
\bibitem [{\citenamefont {Sánchez-García}(2020)}]{SanchezGarcia2020}%
  \BibitemOpen
  \bibfield  {author} {\bibinfo {author} {\bibfnamefont {R.~J.}\ \bibnamefont
  {Sánchez-García}},\ }\href {\doibase 10.1038/s42005-020-0345-z} {\bibfield
  {journal} {\bibinfo  {journal} {Communications Physics}\ }\textbf {\bibinfo
  {volume} {3}},\ \bibinfo {pages} {87} (\bibinfo {year} {2020})}\BibitemShut
  {NoStop}%
\bibitem [{\citenamefont {Bhatta}\ \emph {et~al.}(2022)\citenamefont {Bhatta},
  \citenamefont {Nazerian},\ and\ \citenamefont {Sorrentino}}]{Sorrentino2022}%
  \BibitemOpen
  \bibfield  {author} {\bibinfo {author} {\bibfnamefont {K.}~\bibnamefont
  {Bhatta}}, \bibinfo {author} {\bibfnamefont {A.}~\bibnamefont {Nazerian}}, \
  and\ \bibinfo {author} {\bibfnamefont {F.}~\bibnamefont {Sorrentino}},\
  }\href {\doibase 10.1109/ACCESS.2022.3188392} {\bibfield  {journal} {\bibinfo
   {journal} {IEEE Access}\ }\textbf {\bibinfo {volume} {10}},\ \bibinfo
  {pages} {72658} (\bibinfo {year} {2022})}\BibitemShut {NoStop}%
\bibitem [{\citenamefont {Grabow}\ \emph {et~al.}(2011)\citenamefont {Grabow},
  \citenamefont {Grosskinsky},\ and\ \citenamefont {Timme}}]{Grabow2011}%
  \BibitemOpen
  \bibfield  {author} {\bibinfo {author} {\bibfnamefont {C.}~\bibnamefont
  {Grabow}}, \bibinfo {author} {\bibfnamefont {S.}~\bibnamefont {Grosskinsky}},
  \ and\ \bibinfo {author} {\bibfnamefont {M.}~\bibnamefont {Timme}},\
  }\href@noop {} {\bibfield  {journal} {\bibinfo  {journal} {The European
  Physical Journal B - Condensed Matter and Complex Systems}\ }\textbf
  {\bibinfo {volume} {84}} (\bibinfo {year} {2011})}\BibitemShut {NoStop}%
\bibitem [{\citenamefont {Skardal}\ \emph {et~al.}(2014)\citenamefont
  {Skardal}, \citenamefont {Taylor},\ and\ \citenamefont {Sun}}]{Skardal2014}%
  \BibitemOpen
  \bibfield  {author} {\bibinfo {author} {\bibfnamefont {P.}~\bibnamefont
  {Skardal}}, \bibinfo {author} {\bibfnamefont {D.}~\bibnamefont {Taylor}}, \
  and\ \bibinfo {author} {\bibfnamefont {J.}~\bibnamefont {Sun}},\ }\href@noop
  {} {\bibfield  {journal} {\bibinfo  {journal} {Phys Rev Lett}\ }\textbf
  {\bibinfo {volume} {113}},\ \bibinfo {pages} {144101} (\bibinfo {year}
  {2014})}\BibitemShut {NoStop}%
\bibitem [{\citenamefont {Skardal}\ \emph {et~al.}(2019)\citenamefont
  {Skardal}, \citenamefont {Taylor},\ and\ \citenamefont {Sun}}]{Skardal2019}%
  \BibitemOpen
  \bibfield  {author} {\bibinfo {author} {\bibfnamefont {P.}~\bibnamefont
  {Skardal}}, \bibinfo {author} {\bibfnamefont {D.}~\bibnamefont {Taylor}}, \
  and\ \bibinfo {author} {\bibfnamefont {J.}~\bibnamefont {Sun}},\ }\href
  {\doibase 10.1137/19M1253836} {\bibfield  {journal} {\bibinfo  {journal}
  {SIAM Journal on Applied Mathematics}\ }\textbf {\bibinfo {volume} {79}},\
  \bibinfo {pages} {2409} (\bibinfo {year} {2019})}\BibitemShut {NoStop}%
\bibitem [{\citenamefont {Su}\ \emph {et~al.}(2009)\citenamefont {Su},
  \citenamefont {Wang},\ and\ \citenamefont {Lin}}]{Su2009}%
  \BibitemOpen
  \bibfield  {author} {\bibinfo {author} {\bibfnamefont {H.}~\bibnamefont
  {Su}}, \bibinfo {author} {\bibfnamefont {X.}~\bibnamefont {Wang}}, \ and\
  \bibinfo {author} {\bibfnamefont {Z.}~\bibnamefont {Lin}},\ }\href@noop {}
  {\bibfield  {journal} {\bibinfo  {journal} {Automatica}\ }\textbf {\bibinfo
  {volume} {45}},\ \bibinfo {pages} {2286} (\bibinfo {year}
  {2009})}\BibitemShut {NoStop}%
\bibitem [{\citenamefont {Dörfler}\ \emph {et~al.}(2013)\citenamefont
  {Dörfler}, \citenamefont {Chertkov},\ and\ \citenamefont
  {Bullo}}]{Dorfler2013}%
  \BibitemOpen
  \bibfield  {author} {\bibinfo {author} {\bibfnamefont {F.}~\bibnamefont
  {Dörfler}}, \bibinfo {author} {\bibfnamefont {M.}~\bibnamefont {Chertkov}},
  \ and\ \bibinfo {author} {\bibfnamefont {F.}~\bibnamefont {Bullo}},\
  }\href@noop {} {\bibfield  {journal} {\bibinfo  {journal} {Proceedings of the
  National Academy of Sciences of the United States of America}\ }\textbf
  {\bibinfo {volume} {110}},\ \bibinfo {pages} {2005} (\bibinfo {year}
  {2013})}\BibitemShut {NoStop}%
\bibitem [{\citenamefont {Liu}\ and\ \citenamefont
  {Barab\'asi}(2016)}]{Liu2016}%
  \BibitemOpen
  \bibfield  {author} {\bibinfo {author} {\bibfnamefont {Y.-Y.}\ \bibnamefont
  {Liu}}\ and\ \bibinfo {author} {\bibfnamefont {A.-L.}\ \bibnamefont
  {Barab\'asi}},\ }\href@noop {} {\bibfield  {journal} {\bibinfo  {journal}
  {Rev. Mod. Phys.}\ }\textbf {\bibinfo {volume} {88}},\ \bibinfo {pages}
  {035006} (\bibinfo {year} {2016})}\BibitemShut {NoStop}%
\bibitem [{\citenamefont {Lu}\ \emph {et~al.}(2014)\citenamefont {Lu},
  \citenamefont {Yu}, \citenamefont {Lü},\ and\ \citenamefont {Xue}}]{Lu2014}%
  \BibitemOpen
  \bibfield  {author} {\bibinfo {author} {\bibfnamefont {R.}~\bibnamefont
  {Lu}}, \bibinfo {author} {\bibfnamefont {W.}~\bibnamefont {Yu}}, \bibinfo
  {author} {\bibfnamefont {J.}~\bibnamefont {Lü}}, \ and\ \bibinfo {author}
  {\bibfnamefont {A.}~\bibnamefont {Xue}},\ }\href@noop {} {\bibfield
  {journal} {\bibinfo  {journal} {IEEE Transactions on Neural Networks and
  Learning Systems}\ }\textbf {\bibinfo {volume} {25}},\ \bibinfo {pages}
  {2110} (\bibinfo {year} {2014})}\BibitemShut {NoStop}%
\bibitem [{\citenamefont {Slotine}\ \emph {et~al.}(2004)\citenamefont
  {Slotine}, \citenamefont {Wang},\ and\ \citenamefont {Rifai}}]{Slotine2004}%
  \BibitemOpen
  \bibfield  {author} {\bibinfo {author} {\bibfnamefont {J.-J.}\ \bibnamefont
  {Slotine}}, \bibinfo {author} {\bibfnamefont {W.}~\bibnamefont {Wang}}, \
  and\ \bibinfo {author} {\bibfnamefont {K.}~\bibnamefont {Rifai}},\
  }\href@noop {} {\bibfield  {journal} {\bibinfo  {journal} {Proceedings of the
  16th International Symposium on Mathematical Theory of Networks and Systems}\
  } (\bibinfo {year} {2004})}\BibitemShut {NoStop}%
\bibitem [{\citenamefont {Dörfler}\ and\ \citenamefont
  {Bullo}(2014)}]{Dorfler2014}%
  \BibitemOpen
  \bibfield  {author} {\bibinfo {author} {\bibfnamefont {F.}~\bibnamefont
  {Dörfler}}\ and\ \bibinfo {author} {\bibfnamefont {F.}~\bibnamefont
  {Bullo}},\ }\href@noop {} {\bibfield  {journal} {\bibinfo  {journal}
  {Automatica}\ }\textbf {\bibinfo {volume} {50}},\ \bibinfo {pages} {1539}
  (\bibinfo {year} {2014})}\BibitemShut {NoStop}%
\bibitem [{\citenamefont {Shahal}\ \emph {et~al.}(2020)\citenamefont {Shahal},
  \citenamefont {Wurzberg}, \citenamefont {Sibony}, \citenamefont {Duadi},
  \citenamefont {Shniderman}, \citenamefont {Weymouth}, \citenamefont
  {Davidson},\ and\ \citenamefont {Fridman}}]{Shahal2020}%
  \BibitemOpen
  \bibfield  {author} {\bibinfo {author} {\bibfnamefont {S.}~\bibnamefont
  {Shahal}}, \bibinfo {author} {\bibfnamefont {A.}~\bibnamefont {Wurzberg}},
  \bibinfo {author} {\bibfnamefont {I.}~\bibnamefont {Sibony}}, \bibinfo
  {author} {\bibfnamefont {H.}~\bibnamefont {Duadi}}, \bibinfo {author}
  {\bibfnamefont {E.}~\bibnamefont {Shniderman}}, \bibinfo {author}
  {\bibfnamefont {D.}~\bibnamefont {Weymouth}}, \bibinfo {author}
  {\bibfnamefont {N.}~\bibnamefont {Davidson}}, \ and\ \bibinfo {author}
  {\bibfnamefont {M.}~\bibnamefont {Fridman}},\ }\href@noop {} {\bibfield
  {journal} {\bibinfo  {journal} {Nature Communications}\ }\textbf {\bibinfo
  {volume} {11}},\ \bibinfo {pages} {3854} (\bibinfo {year}
  {2020})}\BibitemShut {NoStop}%
\bibitem [{\citenamefont {Dumas}\ \emph {et~al.}(2010)\citenamefont {Dumas},
  \citenamefont {Nadel}, \citenamefont {Soussignan}, \citenamefont
  {Martinerie},\ and\ \citenamefont {Garnero}}]{Dumas2010}%
  \BibitemOpen
  \bibfield  {author} {\bibinfo {author} {\bibfnamefont {G.}~\bibnamefont
  {Dumas}}, \bibinfo {author} {\bibfnamefont {J.}~\bibnamefont {Nadel}},
  \bibinfo {author} {\bibfnamefont {R.}~\bibnamefont {Soussignan}}, \bibinfo
  {author} {\bibfnamefont {J.}~\bibnamefont {Martinerie}}, \ and\ \bibinfo
  {author} {\bibfnamefont {L.}~\bibnamefont {Garnero}},\ }\href@noop {}
  {\bibfield  {journal} {\bibinfo  {journal} {PLOS ONE}\ }\textbf {\bibinfo
  {volume} {5}},\ \bibinfo {pages} {1} (\bibinfo {year} {2010})}\BibitemShut
  {NoStop}%
\bibitem [{\citenamefont {Estrada}\ \emph {et~al.}(2017)\citenamefont
  {Estrada}, \citenamefont {Gambuzza},\ and\ \citenamefont
  {Frasca}}]{Estrada2017}%
  \BibitemOpen
  \bibfield  {author} {\bibinfo {author} {\bibfnamefont {E.}~\bibnamefont
  {Estrada}}, \bibinfo {author} {\bibfnamefont {L.}~\bibnamefont {Gambuzza}}, \
  and\ \bibinfo {author} {\bibfnamefont {M.}~\bibnamefont {Frasca}},\ }\href
  {\doibase 10.1137/17M1124310} {\bibfield  {journal} {\bibinfo  {journal}
  {SIAM Journal on Applied Dynamical Systems}\ }\textbf {\bibinfo {volume}
  {17}} (\bibinfo {year} {2017}),\ 10.1137/17M1124310}\BibitemShut {NoStop}%
\bibitem [{\citenamefont {Jenkins}(2013)}]{Jenkins2013}%
  \BibitemOpen
  \bibfield  {author} {\bibinfo {author} {\bibfnamefont {A.}~\bibnamefont
  {Jenkins}},\ }\href {\doibase 10.1016/j.physrep.2012.10.007} {\bibfield
  {journal} {\bibinfo  {journal} {Physics Reports}\ }\textbf {\bibinfo {volume}
  {525}},\ \bibinfo {pages} {167} (\bibinfo {year} {2013})}\BibitemShut
  {NoStop}%
\bibitem [{\citenamefont {Tönjes}\ \emph {et~al.}(2020)\citenamefont
  {Tönjes}, \citenamefont {Fiore},\ and\ \citenamefont {Fiore}}]{Fiore2021}%
  \BibitemOpen
  \bibfield  {author} {\bibinfo {author} {\bibfnamefont {R.}~\bibnamefont
  {Tönjes}}, \bibinfo {author} {\bibfnamefont {C.~E.}\ \bibnamefont {Fiore}},
  \ and\ \bibinfo {author} {\bibfnamefont {C.~E.}\ \bibnamefont {Fiore}},\
  }\href@noop {} {\bibfield  {journal} {\bibinfo  {journal} {Nature
  Communications}\ }\textbf {\bibinfo {volume} {12}},\ \bibinfo {pages} {72}
  (\bibinfo {year} {2020})}\BibitemShut {NoStop}%
\bibitem [{Note1()}]{Note1}%
  \BibitemOpen
  \bibinfo {note} {Note that, since the two less negative eigenvalues for
  network B and C are equal, we have to look at the second less negative
  eigenvalues in round brackets.}\BibitemShut {Stop}%
\bibitem [{Note2()}]{Note2}%
  \BibitemOpen
  \bibinfo {note} {Let us notice that, in order ${\protect \bf B}^{-1}$ exists,
  it is necessary that ${\protect \bf L}^2-4{\protect \bf I}$ is non singular,
  that is $\mu _{i}\protect \neq 2, \forall i=1,\protect \dots , n$. However,
  from a computational point of view, it can be proven that this matrix exists
  for $\mu _{i}\to 2$.}\BibitemShut {Stop}%
\bibitem [{\citenamefont {Cartwright}\ and\ \citenamefont
  {Harary}(1956)}]{Harary1956}%
  \BibitemOpen
  \bibfield  {author} {\bibinfo {author} {\bibfnamefont {D.}~\bibnamefont
  {Cartwright}}\ and\ \bibinfo {author} {\bibfnamefont {F.}~\bibnamefont
  {Harary}},\ }\href@noop {} {\bibfield  {journal} {\bibinfo  {journal}
  {Psychological Review}\ }\textbf {\bibinfo {volume} {5}},\ \bibinfo {pages}
  {277–293} (\bibinfo {year} {1956})}\BibitemShut {NoStop}%
\bibitem [{\citenamefont {Altafini}(2012)}]{Altafini2012}%
  \BibitemOpen
  \bibfield  {author} {\bibinfo {author} {\bibfnamefont {C.}~\bibnamefont
  {Altafini}},\ }\href {\doibase 10.1371/journal.pone.0038135} {\bibfield
  {journal} {\bibinfo  {journal} {PLOS ONE}\ }\textbf {\bibinfo {volume} {7}},\
  \bibinfo {pages} {1} (\bibinfo {year} {2012})}\BibitemShut {NoStop}%
\bibitem [{Note3()}]{Note3}%
  \BibitemOpen
  \bibinfo {note} {Credits. Source: Arab Union of Electricity; Author: Haroun
  Beltaifa; Last Updated: April 13, 2020; Country: Syrian Arab Republic;
  License: Creative Commons Attribution 4.0; url:
  https://energydata.info/dataset/syria-electricity-transmission-network-2017;
  Disclaimer: the author did not modify the original data.}\BibitemShut {Stop}%
\bibitem [{Note4()}]{Note4}%
  \BibitemOpen
  \bibinfo {note} {Note that ${\protect \bf G}^{T}{\protect \bf G}\protect \neq
  {\protect \bf G}{\protect \bf G}^{T}$ and ${\protect \bf G}{\protect \bf
  G}^{T}=\left [ \begin {array}{cc} {\protect \bf I} & -{\protect \bf L} \\
  -{\protect \bf L} & {\protect \bf I}+{\protect \bf L}^2 \\ \end {array}
  \right ]$}\BibitemShut {NoStop}%
\end{thebibliography}%

\end{document}